\documentclass[10pt,twoside]{amsart}
\usepackage{}

\usepackage{amsthm}
\usepackage{amsmath}
\usepackage{amssymb}
\usepackage{amsfonts}
\usepackage{latexsym}
\usepackage{enumitem}
\usepackage{mathrsfs}
\usepackage[all]{xy} \SelectTips{eu}{}
\usepackage{hyperref}
\hypersetup{hidelinks} 
\usepackage{eucal}
\usepackage{xcolor}

\usepackage{tikz,tikz-cd}
\usetikzlibrary{arrows}
\usetikzlibrary{matrix}

\newtheorem{intthm}{Theorem}[]

\newtheorem*{intque*}{Question}
\newtheorem*{intexa*}{Example}
\newcommand{\numberseries}{\bfseries}   
\newlength{\thmtopspace}                
\newlength{\thmbotspace}                
\newlength{\thmheadspace}               
\newlength{\thmindent}                  
\newtheoremstyle{bfupright head,slanted body}
{\thmtopspace}{\thmbotspace}
{\slshape}{\thmindent}{\bfseries}{.}{\thmheadspace}
{{\numberseries \thmnumber{#2\;}}\thmnote{#3}}

\newtheoremstyle{bfupright head,upright body}
{\thmtopspace}{\thmbotspace}
{\upshape}{\thmindent}{\bfseries}{.}{\thmheadspace}
{{\numberseries \thmnumber{#2\;}}\thmnote{#3}}

\newtheoremstyle{fixed bf head,slanted body}
{\thmtopspace}{\thmbotspace}{\slshape}
{\thmindent}{\bfseries}{.}{\thmheadspace}
{{\numberseries \thmnumber{#2\;}}\thmname{#1}\thmnote{ (#3)}}

\newtheoremstyle{fixed bf head,upright body}
{\thmtopspace}{\thmbotspace}{\upshape}
{\thmindent}{\bfseries}{.}{\thmheadspace}
{{\numberseries \thmnumber{#2\;}}\thmname{#1}\thmnote{ (#3)}}

\newtheoremstyle{numbered paragraph}
{\thmtopspace}{\thmbotspace}{\upshape}
{\thmindent}{\upshape}{}{\thmheadspace}
{{\numberseries \thmnumber{#2.}}}

\theoremstyle{bfupright head,slanted body}
\newtheorem{res}{}[section]             \newtheorem*{res*}{}

\theoremstyle{bfupright head,upright body}
\newtheorem{bfhpg}[res]{}               \newtheorem*{bfhpg*}{}

\theoremstyle{fixed bf head,slanted body}
\newtheorem{theorem}[res]{Theorem}          \newtheorem*{theorem*}{Theorem}
\newtheorem{proposition}[res]{Proposition}      \newtheorem*{proposition*}{Proposition}
\newtheorem{corollary}[res]{Corollary}        \newtheorem*{corollary*}{Corollary}
\newtheorem{lemma}[res]{Lemma}            \newtheorem*{lemma*}{Lemma}
\theoremstyle{fixed bf head,upright body}
\newtheorem{definition}[res]{Definition}       \newtheorem*{definition*}{Definition}
\newtheorem{remark}[res]{Remark}           \newtheorem*{remark*}{Remark}
\newtheorem{example}[res]{Example}           \newtheorem*{example*}{Example}
           \newtheorem*{question*}{Question}
           \newtheorem*{fact*}{Fact}
           \newtheorem*{notation*}{Notation}
           \newtheorem*{setup*}{Setup}

\theoremstyle{numbered paragraph}
\newtheorem{ipg}[res]{}

\newlength{\thmlistleft}        
\newlength{\thmlistright}       
\newlength{\thmlistpartopsep}   
\newlength{\thmlisttopsep}      
\newlength{\thmlistparsep}      
\newlength{\thmlistitemsep}     

\newcounter{eqc}
\newenvironment{eqc}{\begin{list}{\upshape (\textit{\roman{eqc}})}%
		{\usecounter{eqc}%
			\setlength{\leftmargin}{\thmlistleft}%
			\setlength{\labelwidth}{\thmlistleft}%
			\setlength{\rightmargin}{\thmlistright}%
			\setlength{\partopsep}{\thmlistpartopsep}%
			\setlength{\topsep}{\thmlisttopsep}%
			\setlength{\parsep}{\thmlistparsep}%
			\setlength{\itemsep}{\thmlistitemsep}}}%
	{\end{list}}%

\newcounter{prt}
\newenvironment{prt}{\begin{list}{\upshape (\alph{prt})}%
		{\usecounter{prt}%
			\setlength{\leftmargin}{\thmlistleft}%
			\setlength{\labelwidth}{\thmlistleft}%
			\setlength{\rightmargin}{\thmlistright}%
			\setlength{\partopsep}{\thmlistpartopsep}%
			\setlength{\topsep}{\thmlisttopsep}%
			\setlength{\parsep}{\thmlistparsep}%
			\setlength{\itemsep}{\thmlistitemsep}}}%
	{\end{list}}%

\newcounter{rqm}
\newenvironment{rqm}{\begin{list}{\upshape (\arabic{rqm})}%
		{\usecounter{rqm}%
			\setlength{\leftmargin}{\thmlistleft}%
			\setlength{\labelwidth}{\thmlistleft}%
			\setlength{\rightmargin}{\thmlistright}%
			\setlength{\partopsep}{\thmlistpartopsep}%
			\setlength{\topsep}{\thmlisttopsep}%
			\setlength{\parsep}{\thmlistparsep}%
			\setlength{\itemsep}{\thmlistitemsep}}}%
	{\end{list}}%

\newenvironment{prf*}[1][Proof]{%
	\begin{proof}[\bf #1]
		\setcounter{equation}{0}
		}
	{\end{proof}
}

\newcommand{\pgref}[1]{\ref{#1}}

\renewcommand{\eqref}[1]{(\pgref{eq:#1})}

\numberwithin{equation}{res}
\def\urltilda{\kern -.15em\lower .7ex\hbox{\~{}}\kern .04em}

\DeclareMathOperator{\Ob}{Ob}

\DeclareMathOperator{\Hom}{Hom}
\DeclareMathOperator{\Mod}{-Mod}
\DeclareMathOperator{\Ind}{{\sf In}}
\DeclareMathOperator{\coInd}{{\sf coIn}}
\DeclareMathOperator{\Res}{{\sf Res}}

\newcommand{\coker}{\mbox{\rm coker}}
\newcommand{\Mon}{\mathsf{Mon}}
\newcommand{\Epi}{\mathsf{Epi}}
\newcommand{\UC}{\mathfrak{C}}
\newcommand{\C}{{\mathscr{C}}}

\newcommand{\op}{{\mathrm{op}}}

\newcommand{\End}{{\mathrm{End}}}
\newcommand{\Ext}{{\mathrm{Ext}}}
\newcommand{\Tor}{{\mathrm{Tor}}}

\newcommand{\xra}[2][]{\xrightarrow[#1]{\:#2\:}}

\newcommand{\calA}{\mathcal{A}}
\newcommand{\calB}{\mathcal{B}}
\newcommand{\calC}{\mathcal{C}}
\newcommand{\calD}{\mathcal{D}}
\newcommand{\calE}{\mathcal{E}}
\newcommand{\calF}{\mathcal{F}}

\newcommand{\calGF}{\mathcal{GF}}
\newcommand{\calPGF}{\mathcal{PGF}}
\newcommand{\calCot}{\mathcal{CT}}
\newcommand{\calGI}{\mathcal{GI}}
\newcommand{\calI}{\mathcal{I}}
\newcommand{\calP}{\mathcal{P}}
\newcommand{\calR}{\mathcal{R}}
\newcommand{\calS}{\mathcal{S}}
\newcommand{\calT}{\mathcal{T}}
\newcommand{\calQ}{\mathcal{Q}}
\newcommand{\calW}{\mathcal{W}}

\newcommand{\scrR}{\mathscr{R}}

\newcommand{\sfC}{\mathsf{C}}

\newcommand{\sfF}{\mathsf{F}}
\newcommand{\sfD}{\mathsf{D}}
\newcommand{\sfW}{\mathsf{W}}

\begin{document}
	
	\title[Representations of generalized linear Reedy categories and abelian model structures]{Representations of generalized linear Reedy categories and abelian model structures}
	
	\author[Z.X. Di]{Zhenxing Di}
	\address{Z.X. Di \ School of Mathematical Sciences, Huaqiao University, Quanzhou 362021, China}
	\email{dizhenxing@163.com}
	
	\author[L.P. Li]{Liping Li}
	\address{L.P. Li \ Department of Mathematics, Hunan Normal University, Changsha 410081, China}
	\email{lipingli@hunnu.edu.cn}
	
	\author[L. Liang]{Li Liang}
	\address{L. Liang \ Department of Mathematics, Gansu Center for Fundamental Research in Complex Systems Analysis and Control, Lanzhou Jiaotong University, Lanzhou 730070, China}
	\email{lliangnju@gmail.com}
	\urladdr{https://sites.google.com/site/lliangnju}
	
	\thanks{Z.X. Di was partly supported by NSF of China (Grant No. 12471034),
		Scientific Research Fund of Fujian Province (Grant No. 605-52525002) and Scientific Research Fund of Huaqiao University (Grant No. 605-50Y22050); L.P. Li was partly supported by NSF of China (Grant No. 12171146); L. Liang was partly supported by NSF of China (Grant No. 12271230) and the Foundation for Innovative Fundamental Research Group Project of Gansu Province (Grant No. 25JRRA805).}
	
	\date{\today}
	
	
	\makeatletter
	\@namedef{subjclassname@2020}{\textup{2020} Mathematics Subject Classification}
	\makeatother
	
\begin{abstract}
We study representations of generalized $k$-linear Reedy categories, extending both the linear Reedy categories of Dalezios and \v{S}t'ov\'{\i}\v{c}ek and $k$-linearizations of generalized Reedy categories of Berger and Moerdijk. We establish homological properties of standard modules and classify irreducible representations in terms of simple modules over the local endomorphism algebras. We give criteria for a decomposition of the representation category into a product of local module categories; under the hypotheses of our non-degeneracy criterion, this decomposition holds if and only if the category is Hom-finite over $k$. Applications include span categories and representations of combinatorial categories. Using Grothendieck bifibrations, we then lift families of complete cotorsion pairs on the local module categories to complete cotorsion pairs on the representation category. Under compatibility assumptions, this also lifts families of (hereditary) Hovey triples and yields abelian model structures, with applications to Gorenstein homological algebra.
\end{abstract}
	
	\maketitle
	\vspace*{.3cm}
	\tableofcontents
	\vspace*{-.9cm}
	\enlargethispage{.4cm}
	
	\thispagestyle{empty}
	\section{Introduction}
	
	\subsection{Motivation}
	The simplex category of finite totally ordered sets and order-preserving maps plays a prominent role in various areas including algebraic topology, homotopy theory, model structures, and higher category theory \cite{GJ}. Via extracting the following fundamental combinatorial properties of this category, Kan introduced the notion of Reedy categories $\C$:
	\begin{prt}
		\item there is a degree function from $\Ob(\C)$ to an ordinal;
		\item it has two subcategories $\C^+$ and $\C^-$ sharing the same objects as $\C$ such that non-identity morphisms in $\C^+$ increase the degree and non-identity morphisms in $\C^-$ decrease the degree;
		\item every morphism in $\C$ can be uniquely expressed as a composite of a morphism in $\C^-$ followed by a morphism in $\C^+$.
	\end{prt}
	Reedy categories are widely used to construct model structures on functor categories. In particular, given a model category $\mathscr{M}$, then the functor category $\mathrm{Fun}(\mathscr{C}, \mathscr{M})$ has a natural model structure called the \textit{Reedy model structure}; for details, see \cite{Hir} or \cite{Ho99}.
	
	Recently, Reedy categories have been generalized from different aspects. To overcome the restriction that all automorphisms in Reedy categories are identities, Berger and Moerdijk introduced generalized Reedy categories in \cite{BM} to encompass non-identity automorphisms, and Shulman defined $g$-Reedy categories and $c$-Reedy categories in \cite{Shu} to include non-identity endomorphisms. These generalizations replace strict uniqueness in (c) by uniqueness up to a suitable equivalence relation, and encompass many new examples, such as the category of finite sets and maps, Segal's category $\Gamma$ \cite{Seg}, Connes' cyclic category $\Lambda$ \cite{Connes}, and categories of spans of orbit categories of finite groups \cite{Dress}.
	
	To study enriched Reedy model structures on enriched functor categories, Angeltveit introduced enriched Reedy categories in \cite{Ang}. More recently, to consider $k$-linear Reedy abelian model structures on $k$-linear functor categories where $k$ is a field, Dalezios and \v{S}t\'ov\'{\i}\v{c}ek defined $k$-linear Reedy categories in \cite{GS} via replacing the unique factorization property (c) by the unique factorization via tensors over $k$, and established many interesting representation theoretic properties of these special $k$-linear categories. In particular, they showed that $k$-linear Reedy categories can be viewed as infinite dimensional analogue of quasi-hereditary algebras, and built a bijective correspondence between the set of isomorphism classes of irreducible representations and the set of objects.
	
	Note that $k$-linearizations of classical Reedy categories are $k$-linear Reedy categories, so the work of \cite{GS} recovers many results on representations of Reedy categories. However, since $k$-linearizations of generalized Reedy categories, $g$-Reedy categories, and $c$-Reedy categories are not $k$-linear Reedy categories in the sense of \cite{GS}, one may wonder to introduce a notion of generalized $k$-linear Reedy categories to unify $k$-linear Reedy categories and $k$-linearizations of the above mentioned generalized Reedy categories. Furthermore, we want to extend the representation theoretic results and abelian module structures in \cite{GS} to this generalized framework. Besides, since many combinatorial categories such as
	\begin{rqm}
		\item[$\bullet$] the category of finite sets and partial injections considered in \cite{CEF, Kuhn, LS},
		\item[$\bullet$] the category of finite dimensional vector spaces over a finite field and linear maps studied in \cite{Kuhn, LS},
		\item[$\bullet$] the category of finite sets and all maps considered in \cite{Pow, PV, Gor},
		\item[$\bullet$] Segal's category $\Gamma$ studied in \cite{LS, Pirash, Seg, Slo},
		\item[$\bullet$] Connes' cyclic category $\Lambda$ considered in \cite{Connes, DK, LS, Loday, Slo},
		\item[$\bullet$] the category of spans studied in \cite{CKM, Dress, Hel, Slo}
	\end{rqm}
	are generalized Reedy categories in the sense of \cite{BM}, and their representations have been extensively studied in the literature, we also hope to provide a systematic approach to unify the various proofs of some important results.

	\subsection{Definition and representations}
	Throughout this paper let $k$ be a commutative ring. By reasonably modifying the three axioms listed in the previous subsection, we say that a small skeletal $k$-linear category $\UC$ is a generalized $k$-linear Reedy category if the following conditions are fulfilled:
	\begin{prt}
		\item there is a degree function $d:\Ob(\UC)\to\lambda$ with values in an ordinal $\lambda$;
		\item it has two $k$-linear subcategories $\UC^+$ and $\UC^-$ sharing the same objects as $\UC$ such that nonzero morphisms between distinct objects in $\UC^+$ increase the degree and nonzero morphisms between distinct objects in $\UC^-$ decrease the degree;
		\item for each $x \in \Ob(\UC)$, one has $\UC^+(x, x) = \UC^-(x, x)$, denoted by $A_x^0$;
		\item for all $x, y \in \Ob(\UC)$, the composition map gives an $(A_y^0, A_x^0)$-bimodule isomorphism
		\[
		\bigoplus_{z \in \Ob(\UC)} \UC^+(z, y) \otimes_{A_z^0} \UC^-(x, z) \longrightarrow \UC(x, y).
		\]
	\end{prt}
	It is easy to see that this definition indeed unifies the concept of $k$-linear Reedy categories in \cite{GS} and $k$-linearizations of generalized Reedy categories in \cite{BM}. Denote by $\UC^0$ the intersection of $\UC^+$ and $\UC^-$, which viewed as a $k$-algebra is a coproduct of ``local" endomorphism algebras $A^0_x$'s with $x$ ranging over all objects in $\UC$.
	
	Dalezios and \v{S}t\'ov\'{\i}\v{c}ek have shown in \cite{GS} that $k$-linear Reedy categories are infinite categorical analogue of quasi-hereditary algebras. As a natural generalization, under certain mild conditions we prove that generalized $k$-linear Reedy categories can be viewed as infinite categorical analogue of standardly stratified algebras. That is, for each object $x$ in $\UC$, we construct a special $\UC$-module $\Delta_x$ (called the \textit{standard module}) satisfying the following properties:
	
	\begin{intthm} \label{thm A}
		Let $\UC$ be a generalized $k$-linear Reedy category. Suppose that $\UC^+$ is a right projective $\UC^0$-module. Then:
		\begin{prt}
			\item $\Delta_x \cong \UC \otimes_{\UC^-} A_x^0$ as left $\UC$-modules;
			\item $\End_{\UC} (\Delta_x) \cong (A_x^0)^{\op}$ as $k$-algebras;
			\item $\Hom_{\UC} (\Delta_x, \Delta_y) \neq 0$ only if $x = y$ or $d(x) > d(y)$;
			\item $\Ext_{\UC}^n (\Delta_x, \Delta_y) \neq 0$ for $n \geqslant 1$ only if $d(x) > d(y)$.
		\end{prt}
	\end{intthm}
	Moreover, one can show that every representable functor $\UC(x, -)$ has a (transfinite) filtration whose factors are closely related to these standard modules. For details, please refer to Proposition \ref{filtration} and Remark \ref{a general filtration}.
	
	When $\UC$ is a $k$-linear Reedy category in the sense of \cite{GS}, Dalezios and \v{S}t\'ov\'{\i}\v{c}ek established a bijective correspondence between $\Ob(\UC)$ and the set of isomorphism classes of irreducible $\UC$-modules. This is no longer true for generalized $k$-linear Reedy categories. However, with some extra assumption, we obtain the following result, extending their parametrization.
	
	\begin{intthm} \label{thm B}
		Let $\UC$ be a generalized $k$-linear Reedy category. Suppose that $A^0_x$ is a semi-perfect $k$-algebra for each $x \in \Ob(\UC)$, $\UC^+$ is a right projective $\UC^0$-module, and $\UC^-$ is a left projective $\UC^0$-module. There is a  bijective correspondence between the following sets:
		\begin{rqm}
			\item $\{ \text{isomorphism classes of irreducible $\UC$-modules} \}$;
			\item $\bigsqcup_{x \in \Ob(\UC)} \{ \text{isomorphism classes of irreducible $A_x^0$-modules }\}$.
		\end{rqm}
	\end{intthm}
	
	\subsection{Applications to representations of combinatorial categories}
	As we mentioned before, representation theory of many important combinatorial categories with an underlying generalized Reedy structure has been widely studied in the literature. In particular, based on the work of Kov\'{a}cs \cite{Kov}, Kuhn proved in \cite{Kuhn} the following surprising result: when $q$ is invertible in $k$, the category of functors from the category of finite dimensional vector spaces over a finite field $\mathbb{F}_q$ and $\mathbb{F}_q$-linear maps to the category of $k$-modules is equivalent to the Cartesian product of module categories of $k\mathrm{GL}_n (\mathbb{F}_q)$ for all $n \geqslant 0$. The analogous result for finite sets and partial injections may be viewed informally as the $q=1$ case; see \cite{CEF}. Similar results in a more abstract setup have been established by quite a few authors; see for instances \cite{LS, Slo, St}.
	
	One goal is to formulate structural criteria for these product decompositions; our non-degeneracy criterion also makes Hom-finiteness a necessary condition, as shown in Corollary~\ref{homfinite-iff-decomposition}. The first criterion is motivated by the work of Kuhn \cite{Kuhn}, relying on the existence of certain special central idempotents. Explicitly, we have:
	
	\begin{intthm} \label{thm C}
		Let $\UC$ be a generalized $k$-linear Reedy category satisfying the following conditions: for all $x, y \in \Ob(\UC)$,
		\begin{prt}
			\item $\UC^+(y, x)$ is a right free $A_y^0$-module and $\UC^-(x, y)$ is a left free $A_y^0$-module,

            \item $\UC^+(y, x) = 0$ if and only if $\UC^-(x, y) = 0$,

			\item there is a central idempotent $e_x \in A_x = \UC(x, x)$ such that
			\[
			e_x A_x e_x = \bigoplus_{z \neq x} \UC^+(z, x) \otimes_{A_z^0} \UC^-(x, z).
			\]
		\end{prt}
		Then one has an equivalence
		\[
		\UC \Mod \simeq \prod_{x \in \Ob(\UC)} A_x^0 \Mod.
		\]
	\end{intthm}
	
	The strategy to establish this theorem is a formalization of that in \cite{Kuhn}. We can show with the given conditions that the set $\{\Delta_x \mid x \in \Ob(\UC) \}$ of standard modules forms an orthogonal family of projective generators in $\UC \Mod$, and the conclusion then follows from a categorical version of Morita equivalence.
	
	The above result indeed covers \cite[Theorem 1.1]{Kuhn} and \cite[Theorem 1.7]{CEF}. However, in practice it is not an easy task to find central idempotents satisfying the condition (c) in Theorem \ref{thm C}. In the next theorem we give another criterion, which can be applied to many combinatorial categories more conveniently.
	
	Given $x, y \in \Ob(\UC)$ and a nonzero morphism $g \in \UC^+(x, y)$, we get a $k$-linear map
	\[
	\varphi_g: \UC^-(y, x) \longrightarrow A_x \longrightarrow A_x^0, \quad f \longmapsto fg \longmapsto \overline{fg}.
	\]
	We say that $g$ is \textit{non-degenerate} if $\varphi_g\neq 0$. Here $A_x^0$ is the quotient of $A_x$ by the ideal of morphisms factoring through objects of degree smaller than $d(x)$ (see Section~4). With respect to this notion, one has:
	
	\begin{intthm} \label{thm D}
		Let $k$ be a field and $\UC$ a generalized $k$-linear Reedy category. Suppose that the following conditions hold for all $x, y \in \Ob(\UC)$:
		\begin{prt}
			\item $\UC^+(x, y)$ is a finite dimensional vector space over $k$;
			\item $\UC^+(x, y)$ is a projective right $A_x^0$-module and $\UC^- (y, x)$ is a free left $A_x^0$-module;
			\item $\dim_k \UC^+(x, y) = \dim_k \UC^-(y, x)$;
			\item every nonzero morphism $g \in \UC^+(x, y)$ is non-degenerate;
            \item $\UC$ is Hom-finite over $k$; that is, $\dim_k\UC(x, y)<\infty$ for all $x, y \in\Ob(\UC)$.
		\end{prt}
		Then one has an equivalence
		\[
		\UC \Mod \simeq \prod_{x \in \Ob(\UC)} A_x^0 \Mod.
		\]
	\end{intthm}
	
	In fact, under assumptions (a)--(d), Hom-finiteness is also necessary for the displayed product equivalence; see Corollary~\ref{homfinite-iff-decomposition}. One may use this theorem as well as its dual version Theorem \ref{dual main result for decomposition} to investigate representations of categories of spans. Explicitly, let $\C$ be small skeletal category such that every endomorphism in $\C$ is an automorphism. Then the binary relation $\leqslant$ on $\Ob(\C)$ such that $x \leqslant y$ if and only if $\C(x, y) \neq \emptyset$ is a partial order. We say that $\C$ is an \textit{artinian EI category} if this partial order satisfies the descending chain condition. In the case that $\C$ has pullbacks, one can define the \textit{category $\widehat{\C}$ of spans} as follows: it has the same objects as $\C$; a morphism in $\widehat{\C}(x, y)$ is represented by a diagram
	$\xymatrix{
		x & z \ar[l]_-f \ar[r]^-g & y
	}$
	with $f$ and $g$ morphisms in $\C$, and two pairs $(f, g)$ and $(f', g')$ represent the same morphism in $\widehat{\C}$ if there is an automorphism $\sigma: z \to z$ such that the following diagram commutes:
	\[
	\xymatrix{
		& z \ar@{-->}[d]^-{\sigma} \ar[dr]^-g \ar[dl]_-f \\
		x & z \ar[l]^-{f'} \ar[r]_-{g'} & y
	}
	\]
	Composition of morphisms are induced by pullbacks; see Subsection 5.2 for details. We mention that many interesting combinatorial categories are equivalent to categories of spans, for example, the category of finite sets and partial injections, the category of finite dimensional vector spaces and linear partial injections, etc.
	
	We can deduce the following result from Theorem \ref{thm D}.
	
	\begin{intthm} \label{thm E}
		Let $k$ be a field and $\C$ an artinian EI category satisfying the following conditions:
		\begin{prt}
			\item $\C$ has pullbacks;
			\item $\C(x, y)$ is a finite set for all $x, y \in \Ob(\C)$;
			\item the order of each group $G_x = \C(x, x)$ is invertible in $k$;
			\item every morphism in $\C$ is a monomorphism;
            \item for every $y\in\Ob(\C)$, only finitely many isomorphism classes of objects $x$ admit a morphism $x\to y$.
		\end{prt}
		Then one has an equivalence
		\[
		\underline{\widehat{\C}} \Mod \simeq \prod_{x \in \Ob(\C)} kG_x \Mod,
		\]
		where $\underline{\widehat{\C}}$ is the $k$-linearization of $\widehat{\C}$.
	\end{intthm}
	
	\subsection{Abelian model structures}
	
	A significant motivation to consider Reedy categories and its variants is to construct model structures on functor categories over them. Here we concentrate on the representation category $\UC \Mod = \mathrm{Fun}_k(\UC, k\Mod)$ and construct abelian model structures on this category. By definition (see Hovey \cite{Ho02}), these are model structures compatible with the abelian structure of $\UC\Mod$; that is, cofibrations coincide with monomorphisms with cofibrant cokernels and fibrations coincide with epimorphisms with fibrant kernels. By Hovey's correspondence, this is equivalent to constructing \emph{Hovey triples} $(\calQ, \calW, \calR)$ of classes of objects in $\UC\Mod$ such that both $(\calQ, \calW \cap \calR)$ and $(\calQ \cap \calW, \calR)$ are complete cotorsion pairs, and $\calW$ is thick. Compared to abelian model structures defined in terms of morphisms, Hovey triples can be described by properties of objects, which are much more transparent in practice.
	
	
	We need to introduce a few more notions before stating our main result in the second part of this paper. Let $\UC$ be a generalized $k$-linear Reedy category. Given any family $\calS = \{ \calS_x \}_{x \in \Ob(\UC)}$ of classes of objects in those $A_x^0 \Mod$'s, denote by $\UC\Mod_\calS$ the class of $\UC \Mod$ consisting of all $\UC$-modules $Y$ such that $Y(x)$ belongs to $\calS_x$ for each $x \in \Ob(\UC)$. We also consider two classes $\Phi(\calS)$ and $\Psi(\calS)$ of $\UC$-modules:
	\begin{align*}
		\Phi(\calS)
	&=\left\{ Y \in \UC\Mod \:
	\left|
	\begin{array}{c}
		\text{for each}\ x \in \Ob(\UC) \text{ with the degree}\ \alpha, \\
		l^{\alpha}_Y(x)\ \text{is monic and}\ \coker(l^{\alpha}_Y(x)) \in \calS_x
	\end{array}
	\right.
	\right\},\ \mathrm{and} \\
	\Psi(\calS)
	&=\left\{ Y \in \UC\Mod \:
	\left|
	\begin{array}{c}
		\text{for each}\ x \in \Ob(\UC) \text{ with the degree}\ \alpha, \\
		m^{\alpha}_Y(x)\ \text{is epic and}\ \ker(m^{\alpha}_Y(x)) \in \calS_x
	\end{array}
	\right.
	\right\}.
	\end{align*}
	Here, $l^{\alpha}_Y(x): 1_x \UC \otimes_{\UC_{\alpha}} \Res_{\alpha}Y \to Y(x)$ is the counit of the adjunction
	\[
	\Ind_{\alpha}: \UC_{\alpha}\Mod \rightleftharpoons \UC\Mod :  \Res_{\alpha}
	\]
	evaluated at $Y$ and $x$, and $m^{\alpha}_Y(x): Y(x) \to \Hom_{\UC_{\alpha}} (\UC 1_x, \Res_{\alpha}Y)$ is the unit of the adjunction
	\[
	\Res_{\alpha}: \UC \Mod\rightleftharpoons \UC_{\alpha}\Mod : \coInd_{\alpha}
	\]
	evaluated at $Y$ and $x$, where $\UC_{\alpha}$ is the full subcategory of $\UC$ consisting of objects with degree less than $\alpha$; see Section \ref{lat ma fun} for details.
	
	The next result asserts that under some mild conditions, a family of (hereditary) Hovey triples in those $A_x^0 \Mod$'s induces a (hereditary) Hovey triple in $\UC \Mod$.
	
	\begin{intthm} \label{thm F}
		Let $\UC$ be a generalized $k$-linear Reedy category with $\UC^+$ a projective right $\UC^0$-module and $\UC^-$ a projective left $\UC^0$-module. Suppose that $(\calQ, \calW, \calR) = \{(\calQ_x, \calW_x, \calR_x)\}_{x \in \Ob(\UC)}$ is a family with $(\calQ_x, \calW_x, \calR_x)$ a (hereditary) Hovey triple in $A^0_x\Mod$. If $\{ \calQ_x \cap \calW_x \}_{x \in \Ob(\UC)}$ is cocompatible and $\{ \calW_x \cap \calR_x \}_{x \in \Ob(\UC)}$ is  compatible, then
		\[
		( \Phi(\calQ), \, \UC\Mod_\calW, \, \Psi(\calR) )
		\]
		forms a (hereditary) Hovey triple in $\UC \Mod$.
	\end{intthm}
	
	For definitions of compatible and cocompatible conditions, see Definition \ref{compatible}.
	
	\begin{remark*}
		For the case that the target category of representations is the category of $k$-modules, Theorem \ref{thm F} extends \cite[Theorem 7.2]{GS} from $k$-linear Reedy categories with $k$ a field to generalized $k$-linear Reedy categories. As in \cite{GS}, the key construction relies on a transfinite induction to glue a family of (complete) cotorsion pairs in those $A_x^0 \Mod$'s to obtain a (complete) cotorsion pair in $\UC \Mod$. However, for our convenience, we use a new approach based on Grothendieck bifibrations. In details, the strategy to carry out the induction step can be divided into the three steps:
		\begin{rqm}
			\item We give a transparent description of the fiber ${\UC_{\alpha +1}\Mod}_V$ of a $\UC_{\alpha}$-module $V$, and show that the restriction functor $\Res_{\alpha} : \UC_{\alpha +1}\Mod \to \UC_{\alpha}\Mod$ induced by the embedding $\iota_{\alpha} : \UC_{\alpha} \to \UC_{\alpha + 1}$ is a Grothendieck bifibration; see Theorem \ref{res is Gro bifib}.
			
			\item Given a Grothendieck bifibration $p : \calT \to \calB$, Stanculescu \cite{Stan12} provided a method to construct a weak factorization system in the total category $\calT$ via combining weak factorization systems in the fibers $\calT_A$ for $A \in \Ob(\calB)$ and a weak factorization system in the basis category $\calB$; see also Cagne and Melli\`es \cite[Lemma 2.4]{CM20}. Applying this result, we obtain a weak factorization system in $\UC_{\alpha +1}\Mod$ via gluing weak factorization systems in the fibers ${\UC_{\alpha +1}\Mod}_V$ for all $\UC_{\alpha}$-modules $V$ and a weak factorization system in $\UC_{\alpha}\Mod$.
			
			\item Finally, by the bijective correspondence between complete cotorsion pairs and weak factorization systems established in \cite{Ho02} (see also Positselski and \v{S}\v{t}ov\'{\i}\v{c}ek \cite[Theorem 2.4]{Po-St22}), we obtain desired complete cotorsion pairs in $\UC_{\alpha+1} \Mod$, completing the induction step.
		\end{rqm}
	\end{remark*}
	
	We then give some examples of Theorem \ref{thm F} in Gorenstein homological algebra, and construct several Hovey triples in $\UC \Mod$ for some special kinds of generalized $k$-linear Reedy categories. Explicitly, if $\UC$ is a \textit{generalized $k$-linear direct category} (that is, $\UC = \UC^+$) such that it is projective both as a left and a right $\UC^0$-module (see Example \ref{direct}), then
	\[
	( \Phi(\calPGF), \, \UC\Mod_{\calPGF^\perp}, \, \UC\Mod )\ \quad \text{and}\ \quad
	( \Phi(\calGF), \, \UC\Mod_{\calPGF^\perp}, \, \UC\Mod_{\calCot} )
	\]
	form hereditary Hovey triples in $\UC \Mod$; see Corollaries \ref{PGF} and \ref{GF}. Dually, if $\UC$ is a \textit{generalized $k$-linear inverse category} (that is, $\UC = \UC^-$) such that it is projective both as a left and a right $\UC^0$-module, then
	\[
	( \UC\Mod, \, \UC\Mod_{^{\perp}{\calGI}}, \, \Psi(\calGI) )
	\]
	forms a hereditary Hovey triple in $\UC \Mod$; see Corollary \ref{GI}.
	
	\subsection{Organization}
	
	The paper is divided into two parts, the first one of which handles representation theory of generalized $k$-linear Reedy categories $\UC$. Specifically, in Section 2 we introduce generalized $k$-linear Reedy categories, and study their elementary properties. Standard modules are introduced in Section 3, where we prove Theorem \ref{thm A} and give a parameterization of isomorphism classes of irreducible representations of $\UC$. Theorem \ref{thm C} is proved in Section 4, while in Section 5 we prove Theorem \ref{thm D} and describe a few applications, including representations of categories of spans as well as semisimplicity of Mackey functors.
	
	In the second part we focus on abelian model structures on $\UC \Mod$. For the convenience of the reader, we describe in Section 6 necessary background knowledge on abelian model structures and Grothendieck (op)fibrations. The main technical tools, namely the restriction functor and its adjoints, are considered in Section 7. In Section 8, we further study the restriction functor, revealing its bifibrational structure. In Section 9, we show that a given family of complete cotorsion pairs in $A_x^0 \Mod$'s can be glued to a complete cotorsion pair in $\UC \Mod$. As a direct consequence of this result, we show in the last section that a given family of abelian model structures on $A_x^0 \Mod$'s can be also glued to an abelian model structure on $\UC\Mod$, and describe a few examples.
	
	
	\part*{Part I. Representations of Reedy categories}\label{part:reedy}
	\noindent
	In this part, we introduce generalized $k$-linear Reedy categories, explore their representations, and pursue applications to representation theory of combinatorial categories with an underlying generalized Reedy structure.
	
	\section{Preliminaries on Reedy categories and representations}
	\noindent
	In this section we describe some preliminary knowledge on generalized $k$-linear Reedy categories.
	
	\subsection{Linear categories and representations}
	
	Throughout the paper let $k$ be a commutative ring with identity, and let $\UC$ be a small skeletal $k$-linear category. Given $x, y \in \Ob(\UC)$, denote by $\UC(x, y)$ the set of all morphisms from $x$ to $y$, which has a natural structure of a $k$-module. In particular, when $x = y$, we denote $\UC(x, x)$ by $A_x$, which is a $k$-algebra. It is clear that $\UC(x, y)$ is an $(A_y, A_x)$-bimodule.
	
	Given two $k$-linear subcategories $\mathfrak{D}$ and $\mathfrak{E}$ of $\UC$, define the intersection $\mathfrak{D} \cap \mathfrak{E}$ to be the $k$-linear category: $\Ob(\mathfrak{D} \cap \mathfrak{E}) = \Ob(\mathfrak{D}) \cap \Ob(\mathfrak{E})$; for $x, y \in \Ob(\mathfrak{D} \cap \mathfrak{E})$, set
	$(\mathfrak{D} \cap \mathfrak{E}) (x, y) = \mathfrak{D} (x, y) \cap \mathfrak{E} (x, y)$.
	
	By definition, a \textit{representation} of $\UC$, or a left \textit{$\UC$-module}, is a $k$-linear covariant functor $V$ from $\UC$ to $k \Mod$, the category of all $k$-modules. Fixing an object $x$, we denote by $V_x$ or $V(x)$ the value of $V$ on $x$. Morphisms between two representations of $\UC$ are $k$-linear natural transformations. We say that $\UC$ is \textit{locally finite} if $\UC(x, y)$ is a finitely generated $k$-module for all $x, y \in \Ob(\UC)$. Similarly, a left $\UC$-module $V$ is called \textit{locally finite} if each $V_x$ is a finitely generated $k$-module.
	
	Sometimes it is more convenient to view $\UC$ as a (non-unital) $k$-algebra $A_{\UC}$. More explicitly, as $(k, k)$-bimodules, one has
	\[
	A_{\UC} = \bigoplus_{x, y \in \Ob(\UC)} \UC(x, y),
	\]
	and given two morphisms $f \in \UC(x, y)$ and $g \in \UC(z, w)$, one defines the product $gf$ to be the composite if $y = z$, and to be 0 otherwise. It is easy to check that $A_{\UC}$ is an associative $k$-algebra. Moreover, if $\mathfrak{D}$ and $\mathfrak{E}$ are $k$-linear subcategories of $\UC$, then $A_{\mathfrak{D}}$ and $A_{\mathfrak{E}}$ are (non-unital) $k$-subalgebras of $A_{\UC}$, and
	\[
	A_{\mathfrak{D} \cap \mathfrak{E}} = A_{\mathfrak{D}} \cap A_{\mathfrak{E}}.
	\]
	
	Every left $\UC$-module can be viewed as a left $A_{\UC}$-module in a natural way, but the converse statement only holds if $\Ob(\UC)$ is a finite set. Indeed, a left $A_{\UC}$-module $M$ is a left $\UC$-module if and only if it has the following decomposition
	\[
	M = \bigoplus_{x \in \Ob(\UC)} 1_x M,
	\]
	where $1_x$ is the identity morphism on $x$, and in this case we call $M$ a \textit{graded} left $A_{\UC}$-module. It is well known that $\UC \Mod$, the category of left $\UC$-modules, is isomorphic to the category of graded left $A_{\UC}$-modules, and hence can be viewed as a full subcategory of $A_{\UC} \Mod$. Moreover, for each object $x$, the projective left $A_{\UC}$-module $A_{\UC} 1_x$ is a graded left $A_{\UC}$-module, which is nothing but the representable functor $\UC(x, -)$, and sometimes we also denote it by $\UC 1_x$ by abuse of notations.
	
	\subsection{Generalized $k$-linear Reedy categories}
	
	Now we introduce generalized $k$-linear Reedy categories. Let $\lambda$ be an ordinal, and $d: \Ob(\UC) \to \lambda$ be a \textit{degree function}.
	
	\begin{definition}
		A small skeletal $k$-linear category $\UC$ is called a generalized $k$-linear Reedy category with respect to the degree function $d$ if it has two subcategories $\UC^+$ and $\UC^-$ such that the following conditions hold:
		\begin{prt}
			\item $\Ob(\UC^+) = \Ob(\UC)$ and $(\UC^+(x, y) \neq 0) \Rightarrow (x = y) \vee (d(y) > d(x))$;
			
			\item $\Ob(\UC^-) = \Ob(\UC)$ and $(\UC^-(x, y) \neq 0) \Rightarrow (x = y) \vee (d(x) > d(y))$;
			
			\item $\UC^+(x, x) = \UC^-(x, x)$ for $x \in \Ob(\UC)$;
			
			\item Let $\UC^0 = \UC^+ \cap \UC^-$. Then the composition induces a $(\UC^0, \UC^0)$-bimodule isomorphism
			$$\rho: A_{\UC^+} \otimes_{A_{\UC^0}} A_{\UC^-} \to A_{\UC}.$$
		\end{prt}
	\end{definition}
	
	\begin{remark} \label{ring isomorphism}
		Although $A_{\UC^+} \otimes_{A_{\UC^0}} A_{\UC^-}$ is a tensor of two rings, in general it cannot be equipped with a ring structure given by the tensor product. However, the above isomorphism $\rho$ allows one to define a ring structure on it. Explicitly, given $f^+ \otimes f^-$ and $g^+ \otimes g^-$ in $A_{\UC^+} \otimes_{A_{\UC^0}} A_{\UC^-}$, their product is defined to be
		\[
		(f^+ \otimes f^-) \ast (g^+ \otimes g^-) = \rho^{-1} (\rho(f^+ \otimes f^-) \rho(g^+ \otimes g^-)) = \rho^{-1} (f^+f^-g^+g^-).
		\]
		Via this construction $\rho$ becomes a (non-unital) ring isomorphism, so we identify $A_{\UC^+} \otimes_{A_{\UC^0}} A_{\UC^-}$ with $A_{\UC}$. With respect to this identification, the subring $A_{\UC^+} \otimes_{A_{\UC^0}} A_{\UC^0}$ (resp, $A_{\UC^0} \otimes_{A_{\UC^0}} A_{\UC^-}$ and $A_{\UC^0} \otimes_{A_{\UC^0}} A_{\UC^0}$) of $A_{\UC^+} \otimes_{A_{\UC^0}} A_{\UC^-}$ is identified with the subring $A_{\UC^+}$ (resp., $A_{\UC^-}$ and $A_{\UC^0}$) of $A_{\UC}$.
		
		In the rest of this paper we always equip $A_{\UC^+} \otimes_{A_{\UC^0}} A_{\UC^-}$ the multiplication $\ast$. To simplify notation, from now on we identify $k$-linear categories $\UC$, $\UC^+$, $\UC^-$ and $\UC^0$ with their associated $k$-algebras. Although as pointed out before, not all $A_{\UC}$-modules are $\UC$-modules, this identification does not cause troubles to us since in this paper we only consider graded $A_{\UC}$-modules, which are precisely $\UC$-modules. Hopefully this simplification does not cause confusions for reader to figure out the precisely meaning (a $k$-linear category or a $k$-algebra) of $\UC$ from the context.
		
		By these conventions, one can restate the condition (d) as $(\UC^+ \otimes_{\UC^0} \UC^-, \, \ast) \cong \UC$ as rings.
	\end{remark}
	
	\begin{remark}
		The above definition can extend to non-skeletal $k$-linear categories if we relax the requirement $x = y$ in (a) and (b) to the weaker condition $x \cong y$. To simplify superficial technical issues, we do not pursue this full generality in this paper.
	\end{remark}
	
	Let us describe a few immediate consequences of the above definition. The degree function $d$ does not induce a partial order $\preccurlyeq$ on $\Ob(\UC)$ via setting $x \preccurlyeq y$ when $d(x) \leqslant d(y)$, since distinct objects can have the same degree, and hence the anti-symmetry fails. Instead, we can use morphisms in $\UC^-$ to obtain a satisfactory partial order on $\Ob(\UC)$. Explicitly, given objects $x$ and $y$, define $y \preccurlyeq x$ if there exists a finite sequence of objects $x = x_0, \, x_1, \, \ldots, \, x_n = y$ such that $\UC^-(x_i, x_{i+1}) \neq 0$.
	
	\begin{lemma} \label{partial order}
		The binary relation $\preccurlyeq$ is an artinian partial order on $\Ob(\UC)$.
	\end{lemma}
	
	\begin{prf*}
		It is clear that $\preccurlyeq$ is reflexive and transitive. Suppose that $y \preccurlyeq x$ and $x \preccurlyeq y$. Then we can find two sequences
		\[
		x = x_0, \, x_1, \, \ldots, \, x_m = y \quad \mathrm{and} \quad y = y_0, \, y_1, \, \ldots, \, y_n = x
		\]
		such that $\UC^-(x_i, x_{i+1}) \neq 0$ and $\UC^-(y_i, y_{i+1}) \neq 0$. Accordingly, one has
		\[
		d(x) = d(x_0) \geqslant d(x_1) \geqslant \ldots \geqslant d(x_m) = d(y)
		\quad \mathrm{and} \quad
		d(y) = d(y_0) \geqslant d(y_1) \geqslant \ldots \geqslant d(y_n) = d(x).
		\]
		Thus, all degrees are equal. But $\UC^-(x_i, x_{i+1}) \neq 0$ implies that $x_i = x_{i+1}$ since they share the same degree. Consequently, one has $x = y$, so $\preccurlyeq$ also has the anti-symmetry property. In conclusion, the binary relation $\preccurlyeq$ is a partial order on $\Ob(\UC)$.
		
		We claim that $d(y) < d(x)$ if $y \prec x$. Indeed, since $y \prec x$, we can find a finite sequence $x = x_0, x_1, \ldots, x_n = y$ such that $\UC^-(x_i, x_{i+1}) \neq 0$ and there is at least one $i$ with $x_i \neq x_{i+1}$. But this forces $d(x_{i+1}) < d(x_i)$, so $d(y) < d(x)$, as claimed. Thus, $d$ is an order-preserving map from the poset $(\Ob(\UC), \preccurlyeq)$ to the ordinal $\lambda$, so $\preccurlyeq$ must be artinian.
	\end{prf*}
	
	\begin{remark} \label{another partial order}
		One can define another partial order $\preccurlyeq'$ on $\Ob(\UC)$ via $\UC^+$. Explicitly, given objects $x$ and $y$, define $x \preccurlyeq' y$ if there exists a finite sequence of objects $x = x_0, \, x_1, \, \ldots, \, x_n = y$ such that $\UC^+(x_i, x_{i+1}) \neq 0$. The reader can check that this is also an artinian partial order, but in general we cannot expect that $\preccurlyeq$ and $\preccurlyeq'$ coincide. For example, let $\UC$ be the $k$-linearization of the category whose objects are $[n] = \{1, 2, \ldots, n \}$ and morphisms are injections. Define $d: \Ob(\UC) \to \lambda$ by setting $d([n]) = n$ where $\lambda$ is an arbitrary infinite ordinal. Then $\UC = \UC^+$. It is easy to see that $\preccurlyeq$ is the discrete partial order, while $\preccurlyeq'$ is isomorphic to the usual linear order on $\mathbb{N}$.
		
		Note that if $\UC^-(x, y) \neq 0$, then $y \preccurlyeq x$. Similarly, if $\UC^+(x, y) \neq 0$, then $x \preccurlyeq' y$. But the converse statements fail since the composite of two nonzero morphisms might be 0. Furthermore, in the proof of Lemma \ref{partial order} we have shown that $d$ is an order-preserving map from $(\Ob(\UC), \preccurlyeq)$ to the ordinal $\lambda$. The same conclusion holds while replacing $(\Ob(\UC), \preccurlyeq)$ by $(\Ob(\UC), \preccurlyeq')$. Therefore, we have the following strict implications:
      	\begin{align*}
		\UC^+(x, y) &\neq 0 \Rightarrow x \preccurlyeq' y \Rightarrow (x = y) \vee (d(x) < d(y)),\ \mathrm{and}\\
		\UC^-(x, y) &\neq 0 \Rightarrow y \preccurlyeq x \Rightarrow (x = y) \vee (d(y) < d(x)).
		\end{align*}
	\end{remark}
	
	\textbf{Throughout this paper}, we always use the partial order $\preccurlyeq$ on $\Ob(\UC)$ unless otherwise specified.
	
	The following results can also be easily deduced from the definition.

	\begin{lemma} \label{Reedy factorization}
		Suppose that $\UC$ is a generalized $k$-linear Reedy category. Then for any $x, y \in \Ob(\UC)$,
		\begin{prt}
			\item $\UC^0(x, y) \neq 0$ if and only if $x = y$;
			
			\item the composition $\rho$ induces an isomorphism of $(A_y^0, A_x^0)$-modules
			\[
			\bigoplus_{z \in \Ob(\UC)} \UC^+(z, y) \otimes_{A_z^0} \UC^-(x, z) \cong \UC(x, y),
			\]
			where $A_z^0 = \UC^0(z, z)$;
			
			\item if $d(x) = d(y)$ is minimal but $x \neq y$, then $\UC(x, y) = 0$.
		\end{prt}
	\end{lemma}
	
	\begin{prf*}
		The if direction of (a) is trivial. Conversely, if $\UC^0(x, y) \neq 0$, then both $\UC^+(x, y)$ and $\UC^-(x, y)$ are nonzero. In the case $x \neq y$, one shall have $d(x) < d(y)$ and $d(x) > d(y)$, which is impossible. Thus, $x = y$.
		
		By (a), we have the following decomposition of $k$-algebras
		\[
		\UC^0 = \bigoplus_{x \in \Ob(\UC)} \UC^0(x, x) = \bigoplus_{x \in \Ob(\UC)} A_x^0.
		\]
		Correspondingly, for $x \in \Ob(\UC)$, the left $\UC^0$-module $\UC^- 1_x$ has the following decomposition
		\[
		\UC^- 1_x = \bigoplus_{z \in \Ob(\UC)} 1_z \UC^- 1_x = \bigoplus_{z \in \Ob(\UC)} \UC^-(x, z)
		\]
		since each $\UC^-(x, z)$ is a left $A_z^0$-module, and hence, is a left $\UC^0$-module.
		
		Now we prove (b). Let $1_x$ and $1_y$ be the identity morphisms on $x$ and $y$ respectively. Then
		\begin{align*}
			\UC(x, y) & = 1_y \UC 1_x\\
			& \cong 1_y \UC^+ \otimes _{\UC^0} \UC^- 1_x\\
			& \cong 1_y \UC^+ \otimes_{\UC^0} (\bigoplus_{z \in \Ob(\UC)} 1_z \UC^- 1_x)\\
			& \cong \bigoplus_{z \in \Ob(\UC)} (1_y \UC^+) \otimes_{\UC^0} (1_z \UC^- 1_x)\\
			& = \bigoplus_{z \in \Ob(\UC)} (1_y \UC^+ 1_z) \otimes_{\UC^0} (1_z \UC^- 1_x)\\
			& = \bigoplus_{z \in \Ob(\UC)} \UC^+ (z, y) \otimes_{\UC^0} \UC^-(x, z)\\
			& \cong \bigoplus_{z \in \Ob(\UC)} \UC^+ (z, y) \otimes_{A_z^0} \UC^-(x, z),
		\end{align*}
		where the first isomorphism holds by the condition (d) in the definition of Reedy categories,\footnote{Here we do not require the conclusion that $\rho$ induces a ring isomorphism given in Remark \ref{ring isomorphism}.} and the last isomorphism follows from the observation that $\UC^0$ is a direct sum of $A_z^0$'s, $\UC^+ (z, y)$ is a right $A_z^0$-module, and $\UC^-(x, z)$ is a left $A_z^0$-module.
		
		Statement (c) immediately follows from the isomorphism in (b). Indeed, if $\UC(x, y) \neq 0$, then by the isomorphism in (b), we shall have a certain object $z$ such that $\UC^- (x, z) \neq 0 \neq \UC^+(z, y)$. But since $d(x) = d(y)$ is minimal, one must have $x = z = y$ by the definition of generalized $k$-linear Reedy categories.
	\end{prf*}
	
	\begin{remark} \label{uniqueness of factorization}
		By the statement (b) in the previous lemma, given a morphism $f: x \to y$ in $\UC$, one obtains the decomposition
		$f = f_1^+f_1^- + \ldots + f_n^+ f_n^-$
		with $f_i^+ \in \UC^+(x_i, y)$ and $f_i^- \in \UC^-(x, x_i)$. Furthermore, this expression is unique in the following sense. Let
		$f = g_1^+g_1^- + \ldots + g_m^+ g_m^-$
		be another decomposition with $g_i^+ \in \UC^+(y_i, y)$ and $g_i^- \in \UC^-(x, y_i)$. Without loss of generality, assume that
		\[
		x_1 = x_2 = \ldots = x_{i_1}, \, x_{i_1 + 1} = x_{i_1 +2} = \ldots  = x_{i_2}, \, \ldots, \, x_{i_{r-1} + 1} = \ldots = x_{i_r}
		\]
		with $i_r = n$ and
		\[
		y_1 = y_2 = \ldots = y_{j_1}, \, y_{j_1 + 1} = y_{j_1 +2} = \ldots  = y_{j_2}, \, \ldots, \, y_{j_{s-1} + 1} = \ldots = y_{j_s}
		\]
		with $j_s = m$ such that $x_{i_t}$'s (resp., $y_{j_t}$'s) are distinct for $1 \leqslant t \leqslant r$ (resp., $1 \leqslant t \leqslant s$). Then $r = s$, $y_{i_t} = x_{i_t}$ for $1 \leqslant t \leqslant s$, and
		\[
		g_{i_{t-1} + 1} ^+ \otimes g_{i_{t-1} + 1}^- + \ldots + g_{i_t}^+ \otimes g_{i_t}^- = f_{i_{t-1} + 1} ^+ \otimes f_{i_{t-1} + 1}^- + \ldots + f_{i_t}^+ \otimes f_{i_t}^-
		\]
		in $\UC^+(x_{i_t}, y) \otimes_{A_{x_{i_t}}^0} \UC^-(x, x_{i_t})$ for each $t$. We call this expression the \textit{Reedy factorization} of $f$.
	\end{remark}
	
	The following two examples show that the notion of generalized $k$-linear Reedy categories unifies $k$-linear Reedy categories introduced in \cite{GS} and $k$-linearizations of generalized Reedy categories introduced in \cite{BM}.
	
	\begin{example} \label{k-linear Reedy categories}
		Let $k$ be a field. Recall from \cite[Definitions 3.1 and 3.6]{GS} that $\UC$ is called a $k$-linear Reedy category with respect to a degree function $d: \Ob(\UC) \to \lambda$ if it satisfies the following conditions:
		\begin{prt}
			\item $\UC$ has a subcategory $\UC^+$ sharing the same objects such that $\UC^+(x, x) \cong k$ and $\UC^+(x, y) \neq 0$ only if $x = y$ or $d(x) < d(y)$;
			
			\item $\UC$ has a subcategory $\UC^-$ sharing the same objects such that $\UC^-(x, x) \cong k$ and $\UC^-(x, y) \neq 0$ only if $x = y$ or $d(x) > d(y)$;
			
			\item for $x, y \in \Ob(\UC)$, the composition induces a $k$-linear isomorphism
			\[
			\bigoplus_{z \in \Ob(\UC)} \UC^+(z, y) \otimes_k \UC^-(x, z) \cong \UC(x, y).
			\]
		\end{prt}
		Clearly, a $k$-linear Reedy category is a generalized $k$-linear Reedy category satisfying two extra conditions: $k$ is a field, and $A_x^0 \cong k$ for $x \in \Ob(\UC)$.
	\end{example}
	
	\begin{example} \label{generalized Reedy categories}
		Recall from \cite[Definition 1.1]{BM} that a small skeletal category $\C$ is called a generalized Reedy category with respect to a degree function $d: \Ob(\C) \to \lambda$ if $\C$ satisfies the following conditions:
		\begin{prt}
			\item $\C$ has a subcategory $\C^+$ sharing the same objects such that $\C^+(x, y) \neq \emptyset$ only if $x = y$ or $d(x) < d(y)$;
			
			\item $\C$ has a subcategory $\C^-$ sharing the same objects such that $\C^-(x, y) \neq \emptyset$ only if $x = y$ or $d(x) > d(y)$;
			
			\item for $x \in \Ob(\C)$, $\C^+(x, x) = \C^-(x, x)$ is the group of automorphisms on $x$;
			
			\item every morphism $f$ in $\C$ has a factorization $f = f^+ f^-$, where $f^+$ is a morphism in $\C^+$ and $f^-$ is a morphism in $\C^-$, and this factorization is unique up to automorphisms;
			
			\item $\C^-(y, y)$ acts freely on $\C^-(x, y)$ for $x, y \in \Ob(\C)$.
		\end{prt}
		Let $k$ be a commutative ring. Then the $k$-linearization $\underline{\C}$ is a generalized $k$-linear Reedy category satisfying two extra conditions: $A_x^0$ is a group algebra for $x \in \Ob(\C)$, and $\underline{\C}(x, y)$ is a left free $A_y^0$-module for $x, y \in \Ob(\C)$.
	\end{example}

	\subsection{The induction functor}
	
	Since $\UC^-$ is a subcategory of $\UC$, the inclusion functor $\iota: \UC^- \to \UC$ induces a \textit{restriction functor} $\iota^{\ast}: \UC \Mod \to \UC^- \Mod$. It has a left adjoint given by the left Kan extension, called the \textit{induction functor}, which is nothing but the tensor functor $\UC \otimes_{\UC^-} -$, and a right adjoint given by the right Kan extension, called the \textit{coinduction functor}. In this paper we mainly consider the induction functor, and particularly its exactness.
	
	\begin{lemma} \label{projectivity}
		Let $\UC$ be a generalized $k$-linear Reedy category. One has:
		\begin{prt}
			\item $\UC^-$ is a projective left $\UC^0$-module (resp., projective right $\UC^0$-module) if and only if for all $x, y \in \Ob(\UC)$, $\UC^-(x, y)$ is a projective left $A_y^0$-module (resp., projective right $A_x^0$-module). A similar conclusion holds for $\UC^+$.
			
			\item If $\UC^-$ is a projective left $\UC^0$-module, then $\UC$ as a left $\UC^+$-module is also projective. Dually, if $\UC^+$ is a projective right $\UC^0$-module, then $\UC$ as a right $\UC^-$-module is projective.
		\end{prt}
	\end{lemma}
	
	\begin{prf*}
		(a) Note that $\UC^0$ is a direct sum of $A_x^0$'s for all $x \in \Ob(\UC)$. Consequently, as a left $\UC^0$-module, one has the following decomposition
		\[
		\UC^- = \bigoplus_{x, y \in \Ob(\UC )} \UC^-(x, y).
		\]
		Thus, $\UC^-$ is a projective left $\UC^0$-module if and only if each $\UC^-(x, y)$ is a projective left $A_y^0$-module. The other cases can be proved similarly.
		
		(b) By Remark \ref{ring isomorphism} and Lemma \ref{Reedy factorization}, for every $x \in \Ob(\UC)$, one has
		\[
		\UC(x, -) \cong \bigoplus_{z \in \Ob(\UC)} \UC^+(z, -) \otimes_{A_z^0} \UC^-(x, z)
		\]
		as left $\UC^+$-modules. If $\UC^-$ is a projective left $\UC^0$-module, by (a), $\UC^-(x, z)$ is a projective left $A_z^0$-module, so the tensor product $\UC^+(z, -) \otimes_{A_z^0} \UC^-(x, z)$ is a projective left $\UC^+$-module. Consequently, $\UC(x, -)$ is a projective left $\UC^+$-module for $x \in \Ob(\UC)$, so is $\UC$.
		
		To prove the second statement, we use the isomorphism
		\[
		\UC(-, x) \cong \bigoplus_{z \in \Ob(\UC)} \UC^+(z, x) \otimes_{A_z^0} \UC^-(-, z)
		\]
		of right $\UC^-$-modules.
	\end{prf*}
	
	An immediate consequence of this lemma is:
	
	\begin{corollary} \label{exactness of induction}
		Let $\UC$ be a generalized $k$-linear Reedy category. If $\UC^+$ is a projective right $\UC^0$-module, then the functor $\UC \otimes_{\UC^-} -$ is exact. Dually, if $\UC^-$ is a projective left $\UC^0$-module, then the functor $- \otimes_{\UC^+} \UC$ is exact.
	\end{corollary}
	
	\section{Standard modules and filtrations}
	\noindent
	Standard modules introduced in \cite{GS} play a key role for studying representations of $k$-linear Reedy categories. The main goal of this section is to generalize many important results in \cite{GS} to the framework of generalized $k$-linear Reedy categories.
	
	\begin{setup*}
		In the rest of this part, let $\UC$ be a generalized $k$-linear Reedy category with $d: \Ob(\UC) \to \lambda$ its degree function.
	\end{setup*}
	
	\subsection{Ideals and quotient categories}
	
	Let $\alpha$ be an ordinal with $\alpha \leqslant \lambda$. Define
	\[
	\mathfrak{I}_{\alpha}^- = \bigoplus_{d(y) < \alpha} \UC^-(-, y).
	\]
	This is clearly a right ideal of $\UC^-$. Furthermore, for any $f: x \to y$ in $\mathfrak{I}_{\alpha}^-$ and any morphism $g: y \to z$ in $\UC^-$, the composite $gf: x \to z$ also satisfies $d(z) \leqslant d(y) < \alpha$, so $\mathfrak{I}_{\alpha}^-$ is also a left ideal of $\UC^-$, and hence, a two-sided ideal of $\UC^-$. Consequently,
	\[
	\mathfrak{I}_{\alpha}^- (x, -) = \bigoplus_{d(y) < \alpha} \UC^-(x, y)
	\]
	is a left $\UC^-$-module for each $x \in \Ob(\UC)$.
	
	Similarly, we define
	\begin{align*}
	\mathfrak{I}_{\alpha}
	&= \bigoplus_{d(z) < \alpha} \UC^+(z, -) \otimes_{A_z^0} \UC^-(\bullet, z)\ \mathrm{and} \\
	\mathfrak{I}_{\alpha} (x, -)
    &= \bigoplus_{d(z) < \alpha} \UC^+(z, -) \otimes_{A_z^0} \UC^-(x, z)
	\end{align*}
	for a fixed object $x$.
	
	The following lemma asserts that $\mathfrak{I}_{\alpha}$ is a two-sided ideal of $\UC$,\footnote{To be more precisely, we shall say that the image of $\mathfrak{I}_{\alpha} = \bigoplus_{d(z) < \alpha} \UC^+(z, -) \otimes_{A_z^0} \UC^-(\bullet, z)$ under the composition map $\rho$ is a two-sided ideal of $\UC$; see Remark \ref{ring isomorphism}.}  and can be obtained from $\mathfrak{I}_{\alpha}^-$ via applying the induction functor.
	
	\begin{lemma} \label{2-sided ideals}
		The set $\mathfrak{I}_{\alpha}$ is a two-sided ideal of $\UC$. Furthermore, as left $\UC$-modules, one has
		\[
		\mathfrak{I}_{\alpha} \cong \UC \otimes_{\UC^-} \mathfrak{I}_{\alpha}^-.
		\]
	\end{lemma}
	
	\begin{prf*}
		To establish the first statement, it is enough to show: for $x, y \in \Ob(\UC)$, a morphism
		\[
		f \in \mathfrak{I}_{\alpha} (x, y) = \bigoplus_{d(z) < \alpha} \UC^+(z, y) \otimes_{A_z^0} \UC^-(x, z),
		\]
		and morphisms $g: u \to x$ and $h: y \to v$ in $\UC$, one has $h f \in \mathfrak{I}_{\alpha} (x, v)$ and $fg \in \mathfrak{I}_{\alpha} (u, y)$. We only show $hf \in \mathfrak{I}_{\alpha} (x, v)$ as the second one can be established similarly.
		
		By the Reedy factorization, one can find finitely many objects $z_1, \ldots, z_n$ (which might not be distinct) with $d(z_i) < \alpha$ such that
		$f = f_1^+ f_1^- + \ldots + f_n^+ f_n^-$,
		where $f_i^- \in \UC^-(x, z_i)$ and $f_i^+ \in \UC^+(z_i, y)$. Thus it suffices to show $h f_i^+ f_i^- \in \mathfrak{I}_{\alpha} (x, v)$ for each $i$. Similarly, we have another Reedy factorization
		$h = h_1^+ h_1^- + \ldots + h_m^+ h_m^-$.
		It reduces to show that $h_j^+ h_j^- f_i^+ f_i^- \in \mathfrak{I}_{\alpha} (x, v)$ for each $i$ and $j$. But for each fixed $i$ and $j$ we have
		\[
		h_j^- f_i^+ = l_1^+ l_1^- + \ldots + l_r^+ l_r^-,
		\]
		so it reduces further to show $h_j^+ l_s^+ l_s^- f_i^- \in \mathfrak{I}_{\alpha} (x, v)$ for each $i$, $j$ and $s$. But from the diagram
		\[
		\xymatrix{
			x \ar[r]^-{f_i^-} & z_i \ar[r]^-{l_s^-} & a_s \ar[r]^-{l_s^+} & w_j \ar[r]^-{h_j^+} & v
		}
		\]
		one clearly has $d(a_s) \leqslant d(z_i) < \alpha$, so indeed $h_j^+ l_s^+ l_s^- f_i^- \in \mathfrak{I}_{\alpha} (x, v)$. This finishes the proof.
		
		Now we prove the isomorphism. We just proved that $\mathfrak{I}_{\alpha} (x, -)$ is a left $\UC$-module. Since
		\[
		\mathfrak{I}_{\alpha} = \bigoplus_{x \in \Ob(\UC)} \mathfrak{I}_{\alpha} (x, -) \quad {\rm and} \quad
        \mathfrak{I}^-_{\alpha} = \bigoplus_{x \in \Ob(\UC)} \mathfrak{I}^-_{\alpha} (x, -),
		\]
		it suffices to check the $\UC$-module isomorphism
		$\mathfrak{I}_{\alpha} (x, -) \cong \UC \otimes_{\UC^-} \mathfrak{I}^-_{\alpha} (x, -)$.
		But by Remark \ref{ring isomorphism} we have
		\[
		\UC \otimes_{\UC^-} \mathfrak{I}^-_{\alpha} (x, -) \cong (\UC^+ \otimes_{\UC^0} \UC^-) \otimes_{\UC^-} \mathfrak{I}^-_{\alpha} (x, -) \cong \UC^+ \otimes_{\UC^0} \mathfrak{I}^-_{\alpha} (x, -).
		\]
		It remains to verify the isomorphism
		\[
		\mathfrak{I}_{\alpha} (x, -) = \bigoplus_{d(z) < \alpha} \UC^+(z, -) \otimes_{A_z^0} \UC^-(x, z) \cong \UC^+ \otimes_{\UC^0} \mathfrak{I}^-_{\alpha} (x, -).
		\]
		This can be checked via showing that the images of both sides under the composition map $\rho$ coincide. Here we prove it using the observation that $\UC^0$ is a direct sum of $A_z^0$'s. Indeed, one has
		\begin{align*}
			\UC^+ \otimes_{\UC^0} \mathfrak{I}^-_{\alpha} (x, -) & = \UC^+ \otimes_{\UC^0} (\bigoplus_{d(z) < \alpha} \UC^-(x, z))\\
			& \cong \bigoplus_{d(z) < \alpha} \UC^+(z, -) \otimes_{\UC^0} \UC^-(x, z)\\
			& \cong \bigoplus_{d(z) < \alpha} \UC^+(z, -) \otimes_{A_z^0} \UC^-(x, z).
		\end{align*}
		This finishes the proof.
	\end{prf*}
	
	With respect to $\mathfrak{I}_{\alpha}$, one can define a quotient category $\overline{\UC}_{\alpha}$ whose objects are those $x \in \Ob(\UC)$ with $d(x) \geqslant \alpha$. For two objects $x, y$ with $d(x) \geqslant \alpha$ and $d(y) \geqslant \alpha$, one sets
	\begin{align*}
		\overline{\UC}_{\alpha} (x, y) & = \UC(x, y) / \mathfrak{I}_{\alpha} (x, y)\\
		& \cong (\bigoplus_{z \in \Ob(\UC)} \UC^+(z, y) \otimes_{A_z^0} \UC^-(x, z)) / (\bigoplus_{d(z) < \alpha} \UC^+(z, y) \otimes_{A_z^0} \UC^-(x, z))\\
		& \cong \bigoplus_{d(z) \geqslant \alpha} \UC^+(z, y) \otimes_{A_z^0} \UC^-(x, z).
	\end{align*}
	
	From this construction one immediately has:
	
	\begin{lemma} \label{quotient categories}
		The quotient category $\overline{\UC}_{\alpha}$ is also a generalized $k$-linear Reedy category.
	\end{lemma}
	
	Fix an ordinal $\alpha \leqslant \lambda$. If $d(x) < \alpha$ and $d(y) < \alpha$, then
	\[
	\UC(x, y) = \bigoplus_{z \in \Ob(\UC)} \UC^+(z, y) \otimes_{A_z^0} \UC^-(x, z) = \bigoplus_{d(z) < \alpha} \UC^+(z, y) \otimes_{A_z^0} \UC^-(x, z)
	\]
	since $\UC^-(x, z) = 0$ when $d(z) \geqslant \alpha$. Thus, one can define a full subcategory $\UC_{\alpha}$ of $\UC$ consisting of objects $x$ with $d(x) < \alpha$, which is a generalized $k$-linear Reedy category by the above equality. Similarly, for an object $x$, we can define a full subcategory $\UC_x$ of $\UC$ consisting of objects $y$ with $y \prec x$, which is also a generalized $k$-linear Reedy category.

	\subsection{Standard modules} In this subsection we define standard modules and study their homological properties. These special modules were introduced in \cite{GS}.
	
	\begin{definition}
		For $x \in \Ob(\UC)$ with $d(x) = \alpha$, the \emph{standard} left $\UC$-module $\Delta_x$ is defined as
		\[
		\Delta_x = \UC(x, -) / \mathfrak{I}_{\alpha}(x, -) \cong \UC/\mathfrak{I}_{\alpha} \otimes_{\UC} \UC 1_x = \overline{\UC}_{\alpha} \otimes_{\UC} \UC 1_x,
		\]
		where $\UC$ is viewed as a $k$-algebra in the last two terms; see Remark \ref{ring isomorphism}. Dually, using corresponding right $\UC$-modules, one can define the standard right $\UC$-module
		\[
		\Delta^x = \UC(-, x) / \mathfrak{I}_{\alpha} (-, x) \cong 1_x \UC \otimes_{\UC} \UC/\mathfrak{I}_{\alpha} = 1_x \UC \otimes_{\UC} \overline{\UC}_{\alpha}.
		\]
	\end{definition}
	
	We can use partial orders $\preccurlyeq$ and $\preccurlyeq'$ on $\Ob(\UC)$ to reformulate standard modules. Define
	\begin{align*}
	\mathfrak{I}_x  &= \bigoplus_{z \prec x} \UC^+(z, -) \otimes_{A_z^0} \UC^-(x, z),\ \mathrm{and} \\
	\mathfrak{I}_x' &= \bigoplus_{z \prec' x} \UC^+(z, x) \otimes_{A_z^0} \UC^-(-, z).
    \end{align*}
	
	\begin{lemma} \label{identification of ideals}
		Let $x$ be an object in $\UC$ with $d(x) = \alpha$. Then $\mathfrak{I}_x = \mathfrak{I}_{\alpha} (x, -)$ and $\Delta_x = \UC 1_x / \mathfrak{I}_x$. Similarly, one has $\mathfrak{I}_x' = \mathfrak{I}_{\alpha} (-, x)$ and $\Delta^x = 1_x\UC / \mathfrak{I}_x'$.
	\end{lemma}
	
	\begin{prf*}
		We only prove the first statement. By definitions, it suffices to show
		\[
		\bigoplus_{z \prec x} \UC^+(z, -) \otimes_{A_z^0} \UC^-(x, z) = \bigoplus_{d(z) < \alpha} \UC^+(z, -) \otimes_{A_z^0} \UC^-(x, z).
		\]
		Clearly, if $z \prec x$, then $d(z) < \alpha$, so the left side is contained in the right side. On the other hand, if $\UC^-(x, z) \neq 0$, then $z \preccurlyeq x$. But $d(z) < \alpha = d(x)$, so $z \prec x$, and hence, the right side is also contained in the left side.
	\end{prf*}
	
	\begin{remark} \label{identification of standard modules}
		It is clear that $\Delta_x$ is only supported on objects $y$ with $d(y) \geqslant d(x) = \alpha$. Since the quotient category $\overline{\UC}_{\alpha}$ is also a generalized $k$-linear Reedy category by Lemma \ref{quotient categories}, one may view $\overline{\UC}_{\alpha} \Mod$ as a full subcategory of $\UC \Mod$. With respect to this identification, $\Delta_x$ is isomorphic to $\overline{\UC}_{\alpha} (x, -)$. But $x$ is a minimal object in $\overline{\UC}_{\alpha}$, so one has
		\[
		\overline{\UC}_{\alpha} (x, -) \cong \overline{\UC}_{\alpha}^+ \otimes_{\overline{\UC}_{\alpha}^0} \overline{\UC}_{\alpha}^- (x, -) = \overline{\UC}_{\alpha}^+ \otimes_{\overline{\UC}_{\alpha}^0} \overline{\UC}_{\alpha}^- (x, x) \cong \overline{\UC}_{\alpha}^+(x, -) \cong \UC^+(x, -).
		\]
		Therefore, we may identify $\Delta_x$, $\UC^+(x, -)$ and $\overline{\UC}_{\alpha}(x, -)$ as $\UC_0$-modules. Dually, we can identify $\Delta^x$, $\UC^-(-, x)$ and left $\overline{\UC}_{\alpha}(-, x)$ as right $\UC_0$-modules.
	\end{remark}
	
	It is clear that $A_x^0 \cong \UC^-(x, -) / \mathfrak{I}_{\alpha}^-(x, -)$ when $d(x) = \alpha$, so one can view $A_x^0$ as a left $\UC^-$-module concentrated on the single object $x$. We show under a mild condition that standard modules are induced from $A_x^0$, and have nice homological properties, generalizing \cite[Theorems 4.3 and 4.7]{GS}.
	
	\begin{theorem} \label{properties of standard modules}
		Suppose that $\UC^+$ is a projective right $\UC^0$-module. Then:
		\begin{prt}
			\item $\Delta_x \cong \UC \otimes_{\UC^-} A_x^0$ as left $\UC$-modules;
			
			\item $\End_{\UC} (\Delta_x) \cong (A_x^0)^{\op}$ as $k$-algebras;
			
			\item $\Hom_{\UC} (\Delta_x, \Delta_y) \neq 0$ only if $x = y$ or $d(x) > d(y)$;
			
			\item $\Ext_{\UC}^n (\Delta_x, \Delta_y) \neq 0$ for $n \geqslant 1$ only if $d(x) > d(y)$.
		\end{prt}
	\end{theorem}
	
	\begin{prf*}
		(a) Let $d(x)=\alpha$. Applying the functor $\UC \otimes_{\UC^-} -$, which is exact by Corollary \ref{exactness of induction}, to the exact sequence
		\[
		0 \to \mathfrak{I}^-_{\alpha} (x, -) \to \UC^-(x, -) \to A_x^0 \to 0
		\]
		of left $\UC^-$-modules, we get another short exact
		\[
		0 \to \UC \otimes_{\UC^-} \mathfrak{I}^-_{\alpha} (x, -) \longrightarrow \UC \otimes_{\UC^-} \UC^-(x, -) \longrightarrow \UC \otimes_{\UC^-} A_x^0 \to 0
		\]
		of left $\UC$-modules. But $\UC \otimes_{\UC^-} \UC^-(x, -) \cong \UC(x, -)$ and $\UC \otimes_{\UC^-} \mathfrak{I}^-_{\alpha} (x, -) \cong \mathfrak{I}_{\alpha} (x, -)$ by Lemma \ref{2-sided ideals}. The conclusion follows.
		
		(b) Let $d(x) = \alpha$. Note that $\overline{\UC}_{\alpha}$ is a quotient category of $\UC$, so $\overline{\UC}_{\alpha} \Mod$ is a full subcategory of $\UC \Mod$. Moreover, $\Delta_x$ viewed as a left $\overline{\UC}_{\alpha}$-module coincides with $\overline{\UC}_{\alpha} (x, -)$. Thus, one has
		$$\End_{\UC} (\Delta_x) \cong \End_{\overline{\UC}_{\alpha}} (\overline{\UC}_{\alpha} (x, -)) \cong (1_x \overline{\UC}_{\alpha} 1_x)^{\op} = (A_x^0)^{\op}.$$
		
		(c) We first consider the case that $d(y) = \beta > \alpha = d(x)$. We want to show that $\Hom_{\UC} (\Delta_x, \Delta_y)$ vanishes. Since $\overline{\UC}_{\alpha}$ is again a generalized $k$-linear Reedy category, viewed as $k$-algebras one has
		\[
		\overline{\UC}_{\beta} = \UC / \mathfrak{I}_{\beta} \cong (\UC / \mathfrak{I}_{\alpha}) / (\mathfrak{I}_{\beta} / \mathfrak{I}_{\alpha}).
		\]
		Thus, one may view $\overline{\UC}_{\beta}$ as a quotient algebra of $\overline{\UC}_{\alpha}$, and $\mathfrak{I}_{\beta} / \mathfrak{I}_{\alpha}$ is precisely the two-sided ideal of $\overline{\UC}_{\alpha}$ to define this quotient; see the construction before Lemma \ref{quotient categories}. Moreover, $\overline{\UC}_{\alpha} \Mod$ is a full subcategory of $\UC \Mod$ and it contains $\Delta_x$ and $\Delta_y$, so
		\[
		\Hom_{\UC} (\Delta_x, \Delta_y) = \Hom_{\overline{\UC}_{\alpha}} (\Delta_x, \Delta_y) \cong \Hom_{\overline{\UC}_{\alpha}} (\overline{\UC}_{\alpha} (x, -), \Delta_y),
		\]
		as $\Delta_x \cong \overline{\UC}_{\alpha} (x, -)$ in this full subcategory. But the last term is isomorphic to the value of $\Delta_y$ on $x$, which is clearly 0 since $\Delta_y$ is only supported on objects $z$ with $d(z) \geqslant d(y) = \beta$.
		
		Then we consider the case that $d(y) = d(x)$ but $x \neq y$. As in the previous paragraph, one has
		\[
		\Hom_{\UC} (\Delta_x, \Delta_y) = \Hom_{\overline{\UC}_{\alpha}} (\Delta_x, \Delta_y) \cong \Hom_{\overline{\UC}_{\alpha}} (\overline{\UC}_{\alpha}(x, -), \overline{\UC}_{\alpha}(y, -)) \cong \overline{\UC}_{\alpha}(y, x) = 0
		\]
		by Lemma \ref{Reedy factorization} since $x$ and $y$ are minimal in $\overline{\UC}_{\alpha}$.
		
		(d) By (a), $\Delta_x \cong \UC \otimes_{\UC^-} A_x^0$. By Eckmann-Shapiro's lemma (see \cite[p47]{Ben}, which also holds for rings with a complete family of orthogonal idempotents), one has
		\[
		\Ext_{\UC}^n (\Delta_x, \Delta_y) \cong \Ext_{\UC}^n (\UC \otimes_{\UC^-} A_x^0, \Delta_y) \cong \Ext_{\UC^-}^n (A_x^0, \Delta_y).
		\]
		We want to show that this extension group vanishes whenever $d(x) \leqslant d(y)$.
		
		Construct a projective resolution $\ldots \to P^ 1 \to P^0 \to A_x^0 \to 0$ in $\UC^- \Mod$ with $P^0 = \UC^-(x, -)$. By the structure of $\UC^-$, for $n \geqslant 1$, one can assume that the value of $P^i$ on each object $z$ is 0 if $d(z) \geqslant d(x)$. In other words, $P^n$ is supported on objects $z$ with $d(z) < d(x) \leqslant d(y)$. On the other hand, we know that $\Delta_y$ is supported on objects $z$ with $d(z) \geqslant d(y)$. Consequently, $\Hom_{\UC^-} (P^i, \Delta_y) = 0$ for all $i \geqslant 1$, and the conclusion follows.
	\end{prf*}
	
	\begin{remark} \label{opposite category}
		When $\UC$ is a generalized $k$-linear Reedy category, so is the opposite category $\UC^{\op}$ with $(\UC^{\op})^+ = (\UC^-)^{\op}$ and $(\UC^{\op})^- = (\UC^+)^{\op}$. In this case, the right standard module $\Delta^x$ of $\UC$ is a left standard module over $\UC^{\op}$. Therefore, all results established in this subsection hold for $\Delta^x$ with suitable modifications. In particular, we obtain a dual version of the above theorem.
	\end{remark}
	
	\begin{theorem} \label{dual version}
		Suppose that $\UC^-$ is a projective left $\UC^0$-module. Then one has:
		\begin{prt}
			\item $\Delta^x \cong A_x^0 \otimes_{\UC^+} \UC$ as right $\UC$-modules;
			
			\item $\End_{\UC} (\Delta^x) \cong A_x^0$ as $k$-algebras;
			
			\item $\Hom_{\UC} (\Delta^x, \Delta^y) \neq 0$ only if $x = y$ or $d(x) > d(y)$;
			
			\item $\Ext_{\UC}^n (\Delta^x, \Delta^y) \neq 0$ for $n \geqslant 1$ only if $d(x) > d(y)$.
		\end{prt}
	\end{theorem}
	

	\subsection{Filtrations and irreducible modules}
	
	It was shown in \cite{GS} that a $k$-linear Reedy category is similar to a quasi-hereditary algebra; that is, every representable functor has a filtration by standard modules, see \cite[Theorem 4.5]{GS}. Moreover, standard modules are used to classify all irreducible modules; see \cite[Theorem 4.14]{GS}. In this subsection we extend these results to generalized $k$-linear Reedy categories.
	
	Since a classification of irreducible modules of an arbitrary $k$-algebra in general is not available, in this subsection we impose an extra condition: each $A_x^0$ is a semi-perfect $k$-algebra. Therefore, for each $x \in \Ob(\UC)$, one can fix a finite complete set $E_x$ of orthogonal primitive idempotents in $A_x^0$. We write $e\sim_x e'$ if $A_x^0e\cong A_x^0e'$ as left $A_x^0$-modules. The quotient set $E_x/{\sim_x}$ parametrizes the isomorphism classes of simple left $A_x^0$-modules. In particular, every indecomposable projective left $A_x^0$-module is isomorphic to $A_x^0 e$ for a certain $e \in E_x$. Since $\UC e$ is a projective left $\UC$-module, one can define
	$\Delta_{x, e} = \UC e/ \mathfrak{I}_{\alpha}e$,
	where $\alpha = d(x)$. Consequently, $\Delta_{x, e}$ is a direct summand of $\Delta_x$, and
	$\Delta_x = \bigoplus_{e \in E_x} \Delta_{x, e}$.
	Moreover, $\Delta_{x, e}$ is indecomposable, since
	$\End_{\UC}(\Delta_{x, e}) \cong \End_{\overline{\UC}_{\alpha}} (\overline{\UC}_{\alpha} e) \cong (eA_x^0e)^{\op}$
	is a local ring.
	
	\begin{lemma} \label{filtration of C^-modules}
		Suppose that each $A_x^0$ is a semi-perfect $k$-algebra for $x \in \Ob(\UC)$. Let $V$ be a left $\UC^-$-module such that $V(x)$ is a projective left $A_x^0$-module for each $x \in \Ob(\UC)$. Then $V$ admits a transfinite filtration
		\[
		0 = V_{<0} \subseteq V_{<1} \subseteq \ldots \subseteq V_{<\lambda + 1} = V
		\]
		such that for $0\leqslant\alpha \leqslant \lambda$, each $V_{\alpha} = V_{<\alpha + 1} / V_{<\alpha}$ is isomorphic to a direct sum of members in the family
		\[
		\mathfrak{F}_{\alpha} = \{ A_x^0e \mid d(x) = \alpha, \, e \in E_x \}.
		\]
	\end{lemma}
	
	\begin{prf*}
		For each $\alpha$, define a submodule $V_{<\alpha}$ of $V$ as follows: for $x \in \Ob(\UC)$,
		\[
		V_{< \alpha} (x) = \begin{cases}
			V(x)  & \text{if } d(x) < \alpha;\\
			0  & \text{else.}
		\end{cases}
		\]
		Clearly, one has $V \cong \varinjlim_{\alpha} V_{<\alpha}$. Furthermore, one has
		\[
		V_{\alpha} (x) = V_{<\alpha + 1}(x) / V_{<\alpha}(x) = \begin{cases}
			V(x)  & \text{if } d(x) = \alpha;\\
			0 & \text{else.}
		\end{cases}
		\]
		Note that if $d(x) = \alpha = d(y)$ but $x \neq y$, then $\UC^-(x, y) = 0 = \UC^-(y, x)$. Thus,
		$V_{\alpha} = \bigoplus_{d(x) = \alpha} V(x)$
		as left $\UC^-$-modules. But since $V(x)$ is a projective left $A_x^0$-module, we get a further decomposition
		\[
		V_{\alpha} = \bigoplus_{d(x) = \alpha} V(x) \cong \bigoplus_{d(x) = \alpha} \bigoplus_{e \in E_x} (A_x^0 e)^{c_{x, e}},
		\]
		where $c_{x, e}$ is the multiplicity which can be infinity.
	\end{prf*}
	
	\begin{proposition} \label{filtration}
		Suppose that $\UC^+$ is a projective right $\UC^0$-module,  $\UC^-$ is a projective left $\UC^0$-module, and each $A_x^0$ is a semi-perfect $k$-algebra for $x \in \Ob(\UC)$. Then every $\UC(x, -)$ can be filtered by a coproduct of $\Delta_{z, e}$ with $d(z) \leqslant d(x)$ and $e \in E_z$.
	\end{proposition}
	
	\begin{prf*}
		Let $x\in\Ob(\UC)$ with $d(x)=\beta$ and let $V = \UC^-(x, -)$. Since $\UC^-$ is a projective left $\UC^0$-module, $V(y)$ is a projective left $A_y^0$-module for each $y\in\Ob(\UC)$ by Lemma \ref{projectivity}(a). Applying Lemma \ref{filtration of C^-modules}, we obtain a transfinite filtration
		\[
		0 = V_{<0} \subseteq V_{<1} \subseteq \ldots \subseteq V_{<\beta + 1} = V
		\]
		such that for $\alpha \leqslant \beta$, each $V_{\alpha} = V_{<\alpha + 1} / V_{<\alpha}$ is a coproduct of $A_z^0 e$ with $d(z) = \alpha$ and $e \in E_z$. Since $\UC^+$ is a projective right $\UC^0$-module, $\UC \otimes_{\UC^-} -$ is exact by Corollary \ref{exactness of induction}. Applying it to the above filtration and observing that $\Delta_{z, e} \cong \UC \otimes_{\UC^-} A_z^0 e$ by Theorem \ref{properties of standard modules}, we establish the conclusion.
	\end{prf*}
	
	\begin{remark} \label{a general filtration}
		Suppose that $\UC^+$ is a projective right $\UC^0$-module. Without the assumption that each $A_x^0$ is a semi-perfect $k$-algebra, using the same argument, one can get a filtration of $\UC(x, -)$ whose factors are
		$\bigoplus_{d(z) = \alpha} \UC \otimes_{\UC^-} \UC^-(x, z)$
		for each $\alpha \leqslant d(x)$, where $\UC^-(x, z)$ is viewed as a left $\UC^-$-module concentrated on the object $z$. Moreover, one has the following isomorphisms of left $\UC$-modules:
		\begin{align*}
			\UC \otimes_{\UC^-} \UC^-(x, z) & \cong \UC \otimes_{\UC^-} (A_z^0 \otimes_{A_z^0} \UC^-(x, z))\\
			& \cong (\UC \otimes_{\UC^-} A_z^0) \otimes_{A_z^0} \UC^-(x, z)\\
			& \cong \Delta_z \otimes_{A_z^0} \UC^-(x, z).
		\end{align*}
	\end{remark}
	
	Given a left $\UC$-module $M$, denote by $[M]$ its isomorphism class. The following result gives parameterizations of isomorphism classes of irreducible left $\UC$-modules under some mild assumptions.
	
	\begin{theorem} \label{parameterize irreducibles}
		Suppose that $\UC^+$ is a projective right $\UC^0$-module, $\UC^-$ is a projective left $\UC^0$-module, and each $A_x^0$ is a semi-perfect $k$-algebra. Then there are bijections between the following sets:
		\begin{rqm}
			\item $\{[S] \mid S \text{ is an irreducible left $\UC$-module} \}$;
			
			\item the set of isomorphism classes $\{[\Delta_{x,e}]\mid x\in\Ob(\UC),\ e\in E_x\}$;
			
			\item the disjoint union $\displaystyle\bigsqcup_{x\in\Ob(\UC)}(E_x/{\sim_x})$.
		\end{rqm}
	\end{theorem}
	
	\begin{prf*}
First, for $x,x'\in\Ob(\UC)$ and $e\in E_x$, $e'\in E_{x'}$, we have
\[
\Delta_{x,e}\cong\Delta_{x',e'}
\quad\Longleftrightarrow\quad
x=x'\text{ and }A_x^0e\cong A_x^0e'.
\]
Indeed, an isomorphism of standard summands yields nonzero maps in both directions between $\Delta_x$ and $\Delta_{x'}$. By Theorem~\ref{properties of standard modules}(c), their degrees are equal, say $d(x)=d(x')=\alpha$. Both standard summands are modules over $\overline{\UC}_{\alpha}$. The objects $x$ and $x'$ have minimal degree in this quotient category, so Lemma~\ref{Reedy factorization}(c), applied to $\overline{\UC}_{\alpha}$, shows that there can be no nonzero morphism between the corresponding representables when $x\ne x'$. Consequently $x=x'$. Evaluation at $x$ gives $A_x^0e\cong A_x^0e'$. Conversely, an isomorphism of these projective $A_x^0$-modules induces an isomorphism of $\overline{\UC}_{d(x)}$-representables, hence of the standard summands. This proves the bijection between (2) and (3).

Next, every $\Delta_{x,e}$ has a unique maximal submodule. To see this, view it as the representable $\overline{\UC}_{\alpha}e$ with $\alpha=d(x)$. It is generated by the element $e$ at $x$, and its value at $x$ is the indecomposable projective $A_x^0e$, which has a unique maximal submodule $\operatorname{rad}(A_x^0e)$ since $A_x^0$ is semiperfect. Any proper subrepresentation has its value at $x$ contained in $\operatorname{rad}(A_x^0e)$; otherwise it contains $e$ and generates all of $\Delta_{x,e}$. Thus the sum of all proper subrepresentations is still proper and is the unique maximal one. Write $S_{x,e}$ for the resulting simple quotient. Its value at $x$ is the simple top of $A_x^0e$, while it vanishes on objects of degree smaller than $\alpha$ and on other objects of degree $\alpha$.

For any simple $\UC$-module $S$, choose an object $x$ of minimal degree $\alpha$ with $S(x)\ne0$. A nonzero map $\UC(x,-)\to S$ is surjective. By the Reedy factorization, its kernel contains $\mathfrak I_\alpha(x,-)$ because $S$ vanishes at every object of degree below $\alpha$. Hence it induces an epimorphism $\Delta_x\to S$. Since $\Delta_x=\bigoplus_{e\in E_x}\Delta_{x,e}$, some summand maps onto $S$, and $S\cong S_{x,e}$. This proves surjectivity.

Finally, if $S_{x,e}\cong S_{x',e'}$, their minimal degrees coincide. At that degree the support of $S_{x,e}$ consists of the single object $x$, and likewise for $S_{x',e'}$; hence $x=x'$. Evaluating at $x$ identifies the simple tops of $A_x^0e$ and $A_x^0e'$. Over a semiperfect ring the indecomposable projective covers of isomorphic simple modules are isomorphic, so $e\sim_x e'$. This proves injectivity and completes the proof.
\end{prf*}
	
	{\bf Throughout this paper}, for a small category $\C$, we denote its $k$-linearization by $\underline{\C}$.
	
	\begin{example}
		Let $k$ be a field. Let $\C$ be the category with objects $[n]$ for $n \in \mathbb{N}$ and morphisms the usual maps, whose representation theory has been considered in \cite{Pow, PV, Gor}. It is a generalized Reedy category with $\C^+$ the subcategory of injections and $\C^-$ the subcategory of surjections. By Theorem \ref{parameterize irreducibles}, isomorphism classes of irreducible $\underline{\C}$-modules are parameterized by pairs $(n, [L])$, where $n \in \mathbb{N}$ and $[L]$ is the isomorphism class of a simple $kS_n$-module $L$. In particular, if $k$ is a field of characteristic 0, then each irreducible $\underline{\C}$-module is isomorphic to the top of a standard module
		\[
		\Delta_{n, \mu} = \underline{\C}^+ \otimes_{kS_n} L_{\mu} \cong \underline{\C}^+([n], -) \otimes_{kS_n} L_{\mu},
		\]
		where $\mu$ is a partition of $n$, and $L_{\mu}$ is the corresponding irreducible $kS_n$-module. The reader can compare this result to \cite[Theorem 5.5]{Gor} and \cite[Theorem 2]{Pow}.
	\end{example}
	
	\begin{example}
		Let $k$ be a field. Let $\C$ be the category whose objects are cyclically ordered sets $[n]$ with $n \geqslant1$ elements and whose morphisms are order-preserving maps. This is a skeleton of the cyclic category $\Lambda$ introduced by Connes \cite{Connes}, which has been widely applied in algebraic topology and homology theory. In particular, people have obtained a Dold-Kan type theorem for this category; see for instance \cite{DK, Loday, Slo}. It is easy to check that $\C$ is a generalized Reedy category with $\C^+$ the subcategory of order-preserving injections and $\C^-$ the subcategory of order-preserving surjections.
		
		By Theorem \ref{parameterize irreducibles}, isomorphism classes of irreducible $\underline{\C}$-modules are parameterized by pairs $(n, [L])$, where $n \geqslant 1$ and $[L]$ is the isomorphism class of a simple $kC_n$-module $L$ with $C_n$ the cyclic group of order $n$. If $k$ is an algebraically closed field of characteristic zero, then each irreducible $\underline{\C}$-module is isomorphic to the top of a standard module
		$$\Delta_{n, \xi} = \underline{\C}^+ \otimes_{kC_n} L_{\xi} \cong \underline{\C}^+([n], -) \otimes_{kC_n} L_{\xi},$$
		where $\xi$ is an $n$-th root of unity, and $L_{\xi}$ is the corresponding irreducible $kC_n$-module.
	\end{example}

	\section{Decomposition of module categories I: central idempotents}
	\noindent
	As we mentioned in the introduction, there are many interesting combinatorial categories $\C$ such that $\underline{\C}$ is a generalized $k$-linear Reedy category, and $\underline{\C} \Mod$ is equivalent to $\underline{\C}^0 \Mod$. In this and the next section we give two sufficient criteria such that this equivalence holds for generalized $k$-linear Reedy categories.
	
	Let $\UC$ be a generalized $k$-linear Reedy category. For $x \in \Ob(\UC)$, recall that $A_x^0 = \UC^+(x, x) = \UC^-(x, x)$,
	\[
	A_x = \UC(x, x) = 1_x \UC 1_x \cong \bigoplus_{y \preccurlyeq x} \UC^+(y, x) \otimes_{A_y^0} \UC^-(x, y),
	\]
	where the isomorphism of $k$-algebras is given by the composition $\rho$ (see Remark \ref{ring isomorphism}), and
	\[
	\mathfrak{I}_x = \bigoplus_{y \prec x} \UC^+(y, -) \otimes_{A_y^0} \UC^-(x, y).
	\]
	The short exact sequence
	\begin{equation*}\tag{$\ast$}
		0 \to \mathfrak{I}_x \to \UC 1_x \to \Delta_x \to 0
	\end{equation*}
	of $\UC$-modules induces a short exact sequence
	\begin{equation*}\tag{$\ast\ast$}
		0 \to 1_x \mathfrak{I}_x \to 1_x \UC 1_x = A_x \to 1_x \Delta_x = A_x^0 \to 0.
	\end{equation*}
	It is routine to check that
	\[
	I_x = 1_x \mathfrak{I}_x = \bigoplus_{y \prec x} \UC^+(y, x) \otimes_{A_y^0} \UC^-(x, y)
	\]
	is a two-sided ideal of $A_x$. Thus, we may view $A^0_x$ as a quotient $k$-algebra of $A_x$.
	
	\begin{lemma} \label{split sequences}
	For all $x, y \in \Ob(\UC)$, suppose that the following conditions hold:
\begin{enumerate}
\item $\UC^+(y, x)$ is a right free $A_y^0$-module and $\UC^-(x, y)$ is a left free $A_y^0$-module;

\item $\UC^+(y, x) = 0$ if and only if $\UC^-(x, y) = 0$;

\item there exists a central idempotent $e_x \in A_x$ such that $I_x = e_x A_x e_x$.
\end{enumerate}
Then the short exact sequences $(\ast)$ and $(\ast \ast)$ split. Furthermore, $\mathfrak{I}_x \cong \UC e_x$ and $\Delta_x \cong \UC f_x$ as left $\UC$-modules, where $f_x = 1_x - e_x$.
	\end{lemma}
	
	\begin{prf*}
		Note that $e_x$ is contained in $I_x$ and serves as the identity of $I_x$ as a $k$-algebra. Thus, the short exact sequence $(\ast \ast)$ becomes to
		\[
		0 \to e_xA_xe_x \to (e_xA_xe_x \oplus f_xA_xf_x) \to A_x^0 \to 0
		\]
		which forces $A_x^0 \cong f_xA_xf_x$, and hence, $(\ast \ast)$ splits. We identify $A_x^0$ with $f_x A_x f_x$ and under this identification, the identity of this ring is $f_x$ rather than $1_x$. In particular, the value $\Delta_x(x)$ becomes $f_x A_x f_x$. Thus,
		\[
		\Hom_{\UC} (\UC e_x, \Delta_x) \cong e_x \Delta_x (x) = e_x f_xA_xf_x = 0
		\]
		and similarly $\Hom_{\UC} (\UC f_x, \mathfrak{I}_x) = 0$. Consequently, $\Delta_x$ is a quotient module of $\UC f_x$.
		
		Applying the snake lemma to the commutative diagram of exact rows
		\[
		\xymatrix{
			0 \ar[r] & \UC e_x \ar[r] \ar[d] & \UC 1_x \ar[r] \ar@{=}[d] & \UC f_x \ar[r] \ar[d] & 0\\
			0 \ar[r] & \mathfrak{I}_x \ar[r] & \UC 1_x \ar[r] & \Delta_x \ar[r] & 0,
		}
		\]
		where the third vertical arrow is an epimorphism, we deduce that $\UC e_x \subseteq \mathfrak{I}_x$ from the left vertical map. It suffices to show the reverse inclusion. But
\[
\mathfrak{I}_x \subseteq \UC 1_x = \UC e_x \oplus \UC f_x,
\]
we only need to show that $\mathfrak{I}_x f_x = 0$.

If $\mathfrak{I}_x f_x \neq 0$, then one can find a certain $y \in \Ob(\UC)$ with $y \prec x$ such that $\UC^-(x, y)f_x \neq 0$. In particular, $\UC^-(x, y) \neq 0$, and hence $\UC^+(y, x)$ is a nonzero right free $A_y^0$-module by (1) and (2). Consequently, one has
\[
\UC^+(y, x) \otimes_{A_y^0} \UC^-(x, y) f_x \neq 0.
\]
But $\UC^+(y, x) \otimes_{A_y^0} \UC^-(x, y)$ is contained in $I_x$, so $I_x f_x \neq 0$, which is impossible. The conclusion follows from this contradiction.
	\end{prf*}
	
	Under one more assumption, we show that the family $\{ \Delta_x \mid x \in \Ob(\UC) \}$ is orthogonal.
	
	\begin{lemma} \label{orthogonal family}
		Under the same assumptions as in Lemma \ref{split sequences}, one has $\Hom_{\UC} (\Delta_x, \Delta_y) \neq 0$ if and only if $x = y$.
	\end{lemma}
	
	\begin{prf*}
		Suppose that $\Hom_{\UC} (\Delta_x, \Delta_y) \neq 0$. By Lemma \ref{split sequences}, one can assume that $\Delta_x = \UC f_x$ and $\Delta_y = \UC f_y$, where $f_x$ and $f_y$ are central idempotents in $A_x$ and $A_y$ respectively. Thus,
		\[
		\Hom_{\UC} (\Delta_x, \Delta_y) \cong f_x \UC f_y \cong \bigoplus_{z} f_x \UC^+(z, x) \otimes_{A_z^0} \UC^-(y, z) f_y.
		\]
		If $y \neq z$, then $\UC^-(y, z) \in \mathfrak{I}_y$. But by the proof of Lemma \ref{split sequences}, $\mathfrak{I}_y f_y = \UC e_y f_y = 0$, so one has
		\[
		0 \neq \Hom_{\UC} (\Delta_x, \Delta_y) \cong f_x \UC f_y = f_x \UC^+(y, x) \otimes_{A_y^0} \UC^-(y, y) f_y.
		\]
		In particular, $f_x \UC^+(y, x) \neq 0$. By (2), one has $f_x \UC^+(y, x) \otimes_{A_y^0} \UC^-(x, y) \neq 0$. If $x \neq y$, then
\[
\UC^+(y, x) \otimes_{A_y^0} \UC^-(x, y) \subseteq I_x,
\]
so $f_x I_x = f_x e_x A_x e_x \neq 0$, which is impossible as $f_x e_x = 0$. Thus, $x$ must coincide with $y$.
	\end{prf*}
	
	\begin{remark}
		The conclusion fails without the condition (2) in Lemma \ref{orthogonal family}. For instance, if $\UC = \UC^+$, then $\Delta_x = \UC 1_x$ and $\Delta_y = \UC 1_y$, so $\Hom_{\UC} (\Delta_x, \Delta_y) \cong 1_x \UC 1_y = \UC(y,x)$, which might not be 0 even if $x \neq y$.
	\end{remark}
	
	The following lemma gives us a new family of projective generators in $\UC \Mod$.
	
	\begin{lemma} \label{new generators}
		Under the same assumptions as in Lemma \ref{split sequences}, $\{ \Delta_x \mid x \in \Ob(\UC) \}$ is a set of projective generators in $\UC \Mod$.
	\end{lemma}
	
	\begin{prf*}
		By Remark \ref{a general filtration}, for each $x \in \Ob(\UC)$, the representable functor $\UC(x, -)$ has a filtration with factors
		\[
		\bigoplus_{d(z) = \alpha} \UC \otimes_{\UC^-} \UC^-(x, z) \cong \bigoplus_{d(z) = \alpha} \Delta_z \otimes_{A_z^0} \UC^-(x, z)
		\]
		for $\alpha \leqslant d(x)$. Since $\UC^-(x, z)$ is a free left $A_z^0$-module by Lemma \ref{projectivity}(a), up to isomorphism $\Delta_z \otimes_{A_z^0} \UC^-(x, z)$ is a direct sum of copies of $\Delta_z$. In particular, it is projective. Thus, from the filtration we deduce that
		\[
		\UC(x, -) \cong \bigoplus_{\alpha \leqslant d(x)} \bigoplus_{d(z) = \alpha} \UC \otimes_{\UC^-} \UC^-(x, z) \cong \bigoplus_{\alpha \leqslant d(x)} \bigoplus_{d(z) = \alpha} \Delta_z \otimes_{A_z^0} \UC^-(x, z)
		\]
		as $\UC$-modules, so $\UC(x, -)$ can be generated by the given family. This finishes the proof.
	\end{prf*}

	Now it is ready to prove the main result in this section.
	
	\begin{theorem} \label{main theorem for decomposition}
		Under the same assumptions as in Lemma \ref{split sequences}, one has
		\[
		\UC \Mod \simeq \prod_{x \in \Ob(\UC)} A_x^0 \Mod.
		\]
		In particular, if each $A_x^0$ is a semisimple $k$-algebra, then $\UC \Mod$ is semisimple as well.
	\end{theorem}
	
	\begin{prf*}
		Lemma \ref{split sequences} tells us that each $\Delta_x$ is projective, so by Lemma \ref{new generators} they form a family of projective generators in $\UC \Mod$. Lemma \ref{orthogonal family} asserts that this family is orthogonal. The conclusion then follows from the Morita equivalence and the isomorphism $\End_{\UC} (\Delta_x) \cong (A_x^0)^{\op}$ proved in Theorem \ref{properties of standard modules}.
	\end{prf*}
	
	The following example has been studied in \cite{CEF, Kuhn, Slo}.
	
	\begin{example}
		Let $\mathbb{F}_q$ be a finite field of $q$ elements, and $k$ a commutative ring such that $q$ is invertible in $k$. Let $\C$ be the category consisting of objects $\mathbb{F}_q^n$ with $n \in \mathbb{N}$ and $\mathbb{F}_q$-linear maps. Then $\C$ is a generalized Reedy category with $\C^+$ the subcategory of linear injections and $\C^-$ the subcategory of linear surjections.
		
		For $n \leqslant m$, the general linear group $\mathrm{GL}_n (\mathbb{F}_q)$ acts freely on $\C^+(\mathbb{F}_q^n, \mathbb{F}_q^m)$ from the right side and freely on $\C^-(\mathbb{F}_q^m, \mathbb{F}_q^n)$ from the left side. Furthermore, the existence of central idempotents required in Lemma \ref{split sequences} was proved in \cite{Kov} (see also \cite[Theorem 2.1]{Kuhn}). Thus, by \cite[Theorem 1.1]{Kuhn} one has
		\[
		\underline{\C} \Mod \simeq \prod_{n \in \mathbb{N}} k\mathrm{GL}_n (\mathbb{F}_q) \Mod.
		\]
		
		In the corresponding set-theoretic analogue (informally denoted by $q=1$), let $\C_{\rm set}$ be the category of finite sets and partial injections; equivalently, its morphisms can be represented by basepoint-preserving maps which are injective on the inverse images of non-basepoint elements. See \cite[Section 4]{CEF} and \cite[Section 4]{Kuhn}. In this case, one has
		\[
		\underline{\C_{\rm set}} \Mod \simeq \prod_{n \in \mathbb{N}} kS_n \Mod,
		\]
		where $S_n$ is the symmetric group on $n$ letters; see \cite[Theorem 1.7]{CEF} and \cite[Theorem 4.1]{Kuhn}.
	\end{example}
	
	\section{Decomposition of module categories II: non-degenerate elements}
	\noindent
	The conditions in Theorem \ref{main theorem for decomposition}, in particular the existence of central idempotents, are hard to check in many situations. In this section we provide another criterion for a special family of generalized $k$-linear Reedy categories such that the conclusion of Theorem \ref{main theorem for decomposition} still holds.
	
	\subsection{Another sufficient criterion}
	
	Let $\UC \cong \UC^+ \otimes_{\UC^0} \UC^-$ be a generalized $k$-linear Reedy category. Given $x, y \in \Ob(\UC)$ and a nonzero $g \in \UC^+(x, y)$, the composition $\rho$ gives rise to a $k$-linear map
	\[
	\varphi_g: \UC^-(y, x) \longrightarrow A_x = \UC(x, x) \longrightarrow A_x^0, \quad f \mapsto fg \mapsto \overline{fg},
	\]
where the second arrow is the natural quotient map since $A_x^0$ is a quotient algebra of $A_x$. We say that $g$ is \textit{non-degenerate} if $\varphi_g$ is nonzero. Similarly, for a nonzero $f \in \UC^-(y, x)$, we have a $k$-linear map
	\[
	\psi_f: \UC^+(x, y) \longrightarrow A_x = \UC(x, x) \longrightarrow A_x^0, \quad g \mapsto fg \mapsto \overline{fg},
	\]
	and we say that $f$ is \textit{non-degenerate} if $\psi_f$ is nonzero.
	
	\begin{theorem} \label{main result for decomposition}
		Let $k$ be a field and $\UC$ a generalized $k$-linear Reedy category. Suppose that the following conditions hold for all $x, y \in \Ob(\UC)$:
		\begin{prt}
			\item $\UC^+(x, y)$ is a finite dimensional vector space over $k$;
			
			\item $\UC^+(x, y)$ is a projective right $A_x^0$-module and $\UC^- (y, x)$ is a free left $A_x^0$-module;
			
			\item $\dim_k \UC^+(x, y) = \dim_k \UC^-(y, x)$;
			
			\item every nonzero $g \in \UC^+(x, y)$ is non-degenerate;
            \item $\UC$ is Hom-finite over $k$; that is, $\dim_k\UC(x, y) < \infty$ for all $x, y \in \Ob(\UC)$.
		\end{prt}
		Then one has an equivalence
		\[
		\UC \Mod \simeq \prod_{x \in \Ob(\UC)} A_x^0 \Mod.
		\]
	\end{theorem}
	
	\begin{prf*}
		Similar to the proof of Theorem \ref{main theorem for decomposition}, we show that $\{ \Delta_x \mid x \in \Ob(\UC) \}$ is an orthogonal family of projective generators in $\UC \Mod$. The proof is based on a transfinite induction. By Lemma~\ref{homfinite-support-criterion} below, condition (e) implies that the set of $x$ with $\UC^-(y,x)\ne0$ is finite. This allows us to pass from the vanishing of individual $\Ext^1$ groups to the vanishing against $\mathfrak I_y$.
		
		Firstly, if $x$ and $y$ are distinct objects such that $d(x) = d(y) = 0$, then $\Delta_x = \UC(x, -)$ and $\Delta_y = \UC(y, -)$ are projective $\UC$-modules. Furthermore, these two standard modules are orthogonal by Lemma \ref{Reedy factorization}.
		
		Now we carry out the induction step. Suppose that for
		$$S_{<\alpha} = \{x \in \Ob(\UC) \mid d(x) < \alpha \}$$
		we have proved the following facts:
		\begin{rqm}
			\item $\Hom_{\UC} (\Delta_x, \Delta_y) = 0$ for $x, y \in S_{< \alpha}$ with $x \neq y$;
			\item each $\Delta_x$ is projective.
		\end{rqm}
		We want to show the facts for the subset $S_{\leqslant \alpha} \subseteq \Ob(\UC)$.
		
		\begin{bfhpg*}[Proof of (1)]
			Take distinct $x, y \in S_{\leqslant \alpha}$. Clearly, the conclusion holds if $x, y \in S_{< \alpha}$. Suppose that $d(y) = \alpha$ and $d(x) \leqslant \alpha$. If $d(x) = \alpha$, then $\Hom_{\UC} (\Delta_y, \Delta_x) = 0$ by Theorem \ref{properties of standard modules}, so we may assume that $d(x) < \alpha$. Applying $\Hom_{\UC} (-, \Delta_x)$ to the short exact sequence
			\[
			0 \to \mathfrak{I}_y \to \UC 1_y \to \Delta_y \to 0,
			\]
			one obtains another exact sequence
			\[
			0 \to \Hom_{\UC} (\Delta_y, \Delta_x) \to \Hom_{\UC} (\UC 1_y, \Delta_x) \to \Hom_{\UC} (\mathfrak{I}_y, \Delta_x) \to \Ext_{\UC}^1 (\Delta_y, \Delta_x) \to 0.
			\]
			
			By Remark \ref{a general filtration}, $\mathfrak{I}_y$ has a filtration whose factors are coproducts of the form
			\[
			\UC \otimes_{\UC^-} \UC^-(y, z) \cong \Delta_z \otimes_{A_z^0} \UC^-(y, z)
			\]
			with $d(z) < \alpha$. But $\UC^-(y, z)$ is a free left $A_z^0$-module by the condition (b) and $\Delta_z$ is a projective left $\UC$-module by the induction hypothesis, so every filtration factor is a projective left $\UC$-module. Consequently, we obtain an isomorphism of left $\UC$-modules:
			\[
			\mathfrak{I}_y \cong \bigoplus_{d(z) < \alpha} \Delta_z \otimes_{A_z^0} \UC^-(y, z).
			\]
			Again, the induction hypothesis implies that $\Hom_{\UC} (\Delta_z \otimes_{A_z^0} \UC^-(y, z), \, \Delta_x) = 0$ if $z \neq x$, it follows
			\[
			\Hom_{\UC} (\mathfrak{I}_y, \Delta_x) = \Hom_{\UC} (\Delta_x \otimes_{A_x^0} \UC^-(y, x), \, \Delta_x).
			\]
			Thus the above exact sequence becomes
			\begin{equation*}\label{5.1.1}\tag{$\dag$}
				\begin{gathered}
					0 \to \Hom_{\UC} (\Delta_y, \Delta_x) \to \Hom_{\UC} (\UC 1_y, \Delta_x) \xra{\varphi} \Hom_{\UC} (\Delta_x \otimes_{A_x^0} \UC^-(y, x), \, \Delta_x)\to \Ext_{\UC}^1 (\Delta_y, \Delta_x) \to 0.
				\end{gathered}
			\end{equation*}
			We will show that the middle map $\varphi$ is an isomorphism in the next lemma, so one has
			\[
			\Hom_{\UC} (\Delta_y, \Delta_x) = 0 = \Ext_{\UC}^1 (\Delta_y, \Delta_x).
			\]
		\end{bfhpg*}
		
		\begin{bfhpg*}[Proof of (2)]
We have established $\Ext^1_{\UC}(\Delta_y,\Delta_x)=0$ whenever $d(x)<d(y)=\alpha$. By the Hom-finiteness assumption (e) and Lemma~\ref{homfinite-support-criterion}, only finitely many summands in
\[
\mathfrak I_y\cong\bigoplus_{d(x)<\alpha}
\Delta_x\otimes_{A_x^0}\UC^-(y,x)
\]
are nonzero. Moreover, each $\UC^-(y,x)$ is a free left $A_x^0$-module of finite rank by (a)--(c), so each summand is a finite direct sum of copies of $\Delta_x$. Consequently,
\[
\Ext^1_{\UC}(\Delta_y,\mathfrak I_y)=0.
\]
The short exact sequence $0\to\mathfrak I_y\to\UC1_y\to\Delta_y\to0$ therefore splits. Thus $\Delta_y$ is projective.
\end{bfhpg*}
		
		Now we have proved (1) and (2) for the family $\{ \Delta_x \mid d(x) \leqslant \alpha \}$. A transfinite induction tells us that $\{ \Delta_x \mid x \in \Ob(\UC)\}$ is an orthogonal family of projective generators in $\UC \Mod$. The conclusion then follows.
	\end{prf*}
	
	\begin{lemma} \label{proof of isomorphism}
		Under the assumptions in the previous theorem, the middle map $\varphi$ in the short exact sequence $(\dag)$ is an isomorphism.
	\end{lemma}
	
	\begin{prf*}
By Remark \ref{identification of standard modules}, for each $f \in \UC^+(x, y)$, we denote by $\overline{f}$ the corresponding elements in $\Delta_x(y)$. Via this correspondence,
		\[
		\UC^+(x, y) \cong \Delta_x (y) \cong \Hom_{\UC} (\UC 1_y, \Delta_x)
		\]
		with an explicit isomorphism
		\[
		\phi: \UC^+(x, y) \longrightarrow \Hom_{\UC} (\UC1_y, \Delta_x), \quad g \mapsto (\phi_g: 1_y \mapsto \overline{g})
		\]
		Via this isomorphism, the middle map
		\[
		\varphi: \Hom_{\UC} (\UC 1_y, \Delta_x) \longrightarrow \Hom_{\UC} (\Delta_x \otimes_{A_x^0} \UC^-(y, x), \, \Delta_x)
		\]
		can be described as
		\[
		\phi_g \mapsto [\varphi(\phi_g): \overline{h} \otimes f \mapsto \overline{hfg}, \, \forall \, \overline{h} \in \Delta_x(z), \, f \in \UC^-(y, x)].
		\]
		
		This map is injective. Indeed, if $\phi_g \neq \phi_{g'}$, then $\overline{g} \neq \overline{g'}$, so $g - g' \neq 0$. But $g-g' \in \UC^+(x, y)$ is non-degenerate by the condition (d), so one can find an $f \in \UC^-(y, x)$ such that $\overline{fg} \neq \overline{fg'}$. In particular, for $h = 1_x$, one has
		\[
		\varphi (\phi_g) (\overline{1_x} \otimes f) = \overline{fg} \neq \overline{fg'} = \varphi(\phi_{g'}) (\overline{1_x} \otimes f),
		\]
		so $\varphi (\phi_g) \neq \varphi (\phi_{g'})$; that is, $\varphi$ is injective.
		
		Now we show that $\varphi$ is also surjective. Note that
		$\UC^-(y, x) \cong \bigoplus_{i=1}^n A_x^0$
		for a certain $n \in \mathbb{N}$ since $\UC^-(y, x)$ is a free left $A_x^0$-module of finite rank by the conditions (a) and (b). Consequently, $\Delta_x \otimes_{A_x^0} \UC^-(y, x)$ is isomorphic to the coproduct of $n$ copies of $\Delta_x$, so
		\[
		\Hom_{\UC} (\Delta_x \otimes_{A_x^0} \UC^-(y, x), \, \Delta_x) \cong \bigoplus_{i=1}^n \Hom_{\UC} (\Delta_x, \Delta_x) \cong \bigoplus_{i=1}^n A_x^0,
		\]
		as vector spaces over $k$. By comparing dimensions, we have
		\begin{align*}
		&\dim_k \Hom_{\UC} (\Delta_x \otimes_{A_x^0} \UC^-(y, x), \, \Delta_x)
		= n \dim_k A_x^0 = \dim_k \UC^-(y, x), \ \mathrm{and} \\
		&\dim_k \Hom_{\UC} (\UC 1_y, \Delta_x)
		= \dim_k \UC^+(x, y) = \dim_k \UC^-(y,x),
		\end{align*}
		where the last equality follows from the condition (c). Since $\varphi$ is an injective linear map between two finite dimensional vector spaces with the same dimension, it must be surjective.
	\end{prf*}
	
	Using the right standard module $\Delta^x$ (see Remark \ref{opposite category}), one may obtain the dual version of this theorem.
	
	\begin{theorem} \label{dual main result for decomposition}
		Let $k$ be a field and $\UC$ a generalized $k$-linear Reedy category. Suppose that the following conditions hold for $x, y \in \Ob(\UC)$:
		\begin{prt}
			\item $\UC^+(x, y)$ is a finite dimensional vector space over $k$;
			
			\item $\UC^+ (x, y)$ is a free right $A_x^0$-module and $\UC^-(y, x)$ is a projective left $A_x^0$-module;
			
			\item $\dim_k \UC^+(x, y) = \dim_k \UC^-(y, x)$;
			
			\item every nonzero $f \in \UC^-(y, x)$ is non-degenerate;

            \item $\UC$ is Hom-finite over $k$; that is, $\dim_k\UC(x, y) < \infty$ for all $x, y\in\Ob(\UC)$.
		\end{prt}
		Then one has an equivalence
		\[
		\UC \Mod \simeq \prod_{x \in \Ob(\UC)} A_x^0 \Mod.
		\]
	\end{theorem}

\begin{prf*}
Apply Theorem~\ref{main result for decomposition} to the opposite Reedy category $\UC^{\mathrm{op}}$, with $(\UC^{\mathrm{op}})^+ = (\UC^-)^{\mathrm{op}}$ and $(\UC^{\mathrm{op}})^- = (\UC^+)^{\mathrm{op}}$. The hypotheses above are precisely the corresponding hypotheses for $\UC^{\mathrm{op}}$, so the theorem first gives an equivalence of \emph{right} module categories
\[
\mathrm{Mod}\text{-}\UC \simeq \prod_x\Mod A_x^0.
\]
To pass to left modules, use the Morita theorem for small preadditive categories (equivalently, rings with enough idempotents): a module-category equivalence restricts to an equivalence of their additive idempotent-complete categories of finitely generated projective objects.  The latter determine the associated Morita equivalence; taking the opposite categories yields the Morita equivalence on the other side. In particular, the equivalence of right module categories above induces
\[
\UC \Mod \simeq \prod_x A_x^0 \Mod.
\]
This does not amount to taking the opposite of the entire right-module categories, which in general would not identify them with the corresponding categories of left modules.
\end{prf*}

\begin{lemma}\label{homfinite-support-criterion}
Let $k$ be a field and $\UC$ a generalized $k$-linear Reedy category. Assume that, for all $x, y \in \Ob(\UC)$,
\[
  \dim_k\UC^+(x,y) = \dim_k\UC^-(y,x) < \infty,
\]
and that at least one of the right $A_x^0$-module $\UC^+(x,y)$ and the left $A_x^0$-module $\UC^-(y,x)$ is free. Then the following conditions are equivalent:
\begin{eqc}
\item For every $y\in\Ob(\UC)$, the set $S_y = \{x \in \Ob(\UC) \mid \UC^-(y,x) \ne 0\}$ is finite;
\item $\dim_k\UC(y, y) < \infty$ for every $y \in \Ob(\UC)$;
\item $\UC$ is Hom-finite over $k$.
\end{eqc}
\end{lemma}

\begin{prf*}
By Reedy factorization,
\[
\UC(y,y)\cong
\bigoplus_{x\in\Ob(\UC)}
\UC^+(x,y)\otimes_{A_x^0}\UC^-(y,x).
\]
Under the hypotheses, each tensor factor is finite-dimensional over $k$. Moreover, the equality of dimensions implies that $\UC^+(x,y)$ and $\UC^-(y,x)$ vanish simultaneously. Whenever both are nonzero, their tensor product is nonzero: if $\UC^-(y,x)$ is free, it is a nonzero finite direct sum of copies of $A_x^0$; the other case is symmetric. It follows that $\UC(y,y)$ is finite-dimensional if and only if $S_y$ is finite, proving (i)$\Leftrightarrow$(ii).

For arbitrary objects $u,v$, the decomposition
\[
\UC(u,v)\cong \bigoplus_{x\in S_u} \UC^+(x,v) \otimes_{A_x^0} \UC^-(u,x)
\]
shows that (i) implies (iii), since it is then a finite direct sum of finite-dimensional vector spaces. The implication (iii)$\Rightarrow$(ii) is immediate.
\end{prf*}

\begin{proposition}\label{necessity-of-homfinite}
Let $k$ be a field and $\UC$ a generalized $k$-linear Reedy category. Suppose that $A_x^0$ is finite-dimensional over $k$ for all $x$, and that $\UC^+$ is a projective right $\UC^0$-module. If there is an equivalence of categories
\[
\UC \Mod \simeq \prod_{x \in \Ob(\UC)} A_x^0 \Mod,
\]
then for every $y \in \Ob(\UC)$ the set $S_y = \{ x \in \Ob(\UC) \mid \UC^-(y, x) \ne 0 \}$ is finite.
\end{proposition}

\begin{prf*}
Fix an object $y$ and let $P_y = \UC(y, -)$. Since evaluation at $y$ preserves coproducts, $P_y$ is a compact projective object of $\UC \Mod$. An equivalence preserves compactness and projectivity, so the corresponding object of in $\prod_x A_x^0 \Mod$ is also compact and projective. Thus this object in $\prod_x A_x^0 \Mod$ has only finitely many nonzero components, and each component is a finitely generated projective modules over $A_x^0$. Since all $A_x^0$ are finite-dimensional $k$-algebras, this object has finite length. Thus $P_y$ has finite length as well.

On the other hand, consider the filtration of the left $\UC^-$-module $\UC^-(y,-)$ by degree as in Lemma~\ref{filtration of C^-modules}. The induction functor $\UC\otimes_{\UC^-}-$ is exact by Corollary~\ref{exactness of induction}, and it sends $\UC^-(y,-)$ to $P_y$. Consequently, Remark~\ref{a general filtration} gives a transfinite filtration of $P_y$ whose degree-$\alpha$ subquotient is
\[
\bigoplus_{d(x)=\alpha} \Delta_x \otimes_{A_x^0} \UC^-(y,x).
\]
For each $x$, the summand indexed by $x$ is nonzero whenever $\UC^-(y,x)\ne0$, since its evaluation at $x$ is isomorphic to $\UC^-(y,x)$, using $\Delta_x(x)\cong A_x^0$. If $S_y$ were infinite, either infinitely many degree-subquotients would be nonzero or one subquotient would contain an infinite direct sum of nonzero modules. Both alternatives contradict the finite length of $P_y$. Thus $S_y$ is finite.
\end{prf*}

\begin{corollary}\label{homfinite-iff-decomposition}
Under assumptions (a)--(d) of either Theorem~\ref{main result for decomposition} or Theorem~\ref{dual main result for decomposition}, the following are equivalent:
\begin{eqc}
\item $\UC$ is Hom-finite over $k$;
\item $\UC \Mod \simeq \prod_{x \in \Ob(\UC)} A_x^0 \Mod$.
\end{eqc}
\end{corollary}

\begin{prf*}
The implication (i)$\Rightarrow$(ii) is given by the corresponding decomposition theorem. Conversely, condition (a) ensures that each $A_x^0 = \UC^+(x,x)$ is finite-dimensional, while condition (b) implies that $\UC^+$ is a projective right $\UC^0$-module. The implication (ii)$\Rightarrow$(i) is then a consequence of Proposition~\ref{necessity-of-homfinite} and then Lemma~\ref{homfinite-support-criterion}, noting that the freeness condition in that lemma follows from (b) of the relevant decomposition theorem.
\end{prf*}

\begin{remark} \label{group algebras}
		The condition (b) in the above two theorems can be replaced by other conditions for some special cases. For instance, if each $A_x^0$ is the group algebra $kG_x$ of a finite group $G_x$ whose order is invertible in $k$, then in the proof of Lemma \ref{proof of isomorphism}, by the tensor-hom adjunction one has
		\begin{align*}
			\Hom_{\UC} (\Delta_x \otimes_{kG_x} \UC^-(y, x), \, \Delta_x)
			& \cong \Hom_{kG_x} (\UC^-(y, x), \, \End_{\UC} (\Delta_x))\\
			& \cong \Hom_{kG_x} (\UC^-(y, x), \, kG_x)\\
			& \cong \UC^-(y, x),
		\end{align*}
		as vector spaces over $k$, so the map $\varphi$ is surjective. Note that the assumption that the order of each $G_x$ is invertible in $k$ is essential to guarantee that $\UC^+$ is a projective right $\UC^0$-module and $\UC^-$ is a projective left $\UC^0$-module.
	\end{remark}

	\subsection{Applications}
	
	In this subsection we give some applications of Theorem \ref{main result for decomposition} and its dual version.
	
	Let $\C$ be a small EI category (that is, every endomorphism in $\C$ is an isomorphism) with pullbacks. Without loss of generality we can assume that $\C$ is skeletal since if $\C$ has pullbacks, so does each of its skeletal subcategories. Under this extra assumption it is easy to check that the relation $\preccurlyeq$ on $\Ob(\C)$ defined via setting $x \preccurlyeq y$ if $\C(x, y) \neq \emptyset$ is a partial order on $\Ob(\C)$. Denote by $G_x$ the group of automorphisms on $x$ for $x \in \Ob(\C)$.
	
	\begin{definition}
		A small skeletal EI category $\C$ is called an \emph{artinian} EI category if $(\Ob(\C), \preccurlyeq)$ satisfies the descending chain condition. It is locally finite if $\C(x, y)$ is a finite set for all $x, y \in \Ob(\C)$.
	\end{definition}
	
	From now on, we assume that $\C$ is a locally finite skeletal artinian EI category with pullbacks. Define a degree function $d$ from $\Ob(\C)$ to a sufficiently large ordinal $\lambda$ as follows: for minimal elements in the artinian poset $(\Ob(\C), \preccurlyeq)$, set their degree to be 0. While removing these objects from $\Ob(\C)$, we obtain a new artinian poset, and let minimal elements in the new poset have degree 1. The degree function can be defined via a transfinite process. Let $\C^0$ be the subcategory of $\C$ consisting of all objects and all isomorphisms.
	
	The \textit{category of spans} $\widehat{\C}$ is defined as follows. It has the same objects as $\C$. For $x, y \in \Ob(\C)$, the morphism set is defined to be
	\[
	\widehat{\C} (x, y) = \bigsqcup_{z \in \Ob(\C)} \C (z, y) \times_{G_z} \C^{\op}(x, z).
	\]
	More explicitly, a morphism in $\widehat{\C}(x, y)$ is represented by a diagram
	$\xymatrix{
		x & z \ar[l]_-f \ar[r]^-g & y
	}$
	with $f$ and $g$ morphisms in $\C$, and two pairs $(f, g)$ and $(f', g')$ represent the same morphism if there is an automorphism $\sigma \in G_z$ such that the following diagram commutes:
	\[
	\xymatrix{
		& z \ar@{-->}[d]^-{\sigma} \ar[dr]^-g \ar[dl]_-f \\
		x & z \ar[l]^-{f'} \ar[r]_-{g'} & y
	}
	\]
	The composite of morphisms $(f, g): x \to y$ and $(f', g'): y \to w$ is $(fl, g'h): x \to w$, where the top square is a pullback:
	\[
	\xymatrix{
		v \ar@{-->}[r]^-h \ar@{-->}[d]_-l & u \ar[d]^-{f'} \ar[r]^-{g'} & w\\
		z \ar[r]_-g \ar[d]_-f & y\\
		x
	}
	\]
	For more details, please refer to \cite{Bena} or \cite{CKM}.
	
	It is easy to check that $\widehat{\C}$ is almost a generalized Reedy category in the sense of \cite{BM}. Indeed, by setting $\widehat{\C}^+ = \C$ and $\widehat{\C}^- = \C^{\op}$, we see that all conditions except the last one specified in Example \ref{generalized Reedy categories} hold. Thus, its $k$-linearization $\underline{\widehat{\C}}$ is a generalized $k$-linear Reedy category.
	
	Applying Theorem \ref{main result for decomposition} and Lemma~\ref{homfinite-support-criterion}, we obtain the following result:
	
	\begin{theorem} \label{decomposition for span}
		Let $k$ be a field and $\C$ an artinian EI category satisfying the following conditions:
		\begin{prt}
			\item $\C$ has pullbacks;
			
			\item $\C$ is locally finite;
			
			\item the order of $G_x$ is invertible in $k$ for each $x \in \Ob(\C)$;
			
			\item every morphism in $\C$ is a monomorphism;
            \item for every $y \in \Ob(\C)$, only finitely many isomorphism classes of objects $x$ admit a morphism $x\to y$.
        \end{prt}
		Then one has
		\[
		\underline{\widehat{\C}} \Mod \simeq \prod_{x \in \Ob(\C)} kG_x \Mod.
		\]
	\end{theorem}
	
	\begin{prf*}
		We need to check the conditions (a)-(e) specified in Theorem \ref{main result for decomposition}. Clearly, the condition (a) holds as $\C$ is locally finite, and the condition (b) holds by Remark \ref{group algebras}. Noting that $\widehat{\C}^+ = \C$ and $\widehat{\C}^- = \C^{\op}$. Then the condition (c) automatically holds. Condition (e) follows from the finite-subobject assumption: each span $u \leftarrow z\rightarrow v$ has $z$ among the finitely many objects (up to isomorphism) admitting a morphism to $u$, and the relevant morphism sets in $\C$ are finite, so $\underline{\widehat{\C}}(u,v)$ is finite-dimensional over $k$. Thus it is sufficient to show the condition (d) in Theorem \ref{main result for decomposition}.
		
		Given a nonzero $g \in \underline{\widehat{\C}^+}(x, y) = \underline{\C} (x, y)$, write $g = a_1 g_1 + \ldots + a_n g_n$ with distinct $g_i \in \C(x, y)$ and $0 \neq a_i \in k$. Let $g'_1 \in \widehat{\C}^-(y, x) = \C^{\op}(y, x) = \C(x, y)$ be the morphism corresponding to $g_1$. Then the product $g_1' g_i \in \widehat{\C} (x, x)$ is represented by the pair $(p_i, q_i) $ appeared in the diagram
		\[
		\xymatrix{
			z \ar@{-->}[r]^-{p_i} \ar@{-->}[d]_-{q_i} & x \ar[d]^-{g_1} \ar[r]^-{1_x} & x\\
			x \ar[r]_-{g_i} \ar[d]_{1_x} & y\\
			x
		}
		\]
		When $i = 1$, since $g_1$ is monic, it follows that $z = x$ and $p_i = q_i = 1_x$; see \cite[Proposition 3.12]{Hel}. If $i \neq 1$, then $(p_i, q_i)$ is not equivalent to $(1_x, 1_x)$ since otherwise we have a commutative diagram
		\[
		\xymatrix{
			& x \ar@{-->}[d]^-{\sigma} \ar[dr]^-{1_x} \ar[dl]_-{1_x} \\
			x & x \ar[l]^-{p_i} \ar[r]_-{q_i} & x
		}
		\]
		so $p_i = q_i = \sigma^{-1}$, and hence $g_i = g_1$ from the pullback diagram. Thus, $g_1' g_i$ and $g_1' g_1$ are different morphisms in $\widehat{\C} (x, x)$ if $i \neq 1$. Consequently,
		$g_1' g = a_1 g_1'g_1 + a_2 g_1' g_2 + \ldots + a_n g_1'g_n \in \underline{\widehat{\C}} (x, x)$ is nonzero. Furthermore, its image in the quotient algebra $kG_x$ of $\underline{\widehat{\C}} (x, x)$ is nonzero as well since $g_1'g_1$ is the identity element in $G_x$. Thus, $g$ is non-degenerate.
	\end{prf*}

\begin{remark}\label{span-finiteness-necessary}
The finite-subobject hypothesis in Theorem~\ref{decomposition for span} is also necessary for the asserted product equivalence, under its other assumptions. Indeed, the associated $k$-linear span category satisfies the hypotheses of Proposition~\ref{necessity-of-homfinite} and Lemma~\ref{homfinite-support-criterion}. Hence the product equivalence forces its Hom-finiteness. If an object $y$ admitted infinitely many isomorphism classes of objects $x$ with a morphism $x\to y$, the spans $y\leftarrow x\rightarrow y$ formed using such morphisms would belong to distinct span-isomorphism classes. Their images would give infinitely many linearly independent elements of $\underline{\widehat{\C}}(y,y)$, a contradiction.
\end{remark}

	Theorem \ref{decomposition for span} may apply to many combinatorial categories. We give a few examples, most of which have been studied in literature such as \cite{CEF, LS}. The reader can easily find more.
	
	\begin{example} \label{artinian posets}
		Let $\C$ be an artinian poset such that every two elements in it have the greatest lower bound and every principal lower set $\{x \mid x \leqslant y\}$ is finite. Under these hypotheses all conditions in Theorem \ref{decomposition for span} hold. Therefore, for a field $k$ with any characteristic, one has
		\[
		\underline{\widehat{\C}} \Mod \simeq \prod_{x \in \Ob(\C)} k \Mod.
		\]
	\end{example}
	
	\begin{example} \label{partial linear injections}
		Let $\C$ be the category consisting of objects $\mathbb{F}_q^n$ with $n \in \mathbb{N}$ and $\mathbb{F}_q$-linear injections, and let $\C_{\mathrm{set}}$ be the corresponding category of finite sets and injections. The category $\widehat{\C}$ of spans is equivalent to the category of finite-dimensional spaces over $\mathbb{F}_q$ and partial linear injections. The analogous construction for finite sets and injections gives the
category $\mathrm{FI}\sharp$ of partial injections; it is sometimes described informally as the $q=1$ analogue. One can check that all conditions in Theorem \ref{decomposition for span} hold for $\C$. If $k$ has characteristic zero, the respective vector-space and finite-set conclusions are (see \cite[Corollary 7.7]{LS} and \cite[Theorem 4.1.5]{CEF}):
	\[
		\underline{\widehat{\C}} \Mod \simeq \prod_{n \in \mathbb{N}} k\mathrm{GL}_n(\mathbb{F}_q)\Mod, \quad \underline{\widehat{\C}_{\mathrm{set}}} \Mod \simeq \prod_{n \in \mathbb{N}} kS_n\Mod.
	\]
	\end{example}
	
	\begin{example}
		Let $k$ be a field of characteristic 0. Let $\C$ be the category consisting of objects $\mathbb{F}_q^n$ with $n \in \mathbb{N}$ and $\mathbb{F}_q$-linear surjections. Then $\C^{\op}$ satisfies all conditions in Theorem \ref{decomposition for span}. Consequently, one has
		\[
		\underline{\widehat{\C^{\op}}} \Mod \simeq \prod_{n \in \mathbb{N}} k \mathrm{GL}_n (\mathbb{F}_q) \Mod.
		\]
		Similarly, if we let $\C$ be the category consisting of objects $[n]$ with $n \in \mathbb{N}$ and surjections, then
		\[
		\underline{\widehat{\C^{\op}}} \Mod \simeq \prod_{n \in \mathbb{N}} k S_n \Mod.
		\]
	\end{example}
	
	\begin{remark}
		The conclusions of the above two examples actually hold for any commutative ring $k$ in which $q$ is invertible, as asserted in \cite[Theorem 6.7]{LS}. It follows from the facts that idempotents split and finite direct sums exist in these categories.
	\end{remark}

	\part*{Part II. Abelian model structures}\label{part:model structure}
	\noindent
	In this part we focus on abelian model structures on the representation category $\UC \Mod$ of a generalized $k$-linear Reedy category $\UC$. In particular, using the technique of Grothendieck bifibrations, we construct various abelian model structures on $\UC \Mod$ from families of abelian model structures on $A_x^0 \Mod$'s.
	\begin{setup*}
		Throughout this part, let $\UC$ be a generalized $k$-linear Reedy category with $d: \Ob(\UC) \to \lambda$ its degree function, and recall that $A_x^0=\UC^0(x,x)$ for each $x\in\Ob(\UC)$.
	\end{setup*}

	\section{Preliminaries on model structures and Grothendieck bifibrations} \label{AppA}
	\noindent
	For the convenience of the reader, in this section we give some necessary preliminaries on abelian model structures and Grothendieck bifibrations.
	
	\subsection{Abelian model structures}\label{Abe MS}
	The main content of this subsection includes weak factorization systems, (abelian) model structures, cotorsion pairs as well as various close relations among them.
	
	Let $l: A \to B$ and $r: C \to D$ be morphisms in a category $\calE$. Recall that $l$ has the {\it left lifting property} with respect to $r$ (or $r$ has the {\it right lifting property} with respect to $l$) if for every pair of morphisms $f: A \to C$ and $g: B \to D$ with $rf = gl$, there exists a morphism $t: B\to C$ such that $f = tl$ and $g = rt$. Given a class $\sfC$ of morphisms in $\calE$, let $\sfC^\Box$ denote the class of all morphisms $r$ in $\calE$ having the right lifting property with respect to all morphisms $l$ in $\sfC$. The class $^\Box\sfC$ is defined dually.
	
	Following Bousfield \cite{Bo77}, a pair $(\sfC, \sfD)$ of classes of morphisms in $\calE$ is called a {\it weak factorization system} if $\sfC^\Box = \sfD$ and $\sfC ={^\Box\sfD}$, and every morphism $\alpha$ in $\calE$ can be decomposed as $\alpha = fc$ with $c \in \sfC$ and $f \in \sfD$. Elements in $\sfC$ (resp., in $\sfD$) are called \emph{left maps} (resp., \emph{right maps}) of the weak factorization system.
	
	There are many examples of weak factorization systems. We only recall the following one, which will be used in next sections.
	
	\begin{example}\label{wfs for over/under}
		Given a morphism $f : X \to Y$ in $\calE$, denote by ${_{X}} \backslash^{\calE}_{f} /_{Y}$ the coslice category of $\calE / Y$ under $f$ (or equivalently, the slice category of $X \backslash \calE$ over $f$):
		\begin{prt}
			\item[$\bullet$] objects in ${_{X}} \backslash^{\calE}_{f} /_{Y}$ are morphisms $h$ in $\calE$ such
			that the next diagram commutes:
			\[
			\xymatrix@R=0.5cm@C=0.5cm{
				X
				\ar[rr]^-{h}
				\ar[dr]_-{f}
				&&     \bullet .
				\ar[dl]^-{g}                       \\
				&         Y
			}
			\]
			\item[$\bullet$]
			morphisms in ${_{X}} \backslash^{\calE}_{f} /_{Y}$ from
			$(f \overset{h}\longrightarrow g)$
			to
			$(f \overset{h'}\longrightarrow g')$
			are morphisms $l$ in $\calE$ such that the following diagram commutes:
			\[
			\xymatrix@R=1cm@C=0.75cm{
				\bullet
				\ar@/^2pc/[rrrr]^-{l}
				\ar[drr]_-{g}
				&&          X
				\ar[d]|-{f}
				\ar[rr]^-{h'}
				\ar[ll]_-{h}
				&&      \star .
				\ar[dll]^-{g'}    \\
				&&         Y                 }
			\]
		\end{prt}
		
		It is routine to verify that each weak factorization system $(\sfC, \sfD)$ in $\calE$ induces a weak factorization system $( {_{X}} \backslash^{\sfC}_{f} /_{Y}, \, {_{X}} \backslash^{\sfD}_{f} /_{Y})$ in ${_{X}} \backslash^{\calE}_{f} /_{Y}$, where
		\begin{align*}
		{_{X}} \backslash^{\sfC}_{f} /_{Y}
		&= \{ \ l \in {_{X}} \backslash^{\calE}_{f} /_{Y}
		\, \mid \, \textrm{as a morphism in } \calE, \, l \textrm{ is contained in } \sfC\} \ \mathrm{and} \\
		{_{X}} \backslash^{\sfD}_{f} /_{Y}
		&= \{ \ r \in {_{X}} \backslash^{\calE}_{f} /_{Y}
		\, \mid \, \textrm{as a morphism in } \calE, \, r \textrm{ is contained in } \sfD\}.
		\end{align*}
	\end{example}
	
	The following definition of a model category is adopted from H\"{u}ttemann and R\"{o}ndigs \cite{HR08}. It requires the existence of small limits and colimits, and hence, is a strengthening of Quillen's axioms for a closed model category \cite{Qui67}. It is also slightly more general than the definition given by Hovey \cite{Ho99} in which the factorizations have to be functorial.
	
	\begin{definition}
		Let $\calE$ be a bicomplete category. A {\it model structure} on $\calE$ is a triple $(\sfC, \sfW, \sfF)$ of classes of morphisms satisfying the following conditions:
		\begin{prt}
			\item $(\sfC, \sfW \cap \sfF)$ and $(\sfC \cap \sfW, \sfF)$ are weak factorization systems;
			
			\item $\sfW$ satisfies the 2-out-of-3 property: if two of the three morphisms $\alpha$, $\beta$ and $\beta\alpha$ are contained in $\sfW$, then so is the third one.
		\end{prt}
		
		A morphism in $\sfC$ (resp., $\sfW$, $\sfF$) is said to be a {\it cofibration} (resp., {\it weak equivalence}, {\it fibration}). A morphism in $\sfC \cap \sfW$ (resp., $\sfF \cap \sfW$) is called a {\it trivial cofibration} (resp., {\it trivial fibration}). An object in $\calE$ is said to be \emph{cofibrant} if the map to it from the initial object is a cofibration. Dually, an object in $\calE$ is called \emph{fibrant} if the map from it to the terminal object is a fibration. An object in $\calE$ is \emph{trivial} if the map to it from the initial object is a weak equivalence, or equivalently, the map from it to the terminal object is a weak equivalence.
	\end{definition}
	
	The following definition of an abelian model structure is taken from \cite{Ho02}.
	
	\begin{definition}\label{df of abelian ms}
		Let $\calA$ be a bicomplete abelian category. A model structure on $\calA$ is said to be \emph{abelian} if the following conditions are satisfied:
		\begin{prt}
			\item every cofibration is a monomorphism;
			
			\item a map is a (trivial) fibration if and only if it is an epimorphism with a (trivial) fibrant kernel.
		\end{prt}
	\end{definition}
	
	\begin{remark}\label{df of ams is self-dual}
		By \cite[Proposition 4.2]{Ho02}, a model structure on $\calA$ in which every cofibration is a monomorphism and every fibration is an epimorphism, is abelian if trivial cofibrations are monomorphisms with trivially cofibrant cokernels and cofibrations coincide with monomorphisms with cofibrant cokernels.
	\end{remark}
	
	Weak factorization systems and abelian model structures are closely related to cotorsion pairs, which was first introduced by Salce \cite{S-L} and rediscovered by Enochs and Jenda in \cite{rha}. As an analogue to torsion pairs, cotorsion pairs are defined using the Ext functor instead of the Hom functor. Explicitly, let $\calA$ be an abelian category. Following \cite{rha}, a pair $(\calC,\calD)$ of classes of objects in $\calA$ is called a \emph{cotorsion pair} if $\calC^{\bot} = \calD$ and $^{\bot}\calD =\calC$, where
	\begin{align*}
	\calC^{\bot}
	&= \{\ M \in \calA \, \mid \, \textrm{Ext}_{\calA}^{1}(C, M) = 0
	\textrm{ for all objects } C \in \calC \}  \ \mathrm{and} \\
	^{\bot}\calD
	&= \{\ M \in \calA \, \mid \, \textrm{Ext}_{\calA}^{1}(M, D) = 0
	\textrm{ for all objects } D \in \calD \}.
	\end{align*}
	A cotorsion pair $(\calC,\calD)$ is said to be \emph{complete} if for any object $M$ in $\calA$, there exist short exact sequences
	\begin{center}
		$0 \to D \to C \to M \to 0$ \quad and \quad $0 \to M \to D' \to C' \to 0$
	\end{center}
	in $\calA$ with $D, D' \in \calD$ and $C, C' \in \calC$.
	
	For a class $\calS$ of objects in an abelian category $\calA$, set
\begin{align*}
	\Mon(\calS) & = \{\ f \in \calA \ | \ f \ \text {is a monomorphism}\ \text{with}\
	\coker(f) \in \calS \} \ \mathrm{and} \\
	\Epi(\calS) &= \{\ f \in \calA \ | \ f \ \text {is an epimorphism}\ \text{with}\
	\ker(f) \in \calS \}.
\end{align*}
	The next result, essentially due to Hovey \cite{Ho02} (see also Positselski and \v{S}\v{t}ov\'{\i}\v{c}ek \cite[Theorem 2.4]{Po-St22}), provides a close relation between weak factorization systems and complete cotorsion pairs.
	
	\begin{lemma}\label{bri lem}
		Let $(\calC, \calD)$ be a pair of classes of objects in $\calA$. Then $(\calC,\calD)$ is a complete cotorsion pair if and only if $(\Mon(\calC), \, \Epi(\calD))$ is a weak factorization system in $\calA$.
	\end{lemma}
	
	The following result, now known as \emph{Hovey's correspondence}, is a central result on abelian model structures. Recall that a class of objects in $\calA$ is called {\it thick} if it is closed under direct summands, extensions, kernels of epimorphisms and cokernels of monomorphisms.
	
	\begin{theorem}[Hovey's correspondence]\label{Hovey cor}
		Let $\calA$ be a bicomplete abelian category. Then there exists a bijective correspondence between
		\begin{rqm}
			\item abelian model structures $(\sfC, \sfW, \sfF)$ on $\calA$ and
			
			\item triples $(\calQ, \calW, \calR)$ of classes of objects in $\calA$ such that both $(\calQ, \calW \cap \calR)$ and $(\calQ \cap \calW, \calR)$ are complete cotorsion pairs, and $\calW$ is thick.
		\end{rqm}
		Explicitly, given an abelian model structure $(\sfC, \sfW, \sfF)$ on $\calA$, the corresponding triple of classes of objects in $\calA$ consists of the cofibrant, trivial and fibrant objects. Conversely, given a triple $(\calQ, \calW, \calR)$ satisfying the requirements in $(2)$, the associated abelian model structure is $(\Mon(\calQ), \, \sfW, \, \Epi(\calR))$, where
		\[
		\sfW = \{ w \,|\, w\ \text{can be decomposed as}\ w = fc\ \text{with}\
		c \in \Mon(\calQ \cap \calW)\ \text{and}\ f \in \Epi(\calW \cap \calR) \}.
		\]
	\end{theorem}
	
	By the Hovey's correspondence, an abelian model structure on a bicomplete abelian category $\calA$ can be succinctly represented by a triple of classes of objects in $\calA$ satisfying the conditions in (2). Therefore, one often refers to such a triple as an abelian model structure, and call it a \emph{Hovey triple}.
	
	Recall that a cotorsion pair $(\calC,\calD)$ in $\calA$ is called \textit{resolving} if $\calC$ is closed under kernels of epimorphisms. Dually, one can define \textit{coresolving}  cotorsion pairs. We say that a cotorsion pair $(\calC,\calD)$ is \emph{hereditary} if it is both resolving and coresolving. The following result by Becker \cite{Be14} asserts that the resolving condition and coresolving condition are equivalent when $(\calC,\calD)$ is complete.
	
	\begin{lemma}\label{test hereditary}
		Let $(\calC,\calD)$ be a complete cotorsion pair in an abelian category $\calA$. Then $(\calC,\calD)$ is hereditary if and only if it is resolving, and if and only if it is coresolving.
	\end{lemma}
	
	A Hovey triple  $(\calQ, \calW, \calR)$ in $\calA$ is said to be \emph{hereditary} if both cotorsion pairs $(\calQ, \calW \cap \calR)$ and $(\calQ \cap \calW, \calR)$ are hereditary. The reader may refer to \cite{Gil162} for more on hereditary abelian model structures.

	\subsection{Grothendieck bifibrations}\label{Gro bifib}
	
	The main content of the subsection includes notations, definitions and elementary facts on Grothendieck bifibrations.
	\begin{ipg}
	Let $p : \calT \to \calB$ be a functor. Following Cagne and Melli\`es \cite{CM20}, an object $Y$ in $\calT$ is called \emph{above} an object $V$ in $\calB$ if $p(Y) = V$; similarly, a morphism $f : Y \to Z$ in $\calT$ is called \emph{above} a morphism $u: V \to W$ in $\calB$ if $p(f) = u$. Given an object $V$ in $\calB$, the \emph{fiber} of $V$ with respect to $p$, denoted by $\calT_V$, is the subcategory of $\calT$ with objects $Y$ in $\calT$ such that $p(Y) = V$ and morphisms $f$ in $\calT$ such that $p(f)=1_V$.
	
	Given a morphism $f : Y \to Z$ in $\calT$ which is above a morphism $u : V \to W$ in $\calB$, we say that $f$ is \emph{cartesian} with respect to $p$ if for every morphism $v : U \to V$ in $\calB$ and a morphism $g : X \to Z$ that is above $u \circ v : U \to W$, there is a unique morphism $h : X \to Y$ such that it is above $v$ and $f \circ h = g$, as shown below:
	\[
	\xymatrix@R=0.4cm@C=0.5cm{
		X
		\ar@/^0.85pc/[rrrrrd]^{\forall \, g}
		\ar@.[rrd]^{\exists \, | \,h}
		\ar@{.}[dd]
		\\
		&&        Y
		\ar[rrr]^-{f}
		\ar@{.}[dd]
		&&&       Z
		\ar@{.}[dd]
		\\
		U
		\ar[rrd]^-{\forall \, v}
		\ar@/^0.25pc/@{}[rr]
		&& \ar@/^0.3pc/[rrrd]^(0.15){u \circ v}
		\\
		&&        V
		\ar[rrr]^-{u}
		&&&       W
	}
	\]
	Dually, a morphism $f : Y \to Z$ in $\calT$ that is above a morphism $u : V \to W$ in $\calB$ is called \emph{cocartesian} with respect to $p$ if for every morphism $v : W \to U$ in $\calB$ and a morphism $g : Y \to X$ that is above $v \circ u : V \to U$, there exists a unique morphism $h : Z \to X$ such that it is above $v$ and $h \circ f = g$, as shown below:
	\[
	\xymatrix@R=0.4cm@C=0.5cm{
		&&&&& X
		\ar@{.}[dd]
		\\
		Y
		\ar[rrr]^{f}
		\ar@/^0.85pc/[rrrrru]^{\forall \, g}
		\ar@{.}[dd]
		&&&   Z
		\ar@.[rru]^{\exists \, | \, h}
		\ar@{.}[dd]                          \\
		&&& \ar@/^0.2pc/[rr]
		&&    U                              \\
		V
		\ar[rrr]^{u}
		\ar@/^0.3pc/@{}[rrru]^(0.85){v \circ u}
		&&&   W
		\ar[rru]^{\forall \, v}   }
	\]
	
	It is easy to check that the class of all (co)cartesian morphisms is closed under compositions. We refer the reader to \cite{CM20} for another way to define (co)cartesian morphisms via pullbacks. Recall that a \emph{Grothendieck fibration} is a functor $p : \calT \to \calB$, such that for any morphism $u: V \to W$ in $\calB$ and any object $Z$ in $\calT_{W}$, there is a cartesian morphism $\phi$ with codomain $Z$ that is above $u$; it is called a \emph{cartesian lifting} of $Z$ along $u$. Dually, a \emph{Grothendieck opfibration} is a functor $p : \calT \to \calB$, such that for any morphism $u: V \to W$ in $\calB$ and any object $Y$ in $\calT_V$, there is a cocartesian morphism $\phi$ with domain $Y$ that is above $u$; it is called a \emph{cocartesian lifting} of $Y$ along $u$. A functor $p : \calT \to \calB$ is called a \emph{Grothendieck bifibration} if it is both a Grothendieck fibration and a Grothendieck opfibration.
    \end{ipg}
	
	\begin{definition}\label{df of Grothendieck opfibra}
		Recall from \cite{CM20} that a \emph{cloven Grothendieck bifibration} is a functor $p : \calT \to \calB$ together with
     \begin{prt}
     \item[(G1)] for each morphism $u: V\to W$ in $\calB$ and each object $Z$ in $\calT_W$, an object $u^*(Z)$ in $\calT$ and a cartesian morphism $\rho_{u,Z}: u^*(Z) \to Z$ that is above $u$. Pictorially:
         \[
		\xymatrix@R=0.5cm@C=1.25cm{
			u^*(Z)
			\ar[r]^-{\rho_{u, Z}}
			\ar@{.}[d]
			&  \forall \, Z
			\ar@{.}[d]           \\
			V
			\ar[r]^-{u}
			&  W;}
		\]
     \item[(G2)] for each morphism $u: V\to W$ in $\calB$ and each object $Y$ in $\calT_V$, an object $u_!(Y)$ in $\calT$ and a cocartesian morphism $\lambda_{u, Y}: Y \to u_!(Y)$ that is above $u$. Pictorially:
         \[
		\xymatrix@R=0.5cm@C=1.25cm{
			\forall \, Y
			\ar[r]^-{\lambda_{u, Y}}
			\ar@{.}[d]
			&  u_!(Y)
			\ar@{.}[d]           \\
			V
			\ar[r]^-{u}
			&  W.}
		\]
     \end{prt}
	\end{definition}
	
	In this paper, we always consider \textbf{cloven} Grothendieck bifibrations, and call them simply \emph{Grothendieck bifibrations}. Indeed, if $\calT$ and $\calB$ are small relatively to a universe $\mathbb{U}$ in which the axiom of choice is assumed, then cloven Grothendieck bifibrations coincide with the original notion of Grothendieck bifibrations; see \cite{MoVa20} for details.
	
	\begin{lemma}\label{key factorization}
		Let $p: \calT \to \calB$ be a Grothendieck bifibration and $f: Y \to Z$ a morphism in $\calT$ that is above a morphism $u: V\to W$ in $\calB$.
			Then $f$ can be factored uniquely as the cocartesian morphism $\lambda_{u, Y}: Y \to u_!(Y)$ followed by a morphism $f_\triangleright: u_!(Y)\to Z$ in the fiber $\calT_W$, and $f$ can also be factored uniquely as a morphism $f^{\triangleleft}: Y \to u^*(Z)$ in the fiber $\calT_V$ followed by the cartesian morphism $\rho_{u, Z}: u^*(Z) \to Z$.
		Diagrammatically:
		\[
		\xymatrix@R=1cm@C=1cm{
			Y
			\ar[d]_-{f^{\triangleleft}}
			\ar[rd]|-{f}
			\ar[r]^-{\lambda_{u,Y}}
			& u_!(Y)
			\ar[d]^-{f_{\triangleright}} \\
			u^*(Z)
			\ar[r]_-{\rho_{u,Z}}
			& Z
		}
		\]
	\end{lemma}
	
	\begin{prf*}
		We only prove the first statement as the other one may be proved dually. Consider the following diagram of solid arrows:
		\begin{equation*} \label{factorization 1}
			\tag{\ref{key factorization}.1}
			\begin{gathered}
				\xymatrix@R=0.4cm@C=0.5cm{
					&&&&&
					Z
					\ar@{.}[dd]
					\\
					Y
					\ar[rrr]_{\lambda_{u, Y}}
					\ar@/^0.85pc/[rrrrru]^{f}
					\ar@{.}[dd]
					&&&
					u_!(Y)
					\ar@{.>}[rru]_{f_\triangleright}
					\ar@{.}[dd]
					\\
					&&& \ar@/^0.2pc/[rr] &&
					W
					\\
					V
					\ar[rrr]_{u}
					\ar@/^0.3pc/@{}[rrru]^(0.85){u}
					&&&
					W
					\ar@{=}[rru]}
			\end{gathered}
		\end{equation*}
		By the cocartesian universal property of $\lambda_{u, X}$, there is a unique morphism $f_\triangleright$, which is clearly contained in the fiber $\calT_W$, such that $f_\triangleright\circ\lambda_{u,X}=f$.
	\end{prf*}
	
	\begin{remark}\label{funs u! and u^*}
		Let $p: \calT \to \calB$ be a Grothendieck bifibration and $u: V\to W$ a morphism in $\calB$. Then the assignment $u^*$ appearing in Definition \ref{df of Grothendieck opfibra} forms a functor from $\calT_W$ to $\calT_V$. To specify its action on morphisms, let $k: Z' \to Z$ be a morphism in $\calT_W$, and consider the morphism $k \circ \rho_{u,Z'} : u^*(Z') \to Z$ which is above $u$. By Lemma \ref{key factorization}, the morphism $k \circ \rho_{u, Z'}$ can be factored as $\rho_{u,Z} \circ (k \circ \rho_{u,Z'})^{\triangleleft}$, that is, the following diagram commutes:
		\[
		\xymatrix@R=1cm@C=1cm{
			u^*(Z')
			\ar[d]_-{(k \circ \rho_{u,Z'})^{\triangleleft}} \ar[r]^-{\rho_{u,Z'}}
			& Z'
			\ar[d]^-{k}                                        \\
			u^*(Z)
			\ar[r]^-{\rho_{u,Z}}
			& Z
		}
		\]
		Set $u^*(k)=(k \circ \rho_{u,Z'})^{\triangleleft}$. It is easy to check that this construction is functorial.

 Dually, one may prove that the assignment $u_!$ appearing in Definition \ref{df of Grothendieck opfibra} also forms a functor from $\calT_V$ to $\calT_W$. It sends a morphism $l: Y\to Y'$ in $\calT_V$ to the morphism $(\lambda_{u, Y'} \circ l)_{\triangleright}$ appearing in the following commutative diagram whose existence is guaranteed by Lemma \ref{key factorization}:
		\[
		\xymatrix@R=1cm@C=1cm{
			Y
			\ar[d]_-{l} \ar[r]^-{\lambda_{u,Y}}
			& u_!(Y)
			\ar[d]^-{(\lambda_{u,Y'}\circ l)_\triangleright} \\
			Y'
			\ar[r]^-{\lambda_{u,Y'}}
			& u_!(Y')
		}
		\]
		
		The above functors $u^*$ and $u_!$ are called the \emph{reindexing functors} of $u$. By Lemma \ref{key factorization}, it is routine to check that they form an adjoint pair between $\calT_V$ and $\calT_W$:
		\[
		u_! : \calT_V \rightleftarrows \calT_W : u^*.
		\]
		The reader can refer to \cite{CM20} for more details.
	\end{remark}
	
The following result, taken from \cite[Lemma 2.4]{CM20}, provides a method to construct weak factorization systems in the total category of a Grothendieck bifibration via weak factorization systems in the basis category and all fibers; see also Stanculescu \cite{Stan12}.
	
	\begin{lemma}\label{key lem for con wfs}
		Let $p: \calT \to \calB$ be a Grothendieck bifibration with a weak factorization system $(\sfC_{\calB}, \sfD_{\calB})$ in $\calB$ and a weak factorization system $(\sfC_V, \sfD_V)$ in the fiber $\calT_V$ for each $V \in \Ob(\calB)$. If $u_!(\sfC_V)\subseteq\sfC_W$ for each morphism $u: V\to W$ in $\calB$, then there exists a weak factorization system
		$(\sfC_{\calT}, \sfD_{\calT})$ in $\calT$, which is defined as
		\begin{align*}
		\sfC_{\calT}
		&= \{\ f : Y \to Z \in \calT \, \mid \,
		p(f) \in \sfC_{\calB} \text{ and } f_{\triangleright} \in \sfC_{p(Z)}\} \ \mathrm{and} \\
		\sfD_{\calT}
		&= \{\ f : Y \to Z \in \calT \, \mid \,
		p(f) \in \sfD_{\calB} \text{ and } f^{\triangleleft} \in \sfD_{p(Y)}\ \}.
		\end{align*}
	\end{lemma}
	
	\begin{remark}\label{compatibal-wfs}
		By the proof of \cite[Lemma 2.4]{CM20}, for any morphism $f$ in $\calT$, we have a factorization
		\[
		\xymatrix{
			X
			\ar[rr]^-{f}
			\ar[dr]_-{\tilde{l}}
			&& Y            \\
			&  \bullet
			\ar[ur]_-{\tilde{r}}
		}
		\]
		given by the induced weak factorization system $(\sfC_{\calT}, \sfD_{\calT})$ as well as a factorization
		\[
		\xymatrix{
			p(X)
			\ar[rr]^-{p(f)}
			\ar[dr]_-{l}
			&& p(Y)            \\
			&  p(\bullet)
			\ar[ur]_-{r}
		}
		\]
		of $p(f)$ given by the weak factorization system $(\sfC_{\calB}, \sfD_{\calB})$. Then one has $p(\tilde{l}) = l$ and $p(\tilde{r}) = r$. Moreover,
		the lifting morphism $\tilde{h}$ in the diagram
		\[
		\xymatrix{
			X
			\ar[r]
			\ar[d]
			& V
			\ar[d]   \\
			Y
			\ar[r]
			\ar[ur]|-{\, \tilde{h} \,}
			& W
		}
		\]
		in $\calT$, which is constructed to prove the lifting property of $(\sfC_{\calT}, \sfD_{\calT})$, is above the lifting morphism $h$ in the diagram
		\[
		\xymatrix{
			p(X)
			\ar[r]
			\ar[d]
			& p(V)
			\ar[d]   \\
			p(Y)
			\ar[r]
			\ar[ur]|-{\, h \,}
			& p(W)
		}
		\]
		in $\calB$ given by the lifting property of $(\sfC_{\calB}, \sfD_{\calB})$.
	\end{remark}
	
	\section{Restriction functors and their adjoints}\label{lat ma fun}
	\noindent
	For every non-zero ordinal $\alpha$, define $\UC_{\alpha}$ to be the full subcategory of $\UC$ consisting of objects $x$ with $d(x) < \alpha$. It is easy to see that $\UC_{\alpha}$ inherits a generalized $k$-linear Reedy structure from $\UC$. Consequently, we obtain the chain
	$$\emptyset = \UC_{0}              \subseteq
	\UC_{1}              \subseteq
	\cdots \subseteq
	\UC_{\lambda}         \subseteq
	\UC_{\lambda + 1}      = \UC$$
	of $\UC$ by generalized $k$-linear Reedy subcategories. {\bf Throughout this section}, we let $\beta \leqslant \lambda+1$ be a non-zero ordinal, and consider the generalized $k$-linear Reedy category $\UC_{\beta}$.
	
	\begin{bfhpg}[Adjoint triples]\label{Adjoint tri}
		For each non-zero ordinal $\alpha < \beta$, the fully faithful embedding functor
		$\iota_{\alpha} : \UC_{\alpha} \to \UC_{\beta}$ induces a \textit{restriction functor}
		\[
		\Res_{\alpha} : \UC_{\beta} \Mod \to \UC_{\alpha} \Mod.
		\]
		It admits a left adjoint $\Ind_{\alpha}$ and a right adjoint $\coInd_{\alpha}$ given by the left and right Kan extensions, respectively; see for instance \cite[Section 4]{Mur}. Pictorially:
		\begin{equation*}\label{Ad tri.2}\tag{\ref{lat ma fun}.1}
			\xymatrix@C=5pc{
				\UC_{\beta}\Mod
				\ar[r]|-{ \ \Res_{\alpha} \ }
				&         \UC_{\alpha}\Mod
				\ar@/_1.8pc/[l]_-{\Ind_{\alpha}}
				\ar@/^1.8pc/[l]^-{\coInd_{\alpha}}
			}
		\end{equation*}
		More explicitly, given a left $\UC_{\alpha}$-module $V$ and $x \in \Ob(\UC_{\beta})$, one has
		\begin{align*}
		(\Ind_{\alpha}V)(x)
		&=  1_x \UC_{\beta} \otimes_{\UC_{\alpha}} V =
		\UC_{\beta}(\iota_{\alpha}(-), x) \otimes_{\UC_{\alpha}} V, \ \mathrm{and} \\
		(\coInd_{\alpha}V)(x)
		&= \Hom_{\UC_{\alpha}}(\UC_{\beta} 1_x, V) = \Hom_{\UC_{\alpha}}(\UC_{\beta}(x, \iota_{\alpha}(-)), V).
		\end{align*}
		
		Furthermore, given a left $\UC_{\beta}$-module $Y$ and $x \in \Ob(\UC_{\beta})$ with $d(x) = \alpha$, the counit $l^{\alpha}$ of the adjunction
		$(\Ind_{\alpha}, \Res_{\alpha})$, evaluated at $Y$ and $x$, is the $A^0_x$-homomorphism
		\begin{equation}\label{7.1.1}
			l^{\alpha}_Y(x) : 1_x \UC_{\beta} \otimes_{\UC_{\alpha}} \Res_{\alpha}Y \to Y(x), \quad
			\phi \otimes \xi \mapsto Y(\phi)(\xi),
		\end{equation}
		and the unit $m^{\alpha}$ of the adjunction $(\Res_{\alpha}, \coInd_{\alpha})$, evaluated at $Y$ and $x$, is the $A^0_x$-homomorphism
		\begin{equation}\label{7.1.2}
			m^{\alpha}_Y(x) : Y(x) \to \Hom_{\UC_{\alpha}} (\UC_{\beta} 1_x, \Res_{\alpha}Y), \quad
			\xi \mapsto [f \mapsto Y(f)(\xi)].
		\end{equation}

In literature $1_x \UC_{\beta} \otimes_{\UC_{\alpha}} \Res_{\alpha}Y$ is called the \textit{latching object} of $Y$ at $x$, and $\Hom_{\UC_{\alpha}}(\UC_{\beta} 1_x, \Res_{\alpha}Y)$ is called the \textit{matching object} of $Y$ at $x$.
	\end{bfhpg}

	Now we state the main result of this section, generalizing \cite[Theorem 5.6]{GS}.
	
	\begin{theorem} \label{Cofinality}
		Let $Y$ be a left $\UC_{\beta}$-module and $x \in \Ob(\UC_{\beta})$ with
		$d(x) = \alpha$. Then one has the following natural isomorphisms of left $A_x^0$-modules:
		\begin{prt}
			\item $1_x \UC_{\beta} \otimes_{\UC_{\alpha}} \Res_{\alpha}Y
			\cong
			1_x \UC^+_{\beta} \otimes_{\UC^+_{\alpha}} \Res_{\alpha}Y$;
			
			\item $\Hom_{\UC_{\alpha}}(\UC_{\beta} 1_x, \Res_{\alpha}Y)
			\cong
			\Hom_{\UC^-_{\alpha}}(\UC^-_{\beta} 1_x, \Res_{\alpha}Y).$
		\end{prt}
	\end{theorem}
	
	We need some preparatory work before proving Theorem \ref{Cofinality}. \textbf{In the rest of this section}, for each non-zero ordinal $\alpha < \beta$, we still use the symbol $\mathfrak{I}_{\alpha}$ to denote the following two-sided ideal of $\UC_{\beta}$:
	\[
	\bigoplus_{d(z) < \alpha}
	\UC^+_{\beta}(z, -) \otimes_{A_z^0} \UC^-_{\beta}(\bullet, z).
	\]
	
	\begin{lemma} \label{ind des of I}
		For all $x,y \in \Ob(\UC_{\beta})$ with $d(x) = \alpha$, there is an isomorphism
		\[
		1_x \UC_{\beta}^+ \otimes_{\UC^+_{\alpha}} \UC_{\beta} 1_y
		\cong
		1_x \mathfrak{I}_{\alpha} 1_y.
		\]
	\end{lemma}
	
\begin{prf*}
Let
$N_y=\bigoplus_{d(z)<\alpha}1_z\UC^-_\beta 1_y$,
regarded as a left $\UC^0_\alpha$-module.  For each $t$ with $d(t)<\alpha$, the Reedy
factorization yields an isomorphism
\[
\UC_\beta(y,t)
\cong
\bigoplus_{d(z)<\alpha}
\UC_\alpha^+(z,t)
\otimes_{A_z^0}\UC_\beta^-(y,z).
\]
Indeed, if $\UC_\beta^+(z,t)\ne0$, then
$d(z)\le d(t)<\alpha$, and hence
$\UC_\beta^+(z,t)=\UC_\alpha^+(z,t)$.
These isomorphisms are natural in $t$
with respect to morphisms in $\UC_\alpha^+$.
Consequently, as left $\UC^+_\alpha$-modules, one has
\[
\UC_\beta1_y
 \cong \UC^+_\alpha\otimes_{\UC^0_\alpha}N_y.
\]
Using this isomorphism, we obtain
\begin{align*}
1_x\UC^+_\beta\otimes_{\UC^+_\alpha}\UC_\beta1_y
&\cong
1_x\UC^+_\beta\otimes_{\UC^+_\alpha}
 (\UC^+_\alpha\otimes_{\UC^0_\alpha}N_y)\\
&\cong 1_x\UC^+_\beta\otimes_{\UC^0_\alpha}N_y\\
&\cong \bigoplus_{d(z)<\alpha}
  \UC^+_\beta(z,x)\otimes_{A_z^0}\UC^-_\beta(y,z)\\
&=1_x\mathfrak I_\alpha1_y.
\end{align*}
The resulting isomorphism is induced by composition
in $\UC_\beta$, and is natural in $y$ and
$A_x^0$-linear.
\end{prf*}
	
	The following result will be used frequently in the sequel.
	
	\begin{lemma}\label{T AND T+}
		For each $x \in \Ob(\UC_{\beta})$ with $d(x) = \alpha$, there is an isomorphism
		\[
		1_x \UC_{\beta}^+ \otimes_{\UC^+_{\alpha}} \UC_{\alpha} \cong 1_x \UC_{\beta}
		\]
		of $(A_x^0, \UC_{\alpha})$-bimodules.
	\end{lemma}
	
	\begin{prf*}
		For each $y$ with $d(y)<\alpha$, one has
$\UC_\alpha 1_y=\UC_\beta 1_y$
as left $\UC_\alpha^+$-modules.
		It follows that
		\begin{align*}
			1_x \UC_{\beta}^+ \otimes_{\UC^+_{\alpha}} \UC_{\alpha}
			& \cong 1_x \UC_{\beta}^+ \otimes_{\UC^+_{\alpha}} (\bigoplus_{d(y) < \alpha} \UC_{\alpha} 1_y)       \\
			& =     1_x \UC_{\beta}^+ \otimes_{\UC^+_{\alpha}} (\bigoplus_{d(y) < \alpha} \UC_{\beta} 1_y)  \\ & \cong \bigoplus_{d(y) < \alpha} 1_x \UC_{\beta}^+ \otimes_{\UC^+_{\alpha}} \UC_{\beta} 1_y.
		\end{align*}
		Applying Lemma \ref{ind des of I}, one has
		\[
		\bigoplus_{d(y) < \alpha} 1_x \UC_{\beta}^+ \otimes_{\UC^+_{\alpha}} \UC_{\beta} 1_y  \cong
		\bigoplus_{d(y) < \alpha} 1_x \mathfrak{I}_{\alpha} 1_y.
		\]
		But
		\[
		\bigoplus_{d(y) < \alpha} 1_x \mathfrak{I}_{\alpha} 1_y = \bigoplus_{d(y) < \alpha} \bigoplus_{d(z) < \alpha} \UC^+_{\beta}(z, x)                                \otimes_{A^0_z} \UC^-_{\beta}(y, z) = \bigoplus_{d(y) < \alpha}\UC_{\beta}(y, x),
		\]
		which is precisely the $(A_x^0, \UC_{\alpha})$-bimodule $1_x \UC_{\beta}$. The conclusion follows.
	\end{prf*}
	
	Now, we can give the proof of Theorem \ref{Cofinality}.
	
	\begin{bfhpg}[Proof of Theorem \ref{Cofinality}]
		We prove the isomorphism in (a); the one in (b) may be proved dually. Note that
		\[
			1_x \UC_{\beta}^+ \otimes_{\UC^+_{\alpha}} \Res_{\alpha}Y
			\cong 1_x \UC_{\beta}^+ \otimes_{\UC^+_{\alpha}} (\UC_{\alpha} \otimes_{\UC_{\alpha}} \Res_{\alpha}Y)
			\cong (1_x \UC_{\beta}^+ \otimes_{\UC^+_{\alpha}}\UC_{\alpha}) \otimes_{\UC_{\alpha}} \Res_{\alpha}Y.
		\]
		But $1_x \UC_{\beta} \cong 1_x \UC_{\beta}^+ \otimes_{\UC^+_{\alpha}} \UC_{\alpha}$ as right $\UC_{\alpha}$-modules by Lemma \ref{T AND T+}. Thus, the conclusion follows. \qed
	\end{bfhpg}
	
	The next technical result will be applied in the proofs of Proposition \ref{an des of counit} and Lemma \ref{tilQ = Q cap W}.
	
	\begin{lemma}\label{key for iso}
		Let $Y$ be a left $\UC_{\beta}$-module. Then for each $x \in \Ob(\UC_{\beta})$ with $d(x) = \alpha$, there is an isomorphism
		\[
		1_x\mathfrak{I}_{\alpha} \otimes_{\UC_{\beta}} Y \cong (\Ind_{\alpha} \Res_{\alpha} Y)(x).
		\]
	\end{lemma}
	
	\begin{prf*}
		For all $y \in \Ob(\UC_{\beta})$, by Lemma \ref{ind des of I}, one has
		$1_x \mathfrak{I}_{\alpha} 1_y \cong 1_x \UC_{\beta}^+ \otimes_{\UC^+_{\alpha}} \UC_{\beta} 1_y$.
		On the other hand, by Lemma \ref{T AND T+}, one also has
		\[
		1_x \UC_{\beta} \otimes_{\UC_{\alpha}} \UC_{\beta} 1_y \cong (1_x \UC_{\beta}^+ \otimes_{\UC^+_{\alpha}} \UC_{\alpha}) \otimes_{\UC_{\alpha}} \UC_{\beta} 1_y \cong 1_x \UC_{\beta}^+ \otimes_{\UC^+_{\alpha}} \UC_{\beta} 1_y.
		\]
		Consequently,
		$1_x \mathfrak{I}_{\alpha} 1_y \cong 1_x \UC_{\beta} \otimes_{\UC_{\alpha}} \UC_{\beta} 1_y$,
		and hence,
		$1_x\mathfrak{I}_{\alpha} \cong 1_x\UC_{\beta} \otimes_{\UC_{\alpha}} \UC_{\beta}$
		as right $\UC_{\beta}$-modules. Therefore, one has
		$1_x\mathfrak{I}_{\alpha} \otimes_{\UC_{\beta}} Y \cong (1_x\UC_{\beta} \otimes_{\UC_{\alpha}} \UC_{\beta})\otimes_{\UC_{\beta}} Y \cong 1_x\UC_{\beta} \otimes_{\UC_{\alpha}} \Res_{\alpha} Y=(\Ind_{\alpha} \Res_{\alpha} Y)(x)$.
	\end{prf*}
	
	By Lemma \ref{key for iso} and its dual result, one has the following observation.
	
	\begin{proposition}\label{an des of counit}
		Let $Y$ be a left $\UC_{\beta}$-module. Then for each $x \in \Ob(\UC_{\beta})$ with $d(x) = \alpha$,
		\begin{prt}
			\item the morphism $l^{\alpha}_Y(x) : (\Ind_{\alpha}\Res_{\alpha}Y)(x) \to Y(x)$ in (\ref{7.1.1}) is
			naturally isomorphic to the obvious morphism
			$$1_x\mathfrak{I}_{\alpha} \otimes_{\UC_{\beta}} Y
			\to
			1_x\UC_{\beta} \otimes_{\UC_{\beta}} Y \cong Y(x);$$
			
			\item the morphism $m^{\alpha}_Y(x) : Y(x) \to (\coInd_{\alpha}\Res_{\alpha}Y)(x)$ in (\ref{7.1.2}) is
			naturally isomorphic to the obvious morphism
			\[
			Y(x) \cong \Hom_{\UC_{\beta}} (\UC_{\beta} 1_x, Y)
			\to
			\Hom_{\UC_{\beta}} (\mathfrak{I}_{\alpha} 1_x, Y).
			\]
		\end{prt}
	\end{proposition}
	
	\begin{prf*}
		We only show the first statement as the second one may be proved dually. Since
		$1_x\mathfrak{I}_{\alpha} \otimes_{\UC_{\beta}} Y \cong (\Ind_{\alpha} \Res_{\alpha} Y)(x)$
		by Lemma \ref{key for iso}, the conclusion follows readily as one can obtain the the counit morphism
		$l^{\alpha}_Y(x): (\Ind_{\alpha}\Res_{\alpha}Y)(x) \to Y(x)$ via tensoring the inclusion
		$1_x\mathfrak{I}_{\alpha} \to 1_x\UC_{\beta}$ by $Y$.
	\end{prf*}
	
	As an immediate consequence of Proposition \ref{an des of counit}, one has:
	
	\begin{corollary} \label{cokernel of latching}
		Let $Y$ be a left $\UC_{\beta}$-module and $x \in \Ob(\UC_{\beta})$ with
		$d(x) = \alpha$. Then the following hold.
		\begin{prt}
			\item The cokernel of $l^{\alpha}_Y(x)$ is isomorphic to
			$\Delta^x \otimes_{\UC_{\beta}} Y$, and $l^{\alpha}_Y(x)$ is a monomorphism
			if and only if $\Tor_1^{\UC_{\beta}}(\Delta^x, Y) = 0$;
			
			\item The kernel of $m^{\alpha}_Y(x)$ is isomorphic to
			$\Hom_{\UC_{\beta}}(\Delta_x, Y)$, and $m^{\alpha}_Y(x)$ is an epimorphism
			if and only if $\Ext^1_{\UC_{\beta}}(\Delta_x, Y) = 0$.
		\end{prt}
	\end{corollary}
	
	\begin{prf*}
		As before, we only prove the first statement. Applying $- \otimes_{\UC_{\beta}} Y$ to the short exact sequence
		\[
		0 \to 1_x \mathfrak{I}_{\alpha} \to 1_x \UC_{\beta} \to \Delta^x \to 0
		\]
		of right $\UC_{\beta}$-modules, we deduce the conclusion from Proposition \ref{an des of counit}.
	\end{prf*}
	
	\section{The bifibrational nature of the restriction functor}
	\label{bifibrational nat}
	\noindent
	In this section, we consider the embedding $\iota_{\alpha} : \UC_{\alpha} \to \UC_{\alpha + 1}$ for each non-zero ordinal $\alpha \leqslant \lambda$, and show that the restriction functor $\Res_{\alpha} : \UC_{\alpha + 1}\Mod \to \UC_{\alpha}\Mod$ is a Grothendieck bifibration.
	
	\begin{bfhpg}[A natural transformation]\label{The natural tran}
		Note that $\iota_{\alpha} : \UC_{\alpha} \to \UC_{\alpha + 1}$ is fully faithful. It follows that both the functors $\Ind_{\alpha}$ and $\coInd_{\alpha}$ in (\ref{Ad tri.2}) are also fully faithful; see \cite[Proposition 4.23]{Kel05}. Consequently, both the unit $\eta^{\alpha}$ of $(\Ind_{\alpha}, \Res_{\alpha})$ and the counit $\epsilon^{\alpha}$ of $(\Res_{\alpha}, \coInd_{\alpha})$ are isomorphisms. Thus, there is a natural transformation
		\[
		\tau^{\alpha} : \Ind_{\alpha} \to \coInd_{\alpha},
		\]
		which can be described as the map
		\[
		\Ind_{\alpha}
		\xra{m^{\alpha} \circ \Ind_{\alpha}}
		\coInd_{\alpha} \circ \Res_{\alpha} \circ \Ind_{\alpha}
		\xra{\coInd_{\alpha} \circ (\eta^{\alpha})^{-1}}
		\coInd_{\alpha},
		\]
		or equivalently, the map
		\[
		\Ind_{\alpha}
		\xra{\Ind_{\alpha} \circ (\epsilon^{\alpha})^{-1}}
		\Ind_{\alpha} \circ \Res_{\alpha} \circ \coInd_{\alpha}
		\xra{l^{\alpha} \circ \coInd_{\alpha}}
		\coInd_{\alpha},
		\]
		where $m^{\alpha}$ is the unit of $(\Res_{\alpha}, \coInd_{\alpha})$ and $l^{\alpha}$ is the counit of $(\Ind_{\alpha}, \Res_{\alpha})$; see \cite[Lemmas 2.28]{GS}. Moreover, one has a commutative diagram
		\begin{equation*}\label{Fac of NT}\tag{\ref{bifibrational nat}.1.1}
			\begin{gathered}
				\xymatrix@R=0.75cm@C=0.15cm{
					\Ind_{\alpha}\Res_{\alpha}
					\ar[rr]^-{\tau^{\alpha} \ast 1_{\Res_{\alpha}}}
					\ar[dr]_-{l^{\alpha}}
					&&       \coInd_{\alpha}\Res_{\alpha}                       \\
					&        1_{\UC_{\alpha + 1}\Mod}
					\ar[ur]_-{m^{\alpha}}
				}
			\end{gathered}
		\end{equation*}
		of natural transformations, where ``$\ast$" denotes the Godment product; see \cite[Lemma 2.29]{GS} for details.
	\end{bfhpg}
	
	In the following, we describe the fiber ${\UC_{\alpha +1}\Mod}_V$ of a left $\UC_{\alpha}$-module $V$ with respect to the restriction functor $\Res_{\alpha}: \UC_{\alpha + 1}\Mod \to \UC_{\alpha}\Mod$.
	
	\begin{ipg}\label{1-1 for exten}
		Let $Y$ be in ${\UC_{\alpha +1}\Mod}_V$. By (\ref{Fac of NT}), one gets a factorization
		\[
		\Ind_{\alpha}V = \Ind_{\alpha}\Res_{\alpha}Y
		\xra{l^{\alpha}_Y}
		Y
		\xra{m^{\alpha}_Y}
		\coInd_{\alpha}\Res_{\alpha}Y = \coInd_{\alpha}V
		\]
		of $\tau^{\alpha}_V$. In particular, for all objects $x$ in $\UC_{\alpha +1}$ with
		$d(x) = \alpha$, the above factorization restricts to factorizations
		\[
		\xymatrix@C=1cm{
			(\Ind_{\alpha}V)(x)
			\ar[r]^-{l^{\alpha}_Y(x)}
			&           Y(x)
			\ar[r]^-{m^{\alpha}_Y(x)}
			&    (\coInd_{\alpha}V)(x)
		}
		\]
		of $A^0_x$-homomorphism factorizations of those $\tau^{\alpha}_V(x)$.
		
		Conversely, given a family
		$\{(\Ind_{\alpha}V)(x)
		\xra{\varphi(x)}
		M_x
		\xra{\psi(x)}
		(\coInd_{\alpha}V)(x)
		\}_{d(x) = \alpha}$
		of $A^0_x$-homomorphism factorizations of those $\tau^{\alpha}_V(x)$, one may define a left $\UC_{\alpha +1}$-module $Y$ as follows:
		
		For any object $x$ in $\UC_{\alpha +1}$, set
		\begin{equation*}\label{dfn of U.1}
			\tag{\ref{1-1 for exten}.1}
			Y(x) =
			\begin{cases}
				V(x)         & \text{if } d(x) < \alpha, \text{ i.e. } x \in \Ob(\UC_{\alpha});      \\
				M_x          & \text{if } d(x) = \alpha.
			\end{cases}
		\end{equation*}
		In the next three steps we define the values of $Y$ on morphisms in $\UC_{\alpha +1}$.
		
		\begin{bfhpg*}[Step 1] Given a $f^+ \in \UC^+_{\alpha +1}(x, y)$, we have the following cases:
			\begin{rqm}
				\item
				If $d(x) = \alpha$, then $d(y) = \alpha$ as well. Hence, $\UC^+_{\alpha +1}(x, y) = 0$ or $x = y$, and correspondingly $f^+ = 0$ or $f^+ \in A^0_x$. If $f^+ = 0$, we set $Y(f^+) = 0$; otherwise, we set $Y(f^+)$ to be $M_x(f^+)$ by noting that $M_x$ is a left $A^0_x$-module.
				
				\item
				If $d(x) < \alpha$ and $d(y) = \alpha$, then $Y(x) = V(x)$. Consider the following diagram of solid arrows:
				\[
				\xymatrix@R=1cm@C=1.5cm{
					(\Ind_{\alpha}V)(x) = (\Ind_{\alpha}\Res_{\alpha}V)(x)
					\ar[r]^-{l^{\alpha}_{V}(x)}
					\ar[d]_-{(\Ind_{\alpha}V)(f^+)}
					&  V(x) = Y(x)
					\ar@{.>}[d]^-{Y(f^+)}                                  \\
					(\Ind_{\alpha}V)(y)
					\ar[r]^-{\varphi(y)}
					& M_y = Y(y)
				}
				\]
				Since $d(x) < \alpha$, the top morphism $l^{\alpha}_{V}(x)$ is an isomorphism. We define $Y(f^+)$ to be the morphism such that the above diagram commutes, that is,
				\[
				Y(f^+) = \varphi(y) \circ (\Ind_{\alpha}V)(f^+) \circ (l^{\alpha}_{V}(x))^{-1}.
				\]
				\item
				If $d(x) < \alpha$ and $d(y) < \alpha$, then $f^+ \in \UC^+_{\alpha}(x,y)$, and $Y(x) = V(x)$ and $Y(y) = V(y)$. In this case, set $Y(f^+)$ to be $V(f^+)$.
			\end{rqm}
		\end{bfhpg*}
		
		\begin{bfhpg*}[Step 2] Given a $f^- \in \UC^-_{\alpha +1}(x, y)$, we have the following cases:
			\begin{rqm}
				\item[$(1^\prime)$]
				If $d(y) = \alpha$, then $d(x) = \alpha$ as well. Hence, $\UC^-_{\alpha +1}(x, y) = 0$ or $x = y$, and correspondingly $f^- = 0$ or $f^- \in A^0_x$. If $f^- = 0$, then set $Y(f^-) = 0$; otherwise, $f^- \in A^0_x$, so we can define $Y(f^-)$ to be $M_x(f^-)$.
				
				\item[$(2^\prime)$]
				If $d(y) < \alpha$ and $d(x) = \alpha$, then $Y(y) = V(y)$. Since $m^{\alpha}_V(y)$ in the following diagram is an isomorphism
				\[
				\xymatrix@R=1cm@C=1.5cm{
					Y(x) = M_x
					\ar[r]^-{\psi(x)}
					\ar@{.>}[d]_-{Y(f^-)}
					&  (\coInd_{\alpha}V)(x)
					\ar[d]^-{(\coInd_{\alpha}V)(f^-)}
					\\
					Y(y) = V(y)
					\ar[r]^-{m^{\alpha}_V(y)}
					& (\coInd_{\alpha}\Res_{\alpha}V)(y) = (\coInd_{\alpha}V)(y)
				}
				\]
				we set
				\[Y(f^-) = (m^{\alpha}_V(y))^{-1} \circ (\coInd_{\alpha}V)(f^-) \circ \psi(x).
				\]
				
				\item[$(3^\prime)$]
				If $d(y) < \alpha$ and $d(x) < \alpha$, then $f^- \in \UC^-_{\alpha}(x,y)$, and $Y(y) = V(y)$ and $Y(x) = V(x)$. In this case, we set $Y(f^-)$ to be $V(f^-)$.
			\end{rqm}
		\end{bfhpg*}
		
		\begin{bfhpg*}[Step 3] Given an arbitrary $f \in \UC_{\alpha + 1} (x, y)$, we have a Reedy decomposition
			$f = f_1^+ f_1^- + \ldots + f_n^+ f_n^-$,
			so we define
			$Y(f) = Y(f_1^+) \circ Y(f_1^-) + \ldots + Y(f_n^+) \circ Y(f_n^-)$.
		\end{bfhpg*}
	\end{ipg}
	
	\begin{lemma}
		The above construction gives a left $\UC_{\alpha + 1}$-module $Y$ in ${\UC_{\alpha +1}\Mod}_V$.
	\end{lemma}
	
	\begin{prf*}
		It suffices to prove that the construction of $Y(f)$ in Step 3 is independent of the choice of Reedy factorizations of $f$. But as we explained in Remark \ref{uniqueness of factorization}, Reedy factorizations are unique up to tensor product over $\UC^0$. Thus, we only need to check the construction of $Y(f)$ respects the balanced relation $f^+ \otimes_k \phi f^- - f^+ \phi \otimes_k f^-$, where
		$f^+ \in \UC_{\alpha+1}^+(z, y), \phi \in \UC_{\alpha+1}^0(z, z)$ and $f^- \in \UC^-_{\alpha+1}(x, z)$.
		Explicitly, given the following diagram
		\[
		\xymatrix@R=1cm@C=1.5cm{
			x
			\ar[r]^-{f^-}
			&  z
			\ar[r]^-{f^+}
			\ar@(ur,ul)[]_-{\phi}
			&  y
		}
		\]
		of morphisms in $\UC_{\alpha + 1}$, we need to show that the diagram
		\[
		\xymatrix@R=1cm@C=2cm{
			Y(x)
			\ar[r]^-{Y(\phi f^-)}
			\ar[d]_-{Y(f^-)}
			&   Y(z)
			\ar[d]^-{Y(f^+)}                      \\
			Y(z)
			\ar[r]^-{Y(f^+ \phi)}
			&   Y(y)
		}
		\]
		commutes; that is, $Y(f^+) \circ Y(\phi f^-) = Y(f^+ \phi) \circ Y(f^-)$.
		
		If $d(z) = \alpha$, then $d(x) = \alpha = d(y)$. We may suppose that $x = z = y$ because for the other cases one has $f^+ = 0$ or $f^- = 0$, and the desired equality holds clearly. Under this assumption, all of $f^-$, $\phi$ and $f^+$ are contained in $A^0_x$ and $Y(x) = M_x$ is a left $A^0_x$-module. Thus, the desired equality clearly follows from the axiom of modules.
		
		Now, suppose that $d(z) < \alpha$. We establish the conclusion by showing $Y(\phi f^-) = Y(\phi) \circ Y(f^-)$ and $Y(f^+ \phi) = Y(f^+) \circ Y(\phi)$.
		
		If $d(x) = \alpha$, then one has
		\begin{align*}
			\quad Y(\phi f^-)
			& = (m^{\alpha}_V(z))^{-1} \circ (\coInd_{\alpha}V)(\phi f^-) \circ \psi(x)
			& \text{ by } (2^{\prime})                                                    \\
			& = (m^{\alpha}_V(z))^{-1} \circ (\coInd_{\alpha}V)(\phi) \circ
			(\coInd_{\alpha}V)(f^-) \circ \psi(x)
			& \text{ as }  \phi, f^- \in \UC_{\alpha + 1}                                 \\
			& = V(\phi) \circ (m^{\alpha}_V(z))^{-1} \circ  (\coInd_{\alpha}V)(f^-) \circ \psi(x)
			& \text{as}\ m^{\alpha}_V\ \text{is a natural transformation}
			\\
			& = V(\phi) \circ Y(f^-)
			& \text{ by } (2^{\prime}) \text{ again}                                      \\
			& = Y(\phi) \circ Y(f^-)
			& \text{ by } (3^{\prime})
		\end{align*}
		Otherwise, if $d(x) < \alpha$, then all of $\phi$, $f^-$ and $\phi f^-$ are contained in $\UC^-_{\alpha}$, and so by $(3^{\prime})$, one has
		\[
		Y(\phi f^-) = V(\phi f^-) = V(\phi) \circ V(f^-) = Y(\phi) \circ Y(f^-).
		\]
		Consequently, we always have $Y(\phi f^-) = Y(\phi) \circ Y(f^-)$.
		
		Similarly, if $d(y) = \alpha$, then one has
		\begin{align*}
			Y(f^+ \phi)
			& = \varphi(y) \circ (\Ind_{\alpha}V)(f^+ \phi) \circ (l^{\alpha}_{V}(z))^{-1}
			& \text{ by } (2)                                                             \\
			& = \varphi(y) \circ (\Ind_{\alpha}V)(f^+) \circ
			(\Ind_{\alpha}V)(\phi) \circ (l^{\alpha}_{V}(z))^{-1}
			& \text{ as }  f^+, \phi \in \UC_{\alpha + 1}                                 \\
			& = \varphi(y) \circ (\Ind_{\alpha}V)(f^+) \circ
			(l^{\alpha}_{V}(z))^{-1} \circ V(\phi)
			& \text{as}\ l^{\alpha}_V\ \text{is a natural transformation}
			\\
			& = Y(f^+) \circ V(\phi)
			& \text{ by } (2) \text{ again}                                               \\
			& = Y(f^+) \circ Y(\phi)
			& \text{ by } (3)
		\end{align*}
		Otherwise, if $d(y) < \alpha$, then all of $f^+$, $\phi$, and $f^+ \phi$ are contained in $\UC^+_{\alpha}$, and so by $(3)$, one has
		\[
		Y(f^+ \phi) = V(f^+ \phi) = V(f^+) \circ V(\phi) = Y(f^+) \circ Y(\phi).
		\]
		Consequently, we always have $Y(f^+ \phi) = Y(f^+) \circ Y(\phi)$.

It remains to show that $Y$ preserves compositions.
By Steps~1 and~2, it suffices to verify
\[
Y(g^-f^+)=Y(g^-)Y(f^+)
\]
for $f^+\in\UC_{\alpha+1}^+(a,b)$ and
$g^-\in\UC_{\alpha+1}^-(b,c)$.
If $d(b)<\alpha$, the equality follows from
the functoriality of $V$. If $d(b)=\alpha$
and either $d(a)=\alpha$ or $d(c)=\alpha$,
it follows from the $A_b^0$-linearity of
$\varphi(b)$ and $\psi(b)$.

Now suppose that $d(a),d(c)<\alpha=d(b)$.
Since $\psi(b)\varphi(b)=\tau_V^\alpha(b)$,
one has
\[
\begin{aligned}
Y(g^-)Y(f^+)
&=(m_V^\alpha(c))^{-1}
  (\coInd_\alpha V)(g^-)
  \tau_V^\alpha(b)
  (\Ind_\alpha V)(f^+)
  (l_V^\alpha(a))^{-1}
=V(g^-f^+)
=Y(g^-f^+),
\end{aligned}
\]
where the second equality follows from
the naturality of $\tau_V^\alpha$ and
the adjunction isomorphisms on
$\UC_\alpha$.
Thus $Y$ preserves compositions and hence
defines a left $\UC_{\alpha+1}$-module.

	\end{prf*}

	One can easily deduce the following result from \ref{1-1 for exten} and the above lemma.
	
	\begin{theorem}\label{key 1-1 core}
		Let $V$ be a left $\UC_{\alpha}$-module. Then there exists a one-to-one correspondence between modules in the fiber ${\UC_{\alpha +1}\Mod}_V$ and families
		\[
		\{
		(\Ind_{\alpha}V)(x) \longrightarrow \bullet \longrightarrow (\coInd_{\alpha}V)(x)
		\}_{d(x) = \alpha}
		\]
		of left $A^0_x$-homomorphism factorizations of $\tau^{\alpha}_V(x)$.
	\end{theorem}
	
	\begin{remark}\label{key 1-1 core for mor}
		Let $V$ and $W$ be left $\UC_{\alpha}$-modules, and $Y\in\UC_{\alpha+1}\Mod_{V}$ and $Y'\in\UC_{\alpha+1}\Mod_{W}$. By Theorem \ref{key 1-1 core}, the commutative diagram (\ref{Fac of NT}) of natural transformations, and an argument similar to that in \cite[Observation 3.11]{RV14}, one concludes that an extension of a morphism $u: V\to W$ in $\UC_{\alpha}\Mod$ to a morphism $f : Y \to Y'$ in $\UC_{\alpha +1}\Mod$ that is above $u$ uniquely corresponds to a family
		$\{
		Y(x)
		\xra{\delta(x)}
		Y'(x)
		\}_{d(x) = \alpha}$
		of morphisms of left $A^0_x$-modules such that the following diagram commutes:
		\[
		\xymatrix@R=1cm@C=0.5cm{
			(\Ind_{\alpha}V)(x)
			\ar[rr]
			\ar[d]_-{(\Ind_{\alpha}u)(x)}
			&&                    Y(x)
			\ar[rr]
			\ar[d]|-{\delta(x)}
			&&        (\coInd_{\alpha}V)(x)
			\ar[d]^-{(\coInd_{\alpha}u)(x)}               \\
			(\Ind_{\alpha}W)(x)
			\ar[rr]
			&&                   Y'(x)
			\ar[rr]
			&&        (\coInd_{\alpha}W)(x).
		}
		\]
We mention that for each $x\in \Ob(\UC_{\alpha+1})$, the morphism $f(x): Y(x)\to Y'(x)$ is as follows:
\begin{equation*}
				f(x) =
				\begin{cases}
					u(x)              & \text{if } d(x) < \alpha, \text{ i.e. } x \in \Ob(\UC_{\alpha});      \\
					\delta(x)         & \text{if } d(x) = \alpha.
				\end{cases}
			\end{equation*}
	\end{remark}
	
	\begin{ipg}
		Let $V$ be a left $\UC_{\alpha}$-module. We give a description for the category ${\UC_{\alpha +1}\Mod}_V$. For each $x \in \Ob(\UC_{\alpha +1})$ with $d(x) = \alpha$, recall from Example \ref{wfs for over/under} the category
		${_{(\Ind_{\alpha}V)(x)}}
		\backslash{\substack{A^0_x\Mod \\ \tau^{\alpha}_V(x)}}
		/_{(\coInd_{\alpha}V)(x)}$. Define a rule
		\[
		\scrR: {\UC_{\alpha +1}\Mod}_V
		\longrightarrow
		\prod_{d(x) = \alpha}
		{_{(\Ind_{\alpha}V)(x)}}
		\backslash{\substack{A^0_x\Mod \\ \tau^{\alpha}_V(x)}}
		/_{(\coInd_{\alpha}V)(x)}
		\]
		as follows:
		\begin{prt}
			\item[$\bullet$]
			For any $Y \in {\UC_{\alpha +1}\Mod}_V$, set $\scrR(Y)$ to be the family
			\[
			\quad  \quad  \quad
			\{
			(\Ind_{\alpha}V)(x)
			\xra{l^{\alpha}_Y(x)}
			Y(x)
			\xra{m^{\alpha}_Y(x)}
			(\coInd_{\alpha}V)(x)
			\}_{d(x) = \alpha}
			\]
			of left $A^0_x$-homomorphism factorizations of those $\tau^{\alpha}_V(x)$ given in \ref{1-1 for exten}.
			
			\item[$\bullet$]
			Let $f : Y \to Y'$ be a morphism in ${\UC_{\alpha +1}\Mod}_V$. Then for each $x \in \Ob(\UC_{\alpha +1})$ with $d(x) = \alpha$, $f(x)$ makes the diagram
			\[
			\xymatrix@R=0.5cm@C=0.75cm{
				(\Ind_{\alpha}\Res_{\alpha}Y)(x)
				\ar[rr]^-{l^{\alpha}_Y(x)}
				\ar@{=}[d]
				&&                Y(x)
				\ar[rr]^-{m^{\alpha}_Y(x)}
				\ar[dd]|-{f(x)}
				&& (\coInd_{\alpha}\Res_{\alpha}Y)(x)
				\ar@{=}[d]                                              \\
				(\Ind_{\alpha}V)(x)
				\ar@{=}[d]
				&&&&      (\coInd_{\alpha}V)(x)
				\ar@{=}[d]                                                              \\
				(\Ind_{\alpha}\Res_{\alpha}Y')(x)
				\ar[rr]^-{l^{\alpha}_{Y'}(x)}
				&&               Y'(x)
				\ar[rr]^-{m^{\alpha}_{Y'}(x)}
				&&    (\coInd_{\alpha}\Res_{\alpha}Y)(x)
			}
			\]
			commutes; see Remark \ref{key 1-1 core for mor}. Therefore, $f(x)$ is a morphism in the category
			$            {_{(\Ind_{\alpha}V)(x)}}
			\backslash{\substack{A^0_x\Mod \\ \tau^{\alpha}_V(x)}}
			/_{(\coInd_{\alpha}V)(x)}$.
			Set $\scrR(f)$ to be the family $\{ f(x) \}_{d(x) = \alpha}$.
		\end{prt}
		It is routine to check that $\scrR$ is a functor.
	\end{ipg}

The next result can be easily deduced from Theorem \ref{key 1-1 core}.

\begin{corollary}\label{the iso of fiber}
Let $V$ be a left $\UC_\alpha$-module. Then $\scrR: {\UC_{\alpha +1}\Mod}_V
		\to
		\prod_{d(x) = \alpha}
		{_{(\Ind_{\alpha}V)(x)}}
		\backslash{\substack{A^0_x\Mod \\ \tau^{\alpha}_V(x)}}
		/_{(\coInd_{\alpha}V)(x)}$ is an isomorphism of categories.
\end{corollary}
	
	Now we state the main result in this section, which exhibits the bifibrational nature of the restriction functor $\Res_{\alpha}$.
	
	\begin{theorem}\label{res is Gro bifib}
		The restriction functor
		$\Res_{\alpha} : \UC_{\alpha +1}\Mod \to \UC_{\alpha}\Mod$
		with respect to the embedding $\iota_{\alpha} : \UC_{\alpha} \to \UC_{\alpha + 1}$
		is a Grothendieck bifibration.
	\end{theorem}
	
	\begin{prf*}
		We first prove (G2) in Definition \ref{df of Grothendieck opfibra}. To this end, given a morphism $u : V \to W$ in $\UC_{\alpha}\Mod$ and a module $Y$ in the fiber ${\UC_{\alpha+1}\Mod}_{V}$, we have to construct a $\UC_{\alpha +1}$-module $u_!(Y)$ in the fiber ${\UC_{\alpha+1}\Mod}_{W}$ and a cocartesian morphism $\lambda_{u,Y} : Y \to u_!(Y)$ that is above $u$.
		
		\begin{bfhpg*}[Step 1: Constructions of $u_!(Y)$ and $\lambda_{u,Y}$] By Theorem \ref{key 1-1 core}, the module $Y$ gives rise to a family
			\[
			\{               (\Ind_{\alpha}V)(x)\xra{l^{\alpha}_Y(x)}
			Y(x)
			\xra{m^{\alpha}_Y(x)}
			(\coInd_{\alpha}V)(x)
			\}_{d(x) = \alpha}
			\]
			of left $A^0_x$-homomorphism factorizations of those $\tau^{\alpha}_V(x)$ (see
			\ref{1-1 for exten} for details.) For each $x \in \Ob(\UC_{\alpha +1})$ with $d(x) = \alpha$, the pushout
			\begin{equation*}\label{pushout.1}\tag{\ref{res is Gro bifib}.1}
				\begin{gathered}
					\xymatrix@R=1cm@C=0.35cm{
						(\Ind_{\alpha}V)(x)
						\ar@{}[rd]^(0.75)>>{\lrcorner}
						\ar[r]^-{(\Ind_{\alpha}u)(x) }
						\ar[d]_-{l^{\alpha}_Y(x)}
						&  (\Ind_{\alpha}W)(x)
						\ar[d]^-{\varphi(x)}
						\\
						Y(x)
						\ar[r]^-{\theta(x)}
						&  (\Ind_{\alpha}W)(x) \sqcup_{(\Ind_{\alpha}V)(x)} Y(x) }
				\end{gathered}
			\end{equation*}
			in $A^0_x\Mod$ fits into the following diagram of solid arrows:
			\[
			\xymatrix@R=1cm@C=0.5cm{
				(\Ind_{\alpha}V)(x)
				\ar@/_7.8pc/[dd]_-{\tau^{\alpha}_V(x)}
				\ar@{}[rd]^(0.55)>>{\lrcorner}
				\ar[r]^-{(\Ind_{\alpha}u)(x)}
				\ar[d]_-{l^{\alpha}_Y(x)}
				&  (\Ind_{\alpha}W)(x)
				\ar@/^7.8pc/[dd]^-{\tau^{\alpha}_W(x)}
				\ar[d]^-{\varphi(x)}                                   \\
				Y(x)
				\ar[r]^-{\theta(x)}
				\ar[d]_-{m^{\alpha}_Y(x)}
				&  (\Ind_{\alpha}W)(x) \sqcup_{(\Ind_{\alpha}V)(x)} Y(x)
				\ar@{.>}[d]^-{\psi(x)}                                 \\
				(\coInd_{\alpha}V)(x)
				\ar[r]^-{(\coInd_{\alpha}u)(x)}
				&  (\coInd_{\alpha}W)(x).
			}
			\]
			Since
			$(\coInd_{\alpha}u)(x) \circ \tau^{\alpha}_V(x)
			= \tau^{\alpha}_W(x) \circ (\Ind_{\alpha}u)(x)$
			and
			$m^{\alpha}_Y(x) \circ l^{\alpha}_Y(x)
			= \tau^{\alpha}_V(x)$, one has
			\[
			(\coInd_{\alpha}u)(x) \circ m^{\alpha}_Y(x) \circ l^{\alpha}_Y(x)
			= \tau^{\alpha}_W(x) \circ (\Ind_{\alpha}u)(x).
			\]
			By the universal property of the pushout (\ref{pushout.1}), there exists a unique morphism $\psi(x)$ such that the whole diagram commutes. Consequently, one obtains a family
			\[
			\{                 (\Ind_{\alpha}W)(x)
			\xra{\varphi(x)}
			(\Ind_{\alpha}W)(x) \sqcup_{(\Ind_{\alpha}V)(x)} Y(x)
			\xra{\psi(x)}
			(\coInd_{\alpha}W)(x)
			\}_{d(x) = \alpha}
			\]
			of left $A^0_x$-homomorphism factorizations of those $\tau^{\alpha}_W(x)$ as well as a family
			\[
			\{
			Y(x)
			\xra{\theta(x)}
			(\Ind_{\alpha}W)(x) \sqcup_{(\Ind_{\alpha}V)(x)} Y(x)
			\}_{d(x) = \alpha}
			\]
			of left $A^0_x$-modules homomorphisms such that the following diagram commutes:
			\[
			\xymatrix@R=1cm@C=0.5cm{
				(\Ind_{\alpha}V)(x)
				\ar[rr]^-{l^{\alpha}_Y(x)}
				\ar[d]_-{(\Ind_{\alpha}u)(x)}
				&&                      Y(x)
				\ar[rr]^-{m^{\alpha}_Y(x)}
				\ar[d]|-{\theta(x)}
				&&               (\coInd_{\alpha}V)(x)
				\ar[d]^-{(\coInd_{\alpha}u)(x)}               \\
				(\Ind_{\alpha}W)(x)
				\ar[rr]^-{\varphi(x)}
				&&  (\Ind_{\alpha}W)(x) \sqcup_{(\Ind_{\alpha}V)(x)} Y(x)
				\ar[rr]^-{\psi(x)}
				&&                (\coInd_{\alpha}W)(x).
			}
			\]
			Consequently, one gets a left $\UC_{\alpha+1}$-module $u_!(Y)$ in ${\UC_{\alpha+1}\Mod}_{W}$ by Theorem \ref{key 1-1 core} and a morphism $\lambda_{u,Y} : Y \to u_!(Y)$ in $\UC_{\alpha + 1}\Mod$ that is above $u$ by Remark \ref{key 1-1 core for mor}. More explicitly, for any object $x$ in $\UC_{\alpha +1}$, the morphism $\lambda_{u,Y}(x): Y(x)\to u_!(Y)(x)$ is as follows:
			\begin{equation*}
				\lambda_{u,Y}(x) =
				\begin{cases}
					u(x)              & \text{if } d(x) < \alpha, \text{ i.e. } x \in \Ob(\UC_{\alpha});      \\
					\theta(x)         & \text{if } d(x) = \alpha.
				\end{cases}
			\end{equation*}
		\end{bfhpg*}
		
		\begin{bfhpg*}[Step 2: $\lambda_{u, Y}$ is cocartesian] Consider the following diagram of solid arrows with $f$ above $v \circ u$:
			\[
			\xymatrix@R=1cm@C=1cm{
				&&  Z                                        \\
				Y
				\ar[r]_-{\lambda_{u, Y}}
				\ar@/^1pc/[rru]^(0.43){f}
				&   u_!(Y)
				\ar@{.>}[ru]_-{g}
			}
			\xymatrix@R=1cm@C=1cm{
				\\ \quad \overset{\Res_{\alpha}} \longmapsto \quad}
			\xymatrix@R=1cm@C=1cm{
				&&  U                                        \\
				V
				\ar[r]_-{u}
				\ar@/^1pc/[rru]^(0.43){v \circ u}
				&   W
				\ar[ru]_-{v}
			}
			\]
			 We have to construct a unique morphism $g$ that is above $v$ such that the left triangle in $\UC_{\alpha + 1}\Mod$ commutes. To this end, apply the commutative diagram
			(\ref{Fac of NT}) of natural transformations to $f$, one obtains the following commutative diagram
			\[
			\xymatrix@R=1cm@C=0.5cm{
				\Ind_{\alpha}V = \Ind_{\alpha}\Res_{\alpha}Y
				\ar[rr]^-{l^{\alpha}_Y}
				\ar[d]_-{\Ind_{\alpha}(v \circ u) = \Ind_{\alpha}\Res_{\alpha}f}
				&&                              Y
				\ar[rr]^-{m^{\alpha}_Y}
				\ar[d]|-{f}
				&&              \coInd_{\alpha}\Res_{\alpha}Y = \coInd_{\alpha}V
				\ar[d]^-{\coInd_{\alpha}\Res_{\alpha}f = \coInd_{\alpha}(v\circ u)}               \\
				\Ind_{\alpha}U = \Ind_{\alpha}\Res_{\alpha}Z
				\ar[rr]^-{l^{\alpha}_Z}
				&&                              Z
				\ar[rr]^-{m^{\alpha}_Z}
				&&               \coInd_{\alpha}\Res_{\alpha}Z = \coInd_{\alpha}U
			}
			\]
			in $\UC_{\alpha + 1}\Mod$.
			Restricting further to each object $x$ in $\UC_{\alpha +1}$ with $d(x) = \alpha$, one obtains the following commutative diagram
			\begin{equation*}\label{back commu.1}\tag{\ref{res is Gro bifib}.2}
				\begin{gathered}
					\xymatrix@R=1cm@C=0.5cm{
						(\Ind_{\alpha}V)(x)
						\ar[rr]^-{l^{\alpha}_Y(x)}
						\ar[d]_-{(\Ind_{\alpha}(v\circ u))(x)}
						&&                      Y(x)
						\ar[rr]^-{m^{\alpha}_Y(x)}
						\ar[d]|-{f(x)}
						&&               (\coInd_{\alpha}V)(x)
						\ar[d]^-{(\coInd_{\alpha}(v\circ u))(x)}             \\
						(\Ind_{\alpha}U)(x)
						\ar[rr]^-{l^{\alpha}_Z(x)}
						&&                   Z(x)
						\ar[rr]^-{m^{\alpha}_Z(x)}
						&&                (\coInd_{\alpha}U)(x)
					}
				\end{gathered}
			\end{equation*}
			in $A^0_x\Mod$, which fits into the following diagram of solid arrows:
			\[
			\xymatrix@R=0.35cm@C=0.25cm{
				&&&&&
				(\Ind_{\alpha}U)(x)
				\ar[dd]^(0.6){l^{\alpha}_Z(x)}
				\\
				(\Ind_{\alpha}V)(x)
				\ar[rrr]_(0.45){(\Ind_{\alpha}u)(x)}
				\ar@/^1pc/[rrrrru]^(0.7){(\Ind_{\alpha}{v \circ u})(x)}
				\ar[dd]_(0.7){l^{\alpha}_Y(x)}
				&&&
				(\Ind_{\alpha}W)(x)
				\ar[rru]_{(\Ind_{\alpha}v)(x)}
				\ar[dd]_(0.7){\varphi(x)}
				\\
				&&&   \ar@/^0.25pc/^(0.35){f(x)}[rr]
				&&
				Z(x)
				\ar[dd]^(0.6){m^{\alpha}_Z(x)}                                   \\
				Y(x)
				\ar[rrr]_(0.3){\theta(x)}
				\ar@/^0.45pc/@{}[rrru]
				\ar[dd]_(0.6){m^{\alpha}_Y(x)}
				&&&
				(\Ind_{\alpha}W)(x) \sqcup_{(\Ind_{\alpha}V)(x)} Y(x)
				\ar@{.>}[rru]_-{\delta(x)}
				\ar[dd]_(0.7){\psi(x)}                                             \\
				&&&&&            (\coInd_{\alpha}U)(x)
				\\
				(\coInd_{\alpha}V)(x)
				\ar[rrr]_(0.5){(\coInd_{\alpha}u)(x)}
				\ar@/^1pc/[rrrrru]^(0.7){(\coInd_{\alpha}{v \circ u})(x)}
				&&&              (\coInd_{\alpha}W)(x)
				\ar[rru]_-{(\coInd_{\alpha}v)(x)}
			}
			\]
			By the commutativity of the left square of (\ref{back commu.1}) and the universal property of the pushout (\ref{pushout.1}), there exists a unique morphism $\delta(x)$ such that the upper half of the above diagram commutes. It is routine to check that the equality
			$$m^{\alpha}_Z(x) \circ \delta(x) = (\coInd_{\alpha}v)(x) \circ \psi(x)$$
			holds by the universal property of the pushout (\ref{pushout.1}) again. Thus, one gets a family
			\[
			\{
			(\Ind_{\alpha}W)(x) \sqcup_{(\Ind_{\alpha}V)(x)} Y(x)
			\xra{\delta(x)}
			Z(x)
			\}_{d(x) = \alpha}
			\]
			of left $A^0_x$-modules homomorphisms such that the following diagram commutes:
			\[
			\xymatrix@R=1cm@C=0.35cm{
				(\Ind_{\alpha}W)(x)
				\ar[rr]^-{\varphi(x)}
				\ar[d]_-{(\Ind_{\alpha}v)(x)}
				&&   (\Ind_{\alpha}W)(x) \sqcup_{(\Ind_{\alpha}V)(x)} Y(x)
				\ar[rr]^-{\psi(x)}
				\ar[d]|-{\delta(x)}
				&&               (\coInd_{\alpha}W)(x)
				\ar[d]^-{(\coInd_{\alpha}v)(x)}               \\
				(\Ind_{\alpha}U)(x)
				\ar[rr]^-{l^{\alpha}_Z(x)}
				&&                       Z(x)
				\ar[rr]^-{m^{\alpha}_Z(x)}
				&&               (\coInd_{\alpha}U)(x)
			}
			\]
			According to Remark \ref{key 1-1 core for mor}, one obtains a (unique) morphism
			$g : u_!(Y) \to Z$ in $\UC_{\alpha + 1}\Mod$ that is above $v$. More explicitly,
			for any object $x$ in $\UC_{\alpha +1}$, the morphism $g(x): u_!(Y)(x) \to Z(x)$ is as follows:
			\begin{equation*}
				g(x) =
				\begin{cases}
					v(x)        & \text{if } d(x) < \alpha, \text{ i.e. } x \in \Ob(\UC_{\alpha});      \\
					\delta(x)   & \text{if } d(x) = \alpha.
				\end{cases}
			\end{equation*}
			Now it is routine to check that $g \circ \lambda_{u, Y} = f$, as desired.
		\end{bfhpg*}

To verify (G1) in Definition
\ref{df of Grothendieck opfibra},
let $u:V\to W$ and let $Z$ lie over $W$.
For each $x$ with $d(x)=\alpha$, form the
pullback
$P_x=Z(x)\times_{(\coInd_\alpha W)(x)}
(\coInd_\alpha V)(x)$,
with projections $p_x$ and $q_x$.
Since
\[
m_Z^\alpha(x)l_Z^\alpha(x)
(\Ind_\alpha u)(x)
=(\coInd_\alpha u)(x)\tau_V^\alpha(x),
\]
there is a unique morphism
$h_x:(\Ind_\alpha V)(x)\to P_x$
such that
$p_xh_x=l_Z^\alpha(x)(\Ind_\alpha u)(x)$ and
$q_xh_x=\tau_V^\alpha(x)$.
By Theorem~\ref{key 1-1 core}, the
factorizations
\[
(\Ind_\alpha V)(x)\xrightarrow{h_x}
P_x\xrightarrow{q_x}
(\coInd_\alpha V)(x)
\]
define a module $u^*(Z)$ over $V$.
The projections $p_x$, together with $u$,
induce a morphism
$\rho_{u,Z}:u^*(Z)\to Z$ over $u$.
We claim that $\rho_{u,Z}$ is cartesian.
Indeed, let $v:U\to V$ and let
$f:X\to Z$ lie over $uv$.
For each $x$ with $d(x)=\alpha$, the maps
$f(x):X(x)\to Z(x)$
and
$(\coInd_\alpha v)(x)m_X^\alpha(x):
X(x)\to(\coInd_\alpha V)(x)$
induce a unique morphism $X(x)\to P_x$.
By the naturality of the latching and
matching maps, these morphisms are
compatible with the above factorizations.
Consequently, they define a unique
morphism $X\to u^*(Z)$ over $v$
whose composite with $\rho_{u,Z}$
is $f$. Thus $\rho_{u,Z}$ is cartesian.

Together with (G2), this shows that
$\Res_\alpha$ is a Grothendieck
bifibration.
	\end{prf*}
	
	\begin{remark}\label{des of fiber mor}
		Let $f : Y \to Z$ be a morphism in $\UC_{\alpha +1}\Mod$ that is above the morphism
		$u : V \to W$ in $\UC_{\alpha}\Mod$. Then by Lemma \ref{key factorization}, $f$ can be factored uniquely as
		\begin{prt}
			\item[$\bullet$]
			the cocartesian morphism $\lambda_{u, Y}: Y \to u_!(Y)$ followed by a morphism $f_\triangleright: u_!(Y)\to Z$ in the fiber $\UC_{\alpha +1}\Mod_W$, and
			\item[$\bullet$]
			a morphism $f^{\triangleleft}: Y \to u^*(Z)$ in the fiber $\UC_{\alpha +1}\Mod_V$ followed by the cartesian morphism $\rho_{u, Z}: u^*(Z) \to Z$.
		\end{prt}
		Under the identification given by the isomorphism functor $\scrR$ established in Corollary \ref{the iso of fiber},
		\begin{prt}
			\item[$\bullet$]
			the fiber morphism $f_{\triangleright}$ is the family
			\[
			\{
			(\Ind_{\alpha}W)(x) \sqcup_{(\Ind_{\alpha}V)(x)} Y(x)
			\xra{\sigma(x)}
			Z(x)
			\}_{d(x) = \alpha}
			\]
			of morphisms of left $A^0_x$-modules, where $\sigma(x)$ is the unique morphism, obtained by the universal property of the pushout (\ref{pushout.1}), such that the following diagram in $A^0_x\Mod$ commutes:
			\begin{equation*}\label{the fiber}\tag{\ref{des of fiber mor}.1}
				\begin{gathered}
					\xymatrix@R=0.35cm@C=0.25cm{
						&&&&&
						(\Ind_{\alpha}W)(x)
						\ar[dd]^(0.6){l^{\alpha}_Z(x)}
						\\
						(\Ind_{\alpha}V)(x)
						\ar[rrr]_(0.45){(\Ind_{\alpha}u)(x)}
						\ar@/^1pc/[rrrrru]^(0.7){(\Ind_{\alpha}u)(x)}
						\ar[dd]_(0.7){l^{\alpha}_Y(x)}
						&&&               (\Ind_{\alpha}W)(x)
						\ar@{=}[rru]
						\ar[dd]_(0.7){\varphi(x)}                                     \\
						&&&
						\ar@/^0.25pc/^(0.35){f(x)}[rr]
						&&                       Z(x)
						\ar[dd]^(0.6){m^{\alpha}_Z(x)}                                \\
						Y(x)
						\ar[rrr]_(0.3){\theta(x)}
						\ar@/^0.45pc/@{}[rrru]
						\ar[dd]_(0.6){m^{\alpha}_Y(x)}
						&&&  (\Ind_{\alpha}W)(x) \sqcup_{(\Ind_{\alpha}V)(x)} Y(x)
						\ar@{.>}[rru]_-{\sigma(x)}
						\ar[dd]_(0.7){\psi(x)}                                         \\
						&&&   \ar@/^0.25pc/^(0.35){(\coInd_{\alpha}u)(x)}[rr]
						&&              (\coInd_{\alpha}W)(x).
						\\
						(\coInd_{\alpha}V)(x)
						\ar[rrr]_(0.5){(\coInd_{\alpha}u)(x)}
						\ar@/^0.45pc/@{}[rrru]
						&&&
						(\coInd_{\alpha}W)(x)
						\ar@{=}[rru]
					}
				\end{gathered}
			\end{equation*}
			
			\item[$\bullet$]
			dually, the fiber morphism $f^{\triangleleft}$ is the family
			\[
			\{
			Y(x)
			\xra{\varsigma(x)}
			Z(x)\times_{(\coInd_{\alpha}W)(x)}(\coInd_{\alpha}V)(x)
			\}_{d(x) = \alpha}
			\]
			of morphisms of left $A^0_x$-modules, where $\varsigma(x)$ is the unique morphism, obtained by the universal property of the pullback
			\begin{equation*}\label{pullback1}\tag{\ref{des of fiber mor}.2}
				\begin{gathered}
					\xymatrix@R=1cm@C=0.35cm{
						Z(x)\times_{(\coInd_{\alpha}W)(x)}(\coInd_{\alpha}V)(x)
						\ar[r]^-{ }
						\ar[d]_-{ }
						&  Z(x)
						\ar[d]^-{m^{\alpha}_Z(x)}
						\\
						(\coInd_{\alpha}V)(x)
						\ar[r]^-{(\coInd_{\alpha}u)(x)}
						&  (\coInd_{\alpha}W)(x)
						\ar@{}[ul]^>>{\ulcorner}
					}
				\end{gathered}
			\end{equation*}
			
			such that the following diagram in $A^0_x\Mod$ commutes:
			\[
			\xymatrix@R=0.35cm@C=0.1cm{
				(\Ind_{\alpha}V)(x)
				\ar@{=}[rrd]
				\ar@/^1pc/[rrrrrrd]^(0.7){}
				\ar[dd]                                     \\
				&&     (\Ind_{\alpha}V)(x)
				\ar[rrrr]_{}
				\ar[dd]
				&&&&   (\Ind_{\alpha}W)(x)
				\ar[dd]                                     \\
				Y(x)
				\ar[dd]
				\ar@/^0.3pc/@{}[rr]
				\ar@{.>}[rrd]_(0.25){\varsigma(x)}
				&&        \ar@/^0.3pc/[rrrrd]               \\
				&&  Z(x)\times_{(\coInd_{\alpha}W)(x)}(\coInd_{\alpha}V)(x)
				\ar[rrrr]^{}
				\ar[dd]
				&&&&          Z(x)
				\ar[dd]
				\\
				(\coInd_{\alpha}V)(x)
				\ar@/^0.3pc/@{}[rr]
				\ar@{=}[rrd]
				&&        \ar@/^0.3pc/[rrrrd]
				\\
				&&     (\coInd_{\alpha}V)(x)
				\ar[rrrr]^{}
				&&&&   (\coInd_{\alpha}W)(x).
			}
			\]
		\end{prt}
	\end{remark}
	
	\begin{remark}\label{des of u(k,l)}
		Let $u : V \to W$ be a morphism in $\UC_{\alpha}\Mod$. Under the identification given by $\scrR$ established in Corollary \ref{the iso of fiber}, we can give explicit descriptions of $u_!(l)$ for any morphism $l: Y \to Y'$ in $\UC_{\alpha + 1}\Mod_V$ as well as $u^*(k)$ for any morphism $k: Z' \to Z$ in $\UC_{\alpha + 1}\Mod_W$.
		
		Recall from Remark \ref{funs u! and u^*} that $u_!(l)$ is the fiber morphism
		$(\lambda_{u,Y'} \circ l)_\triangleright$ in $\UC_{\alpha +1}\Mod_W$, which modulo identification $\scrR$, is the family
		\[
		\{
		(\Ind_{\alpha}W)(x) \sqcup_{(\Ind_{\alpha}V)(x)} Y(x)
		\xra{\varrho(x)}
		(\Ind_{\alpha}W)(x) \sqcup_{(\Ind_{\alpha}V)(x)} Y'(x)
		\}_{d(x) = \alpha}
		\]
		of morphisms of left $A^0_x$-modules, where $\varrho(x)$ is the unique morphism given by the universal property of the pushout (\ref{pushout.1}) such that the following diagram in $A^0_x\Mod$ commutes (see Remark \ref{des of fiber mor}):
		\[
		\xymatrix@R=0.35cm@C=0.01cm{
			&    (\Ind_{\alpha}V)(x)
			\ar[rr]^{}
			\ar@/^0pc/@{}[d]
			&&   (\Ind_{\alpha}W)(x)
			\ar[dd]^{}                                                 \\
			(\Ind_{\alpha}V)(x)
			\ar[dd]^{}
			\ar[rr]^{}
			\ar@{=}[ru]^{}
			&
			\ar[d]^{}
			&    (\Ind_{\alpha}W)(x)
			\ar@{=}[ru]^{}
			\ar[dd]^{}                                                   \\
			&    Y'(x)
			\ar@/^0pc/@{}[d]
			\ar@/^0pc/@{}[r]^(0.5){}
			&
			\ar[r]
			&    (\Ind_{\alpha}W)(x) \sqcup_{(\Ind_{\alpha}V)(x)} Y'(x)
			\ar[dd]^{}                                                     \\
			Y(x)
			\ar[ru]^{l(x)}
			\ar[rr]^(0.3){}
			\ar[dd]^{}
			&
			\ar[d]^{}
			&    (\Ind_{\alpha}W)(x) \sqcup_{(\Ind_{\alpha}V)(x)} Y(x)
			\ar@{.>}[ru]^(0.4){\varrho(x)}
			\ar[dd]^{}                                                       \\
			&    (\coInd_{\alpha}V)(x)
			\ar@/^0pc/@{}[r]
			&
			\ar[r]
			&    (\coInd_{\alpha}W)(x)                                        \\
			(\coInd_{\alpha}V)(x)
			\ar[rr]^{}
			\ar@{=}[ru]^{}
			&&   (\coInd_{\alpha}W)(x)
			\ar@{=}[ru]^{}
		}
		\]
		Note that both upper squares in the front and back of the above diagram are pushouts. By the pasting lemma, the square
		\begin{equation*}\label{ind pushout}\tag{\ref{des of u(k,l)}.1}
			\begin{gathered}
				\xymatrix@R=1cm@C=0.75cm{
					Y(x)
					\ar@{}[rd]^(0.75)>>{\lrcorner}
					\ar[r]^-{}
					\ar[d]_-{l(x)}
					&  (\Ind_{\alpha}W)(x) \sqcup_{(\Ind_{\alpha}V)(x)} Y(x)
					\ar[d]^-{\varrho(x)\, = \,(u_!(l))(x)}                             \\
					Y'(x)
					\ar[r]^-{}
					&  (\Ind_{\alpha}W)(x) \sqcup_{(\Ind_{\alpha}V)(x)} Y'(x)
				}
			\end{gathered}
		\end{equation*}
		is also a pushout.
		
		Dually, by Remark \ref{funs u! and u^*}, $u^*(k)$ is the fiber morphism
		$(k \circ \rho_{u,Z'})^{\triangleleft} \in \UC_{\alpha +1}\Mod_V$,
		which modulo identification $\scrR$, is the family
		\[
		\{
		Z'(x)\times_{(\coInd_{\alpha}W)(x)}(\coInd_{\alpha}V)(x)
		\xra{\zeta(x)}
		Z(x)\times_{(\coInd_{\alpha}W)(x)}(\coInd_{\alpha}V)(x)
		\}_{d(x) = \alpha}
		\]
		of morphisms of left $A^0_x$-modules, where $\zeta(x)$ is the unique morphism given by the universal property of the pullback (\ref{pullback1}), such that the following diagram in $A^0_x\Mod$ commutes:
		\[
		\xymatrix@R=0.35cm@C=1mm{
			\Ind_{\alpha}V(x)
			\ar[rr]^-{}
			\ar@{=}[dr]
			\ar[dd]^-{}
			&&                 \Ind_{\alpha}W(x)
			\ar@{=}[dr]
			\ar@/^0pc/@{}[d]                                                        \\
			&                  \Ind_{\alpha}V(x)
			\ar[rr]^-{}
			\ar[dd]^-{}
			&    \ar[d]
			&                  \Ind_{\alpha}W(x)
			\ar[dd]
			\\
			Z'(x)\times_{(\coInd_{\alpha}W)(x)}(\coInd_{\alpha}V)(x)
			\ar[dd]^-{}
			\ar@/^0pc/@{}[r]
			\ar@{.>}[dr]^-{\zeta(x)}
			&    \ar[r]
			&                  Z'(x)
			\ar@/^0pc/@{}[d]
			\ar[dr]^-{k(x)}                                                            \\
			&                  Z(x)\times_{(\coInd_{\alpha}W)(x)}(\coInd_{\alpha}V)(x)
			\ar[dd]^-{}
			\ar[rr]^-{}
			&    \ar[d]
			&                  Z(x)
			\ar[dd]^-{}                                                                 \\
			\coInd_{\alpha}V(x)
			\ar@/^0pc/@{}[r]
			\ar@{=}[dr]
			&    \ar[r]
			&                  \coInd_{\alpha}W(x)
			\ar@{=}[dr]                                                                 \\
			&                  \coInd_{\alpha}V(x)
			\ar[rr]^-{}
			&&                 \coInd_{\alpha}W(x)
		}
		\]
		Since both lower squares in the front and back of the above diagram are pullbacks, the square
		\[
		\xymatrix@R=1cm@C=0.75cm{
			Z'(x)\times_{(\coInd_{\alpha}W)(x)}(\coInd_{\alpha}V)(x)
			\ar[r]^-{}
			\ar[d]_-{(u^*(k))(x) \, = \, \zeta(x)}
			&                       Z'(x)
			\ar[d]^-{k(x)}                                                         \\
			Z(x)\times_{(\coInd_{\alpha}W)(x)}(\coInd_{\alpha}V)(x)
			\ar[r]^-{}
			&                       Z(x)
			\ar@{}[ul]^>>{\ulcorner}
		}
		\]
		is also a pullback.
	\end{remark}
	
	\section{Lifting of complete cotorsion pairs}\label{lift CCP}
	\noindent
	In this section, our main goal is to show that under some mild assumptions, a family
	$\{(\calC_x, \calD_x) \}_{x \in \Ob(\UC)}$ of complete cotorsion pairs in those $A_x^0 \Mod$'s induces a complete cotorsion pair in $\UC \Mod$.
	
	We begin by introducing the following classes of objects, which constitute the key components of the above mentioned induced complete cotorsion pair as well as the induced Hovey triple given in the next section. We refer the reader to \ref{Adjoint tri} for the adjoint triple induced by the embedding $\iota_{\alpha} : \UC_{\alpha} \to \UC_{\beta}$ with
	$\alpha < \beta \leqslant \lambda+1$ and for some relevant facts.
	
	\begin{definition}\label{Phi and Psi}
		Let $\calS = \{ \calS_x \}_{x \in \Ob(\UC)}$ be a family with each $\calS_x$ a class of objects in $A_x^0 \Mod$. For each ordinal $\beta$ with $\beta \leqslant \lambda+1$, we define the following classes of objects in $\UC_{\beta}\Mod$:
		\begin{align*}
			\Phi_{\beta}(\calS)
			& =\left\{ Y \in \UC_{\beta}\Mod \:
		\left|
		\begin{array}{c}
			\text{for each}\ x \in \Ob(\UC_{\beta}) \text{ with } d(x)= \alpha, \\
			l^{\alpha}_Y(x)\ \text{is monic and}\ \coker(l^{\alpha}_Y(x)) \in \calS_x
		\end{array}
		\right.
		\right\}, \ \mathrm{and} \\
		\Psi_{\beta}(\calS)
		&=\left\{ Y \in \UC_{\beta}\Mod \:
		\left|
		\begin{array}{c}
			\text{for each}\ x \in \Ob(\UC_{\beta}) \text{ with } d(x)= \alpha, \\
			m^{\alpha}_Y(x)\ \text{is epic and}\ \ker(m^{\alpha}_Y(x)) \in \calS_x
		\end{array}
		\right.
		\right\}.
		\end{align*}
		Alternatively, by Corollary \ref{cokernel of latching}, $\Phi_{\beta}(\calS)$ and $\Psi_{\beta}(\calS)$ can be described as follows:
		\begin{align*}
		\Phi_{\beta}(\calS)
		&= \{
		Y \in \UC_{\beta}\Mod \,| \
		\forall x \in \Ob(\UC_{\beta}),  \,
		\Tor_1^{\UC_{\beta}}(\Delta^x, Y) = 0,
		\Delta^x \otimes_{\UC_{\beta}} Y \in \calS_x
		\}, \ \mathrm{and} \\
		\Psi_{\beta}(\calS)
		&= \{
		Y \in \UC_{\beta}\Mod \,| \
		\forall x \in \Ob(\UC_{\beta}),  \,
		\Ext^1_{\UC_{\beta}}(\Delta_x, Y) = 0,
		\Hom_{\UC_{\beta}}(\Delta_x, Y) \in \calS_x
		\}.
		\end{align*}
		If $\beta = \lambda+1$, then $\UC_{\beta}\Mod = \UC \Mod$ obviously. In this case, we use the symbols $\Phi(\calS)$ and $\Psi(\calS)$ for short instead of $\Phi_{\beta}(\calS)$ and $\Psi_{\beta}(\calS)$.
	\end{definition}

The next result will be used frequently in the rest of the paper, which can be proved easily.

\begin{proposition}\label{useful}
Let $\Res_{\alpha} : \UC_{\beta}\Mod \to \UC_{\alpha}\Mod$ be the restriction functor
		with respect to the embedding $\iota_{\alpha} : \UC_{\alpha} \to \UC_{\beta}$ with $\alpha < \beta \leqslant \lambda+1$. Then for each family $\calS = \{ \calS_x \}_{x \in \Ob(\UC)}$ of classes of objects in those $A^0_x\Mod$'s, one has $\Res_{\alpha}(\Phi_{\beta}(\calS))\subseteq\Phi_{\alpha}(\calS)$ and $\Res_{\alpha}(\Psi_{\beta}(\calS))\subseteq\Psi_{\alpha}(\calS)$.
\end{proposition}
	
	The main result in this section is as follows:
	
	\begin{theorem}\label{ma re for ccp}
		Suppose that $\UC^+$ is a projective right $\UC^0$-module and $\UC^-$ is a projective left $\UC^0$-module. If $(\calC, \calD) = \{(\calC_x, \calD_x)\}_{x \in \Ob(\UC)}$ is a family of complete cotorsion pairs in those $A^0_x\Mod$'s, then $(\Phi(\calC), \Psi(\calD))$ forms a complete cotorsion pair in $\UC \Mod$.
	\end{theorem}
	
	As one may expect, the proof of Theorem \ref{ma re for ccp} relies on a transfinite induction. To carry out the induction step, we need some preparations. The next result is essentially a dual version of Remark \ref{a general filtration}.
	
	\begin{lemma}\label{filtration1}
		Suppose that $\UC^+$ is a projective right $\UC^0$-module and $\UC^-$ is a projective left $\UC^0$-module. Then for each $x\in\Ob(\UC)$ with $d(x)=\alpha$, the right $\UC_{\alpha}$-module
		$1_x \UC_{\beta} = \UC_{\beta}(\iota_{\alpha}(-), x)$ admits a filtration of right $\UC_{\alpha}$-modules with each factor
		$$\bigoplus_{d(z) = \gamma} \UC^+_{\beta}(z, x) \otimes_{A^0_z} \Delta^z$$
for $0 \leqslant \gamma < \alpha$, where $\Delta^z$ is the standard right $\UC_\alpha$-module.
\end{lemma}
\begin{prf*}
		Consider the right $\UC^+_{\alpha}$-module $1_x \UC^+_{\beta}$. Since $\UC^+$ is a projective right $\UC^0$-module, it follows from Lemma \ref{projectivity}(a) that $\UC^+_{\beta}(z, x)$ is a projective right $A^0_z$-module for each $z \in \Ob(\UC_{\alpha})$. By the dual version of Lemma \ref{filtration of C^-modules} (taking $\UC = \UC_{\alpha}$), the right $\UC^+_{\alpha}$-module $1_x \UC^+_{\beta}$ admits a filtration
		\[
		0 = (1_x \UC^+_{\beta})_{< 0} \subseteq
		(1_x \UC^+_{\beta})_{< 1}       \subseteq\cdots\subseteq
		(1_x \UC^+_{\beta})_{< \alpha}
		= 1_x \UC^+_{\beta}
		\]
		such that each factor is
		$$(1_x \UC^+_{\beta})_{< \gamma+1}/(1_x \UC^+_{\beta})_{< \gamma}=\bigoplus_{d(z) = \gamma} \UC^+_{\beta}(z, x)$$
		for $0 \leqslant \gamma < \alpha$.
		Since $\UC^-$ is a projective left $\UC^0$-module, so is $\UC^-_{\alpha}$. Hence, by
		Corollary \ref{exactness of induction}, $- \otimes_{\UC^+_{\alpha}} \UC_{\alpha}$ is exact. Applying it to the above filtration of $1_x \UC^+_{\beta}$ and noting that
		$1_x \UC_{\beta} \cong 1_x \UC^+_{\beta} \otimes_{\UC^+_{\alpha}} \UC_{\alpha}$
		by Lemma \ref{T AND T+}, one gets a filtration of the right $\UC_{\alpha}$-module
		$1_x \UC_{\beta}$ with each factor
		\begin{align*}
			(\bigoplus_{d(z) = \gamma}
			\UC^+_{\beta}(z, x)) \otimes_{\UC^+_{\alpha}} \UC_{\alpha}
			& \cong
			\bigoplus_{d(z) = \gamma}
			( \UC^+_{\beta}(z, x) \otimes_{A^0_z} (A^0_z \otimes_{\UC^+_{\alpha}} \UC_{\alpha}) )
			\\
			& \cong
			\bigoplus_{d(z) = \gamma}
			( \UC^+_{\beta}(z, x) \otimes_{A^0_z} \Delta^z  )
		\end{align*}
		for $0 \leqslant \gamma < \alpha$, where the second isomorphism holds by Theorem \ref{dual version}(a).
	\end{prf*}
	
	The following technical result is well known; see \cite[Lemma 6.7]{GS} for a proof.
	
	\begin{lemma}\label{fil tor van}
		Let $U$ be a left $\UC$-module and let $W$ be a right $\UC$-module admitting a filtration
		$$0=W_{0}\xra{\iota_{0,1}} W_1\xra{\iota_{1,2}} W_2 \to \cdots \to W_\mu=W$$
		such that $\Tor^{\UC}_1(\coker(\iota_{\omega, \omega + 1}), U) = 0$ for each $\omega < \mu$. Then one has $\Tor^{\UC}_1(W, U) = 0$.
	\end{lemma}
	
	The next result will play a key role in the proof of Theorem \ref{ma re for ccp}.
	
	
	\begin{lemma}\label{ind keep exact}
		Suppose that $\UC^+$ is a projective right $\UC^0$-module and $\UC^-$ is a projective left $\UC^0$-module. Let $x\in\Ob{\UC}$ with $d(x) = \alpha$, and let $0 \to V \to W \to U \to 0$ be an exact sequence in $\UC_{\alpha}\Mod$ such that $\Tor^{\UC_{\alpha}}_1 (\Delta^z, U) = 0$ for all $z\in\Ob(\UC_{\alpha})$. Then there is a short exact sequence in $A^0_x\Mod$:
		\[
		0 \to (\Ind_{\alpha}V)(x) \to (\Ind_{\alpha}W)(x) \to (\Ind_{\alpha}U)(x) \to 0
		\]
	\end{lemma}
	
	\begin{prf*}
		Applying $1_x \UC_{\alpha + 1} \otimes_{\UC_{\alpha}}-$ to the given short exact sequence, one gets an exact sequence
		\[
		\cdots   \to \Tor_1^{\UC_{\alpha}} (1_x \UC_{\alpha + 1}, U)   \to
		1_x \UC_{\alpha + 1} \otimes_{\UC_{\alpha}} V         \to
		1_x \UC_{\alpha + 1} \otimes_{\UC_{\alpha}} W         \to
		1_x \UC_{\alpha + 1} \otimes_{\UC_{\alpha}} U         \to   0.
		\]
		It suffices to prove $\Tor_1^{\UC_{\alpha}} (1_x \UC_{\alpha + 1}, U) = 0$. But the right $\UC_{\alpha}$-module $1_x \UC_{\alpha + 1}$ has a filtration of right $\UC_{\alpha}$-modules with each factor
		$\bigoplus_{d(z) = \gamma} (\UC^+_{\alpha + 1}(z, x) \otimes_{A^0_z} \Delta^z)$ for $0 \leqslant \gamma < \alpha$ (see Lemma \ref{filtration1}), so one only needs to show that $\Tor^{\UC_{\alpha}}_1 (\UC^+_{\alpha + 1}(z, x) \otimes_{A_z^0} \Delta^z, U) = 0$ for each $z \in \Ob(\UC_{\alpha})$ by Lemma \ref{fil tor van}. Since $\UC^+$ is a projective right $\UC^0$-module, it follows from Lemma \ref{projectivity} (a) that $\UC^+_{\alpha + 1}(z, x)$ is a projective right $A^0_z$-module. Hence, one has $\Tor^{\UC_{\alpha}}_1 (\UC^+_{\alpha + 1}(z, x) \otimes_{A_z^0} \Delta^z, U)\cong \UC^+_{\alpha + 1}(z, x) \otimes_{A_z^0} \Tor^{\UC_{\alpha}}_1(\Delta^z, U)= 0$ as $\Tor^{\UC_{\alpha}}_1 (\Delta^z, U) = 0$ by the assumption.
	\end{prf*}
	
\begin{lemma}\label{matching-exact-on-psi}
Suppose that $\UC^+$ is a projective right $\UC^0$-module and $\UC^-$ is a projective left $\UC^0$-module. Let $x\in\Ob{\UC}$ with $d(x) = \alpha$, and let
$0\to K\to V\to W
\to 0$
be an exact sequence in $\UC_\alpha\Mod$ such that
$K\in\Psi_\alpha(\calD)$. Then there is a short exact sequence in $A_x^0\Mod$:
\[
0\to\Hom_{\UC_\alpha}(\UC_{\alpha+1}1_x,K)\to \Hom_{\UC_\alpha}(\UC_{\alpha+1}1_x,V)
\to \Hom_{\UC_\alpha}(\UC_{\alpha+1}1_x,W)\to 0.
\]
\end{lemma}

\begin{prf*}
By the Reedy factorization and the exactness
of induction, the left $\UC_\alpha$-module
$\UC_{\alpha+1}1_x$ admits a transfinite
filtration with factors
$\bigoplus_{d(z)=\gamma}
\Delta_z\otimes_{A_z^0}\UC^-(x,z)$, where
$\gamma<\alpha$.
Since $K\in\Psi_\alpha(\calD)$, one has
$\Ext^1_{\UC_\alpha}(\Delta_z,K)=0$
for all $z$ with $d(z)<\alpha$. Moreover,
$\UC^-(x,z)$ is a projective left
$A_z^0$-module. Thus, by Eklof's lemma, one has
$\Ext^1_{\UC_\alpha}
(\UC_{\alpha+1}1_x,K)=0$.
Applying
$\Hom_{\UC_\alpha}(\UC_{\alpha+1}1_x,-)$
to the given short exact sequence
yields the desired conclusion.
\end{prf*}

	Now we are ready to give a proof of Theorem \ref{ma re for ccp}.
	
	\begin{bfhpg}[Proof of Theorem \ref{ma re for ccp}]
		We prove that for each ordinal $\beta \leqslant \lambda+1$, $(\Phi_{\beta}(\calC), \Psi_{\beta}(\calD))$ forms a complete cotorsion pair in $\UC_{\beta}\Mod$. The conclusion then follows.		
		\begin{bfhpg*}[The case that $\beta = 1$]
			In this case $\UC_1$ only contains objects of degree zero, so
			\[\UC_1\Mod \cong \prod_{x \in \Ob(\UC_1)} A_x^0 \Mod.
			\]
			Meanwhile, for each $Y \in \UC_1\Mod$ and each $x \in \Ob(\UC_1)$, it is clear that
			\[
			1_x \UC_{1} \otimes_{\UC_{0}} \Res_{0}Y = 0 = \Hom_{\UC_{0}}(\UC_{1} 1_x, \Res_{0}Y)
			\]
			as $\UC_{0} = \emptyset$, so both $l^0_Y(x)$ and $m^0_Y(x)$ are zero. Consequently, one has
			\[
			\Phi_1(\calC) \cong \prod_{x \in \Ob(\UC_1)}\calC_x
			\quad \mathrm{and} \quad
			\Psi_1(\calD) \cong \prod_{x \in \Ob(\UC_1)}\calD_x.
			\]
			But each $(\calC_x, \calD_x)$ is a complete cotorsion pair in $A_x^0 \Mod$, so $(\Phi_1(\calC), \Psi_1(\calD))$ is also a complete cotorsion pair in $\UC_1\Mod$.
		\end{bfhpg*}
		
		\begin{bfhpg*}[The case that $\beta = \alpha + 1$ is a successor ordinal] We carry out the induction procedure with the following strategy. By Lemma \ref{bri lem}, it suffices to show that
			$(\Mon(\Phi_{\alpha+1}(\calC)), \ \Epi(\Psi_{\alpha+1}(\calD)))$
			is a weak factorization system in $\UC_{\alpha+1}\Mod$. Since $(\Phi_{\alpha}(\calC), \Psi_{\alpha}(\calD))$ is a complete cotorsion pair in $\UC_{\alpha}\Mod$ by the induction hypothesis, it follows from Lemma \ref{bri lem} that $(\Mon(\Phi_{\alpha}(\calC)), \Epi(\Psi_{\alpha}(\calD))$ is a weak factorization system in $\UC_{\alpha}\Mod$. But the restriction $\Res_{\alpha} : \UC_{\alpha+1}\Mod \to \UC_{\alpha}\Mod$ is a Grothendieck bifibration by Theorem \ref{res is Gro bifib}. Thus, one may use Lemma \ref{key lem for con wfs} to construct an induced weak factorization system $(\sfC_{\UC_{\alpha+1}\Mod}, \ \sfD_{\UC_{\alpha+1}\Mod})$
			in $\UC_{\alpha+1}\Mod$, and show further that
			$\Mon(\Phi_{\alpha+1}(\calC)) = \sfC_{\UC_{\alpha+1}\Mod}$ and $\Epi(\Psi_{\alpha+1}(\calD)) = \sfD_{\UC_{\alpha+1}\Mod}$.
			
			\begin{bfhpg*}
				[Step 1: Construct the weak factorization system
				$(\sfC_{\UC_{\alpha+1}\Mod}, \sfD_{\UC_{\alpha+1}\Mod})$]
				Firstly, we show that the fiber ${\UC_{\alpha+1}\Mod}_V$ admits a weak factorization system for each $V \in \UC_{\alpha}\Mod$. It follows from Corollary \ref{the iso of fiber} that the functor $\scrR$ gives an isomorphism
				\[
				{\UC_{\alpha+1}\Mod}_V \cong \prod_{d(x) = \alpha}
				{_{(\Ind_{\alpha}V)(x)}}
				\backslash{\substack{A^0_x\Mod \\ \tau^{\alpha}_V(x)}}
				/_{(\coInd_{\alpha}V)(x)}.
				\]
				For each
				$x \in \Ob(\UC_{\alpha+1})$ with $d(x) = \alpha$, since $(\calC_x, \calD_x)$ is a complete cotorsion pair in $A_x^0 \Mod$, it follows from Lemma \ref{bri lem} that
				$(\Mon(\calC_x), \Epi(\calD_x))$ is a weak factorization system in $A_x^0 \Mod$. This induces a weak factorization system
				\begin{equation*}\label{ind wfs}\tag{wfs 1}
					\begin{gathered}
						( \
						_{(\Ind_{\alpha}V)(x)}
						\backslash{\substack{\Mon(\calC_x) \\ \tau^{\alpha}_V(x)}}
						/_{(\coInd_{\alpha}V)(x)},                                         \
						_{(\Ind_{\alpha}V)(x)}
						\backslash{\substack{\Epi(\calD_x) \\ \tau^{\alpha}_V(x)}}
						/_{(\coInd_{\alpha}V)(x)}
						\ )
					\end{gathered}
				\end{equation*}
				in the category
				${_{(\Ind_{\alpha}V)(x)}}
				\backslash{\substack{A^0_x\Mod \\ \tau^{\alpha}_V(x)}}
				/_{(\coInd_{\alpha}V)(x)}$; see Example \ref{wfs for over/under}. Thus, under the identification via $\scrR$, one gets the following weak factorization system in ${\UC_{\alpha+1}\Mod}_V$:
				\begin{equation*}\label{ind wfs all}\tag{wfs 2}
					\begin{gathered}
						( \
						\prod_{d(x) = \alpha}
						{_{(\Ind_{\alpha}V)(x)}
							\backslash{\substack{\Mon(\calC_x) \\ \tau^{\alpha}_V(x)}}
							/_{(\coInd_{\alpha}V)(x)}},                               \
						\prod_{d(x) = \alpha}
						{_{(\Ind_{\alpha}V)(x)}}
						\backslash{\substack{\Epi(\calD_x) \\ \tau^{\alpha}_V(x)}}
						/_{(\coInd_{\alpha}V)(x)}
						\ ).
					\end{gathered}
				\end{equation*}
				
				In what follows, in order to simplify the symbols,
				denote the weak factorization system (\ref{ind wfs}) by $(\sfC_V(x), \sfD_V(x))$.
				We write $(\sfC_V, \sfD_V)$ for the system in (\ref{ind wfs all}).
				
				Next, for any morphism $u : V \to W \in \UC_{\alpha}\Mod$, consider the adjunction
				$$u_! : {\UC_{\alpha+1}\Mod}_V \ \rightleftarrows \ {\UC_{\alpha+1}\Mod}_W : u^*.$$
We prove that $u_!$ preserves the left maps in (\ref{ind wfs all}). By the explicit description of the action of $u_!$ on morphisms (see Remark \ref{des of u(k,l)}), it suffices to prove that for any morphism $l \in {\UC_{\alpha+1}\Mod}_V$ and each $x \in \Ob(\UC_{\alpha+1})$ with $d(x) = \alpha$, if $l(x)$ is contained in $\Mon(\calC_x)$, then so is $(u_!(l))(x)$. But this is clear since (\ref{ind pushout}) is a pushout and $(\Mon(\calC_x), \Epi(\calD_x))$ is a weak factorization system.
				
				Finally, applying Lemma \ref{key lem for con wfs}, one gets a weak factorization system
				$(\sfC_{\UC_{\alpha+1}\Mod}, \sfD_{\UC_{\alpha+1}\Mod})$
				in $\UC_{\alpha+1}\Mod$, where
				\begin{align*}
				\sfC_{\UC_{\alpha+1}\Mod}
				&= \{\ f : Y \to Z \in \UC_{\alpha+1}\Mod \, \mid \,
				\Res_{\alpha}f \in \Mon(\Phi_{\alpha}(\calC)) \text{ and } f_{\triangleright} \in \sfC_{\Res_{\alpha}Z}\}, \ \mathrm{and} \\
				\sfD_{\UC_{\alpha+1}\Mod}
				&= \{\ f : Y \to Z \in \UC_{\alpha+1}\Mod \, \mid \,
				\Res_{\alpha}f \in \Epi(\Psi_{\alpha}(\calD)) \text{ and } f^{\triangleleft} \in \sfD_{\Res_{\alpha}Y}\ \}.
				\end{align*}
			\end{bfhpg*}
			
			\begin{bfhpg*}
				[Step 2: Prove $\sfC_{\UC_{\alpha+1}\Mod}=\Mon(\Phi_{\alpha+1}(\calC))$ and $\sfD_{\UC_{\alpha+1}\Mod}=\Epi(\Psi_{\alpha+1}(\calD))$] We prove both equalities, treating the latching and matching cases separately.

	We first prove the equality $\sfC_{\UC_{\alpha+1}\Mod}=\Mon(\Phi_{\alpha+1}(\calC))$.			Let $f : Y \to Z$ be a morphism in $\UC_{\alpha+1}\Mod$ and let $u=\Res_{\alpha}f: V \to W$. Then one has
				\begin{rqm}
					\item $f$ is contained in $\Mon(\Phi_{\alpha+1}(\calC))$ if and only if it satisfies the following conditions:
					\begin{eqc}
						\item $f$ is monic, that is, $f(x)$ is monic for each $x \in \Ob(\UC_{\alpha+1})$;
						
						\item $\coker(f) \in \Phi_{\alpha+1}(\calC)$, that is, $l^{d(x)}_{\coker(f)}(x)$ is monic with its cokernel contained in $\calC_x$ for each $x \in \Ob(\UC_{\alpha+1})$.
					\end{eqc}
					
					\item $f$ is contained in $\sfC_{\UC_{\alpha+1}\Mod}$ if and only if it satisfies the following conditions:
					\begin{prt}
						\item $u$ is monic, that is, $u(x)$ is monic for each $x \in \Ob(\UC_{\alpha})$.
						
						\item $\coker(u) \in \Phi_{\alpha}(\calC)$, that is, $l^{d(x)}_{\coker(u)}(x)$ is monic with its cokernel contained in $\calC_x$ for each $x \in \Ob(\UC_{\alpha})$.
						\item The fiber morphism $f_{\triangleright}$ is in $\sfC_{W}$, that is, up to the identification given by $\scrR$, the family
						$$\{(\Ind_{\alpha}W)(x) \sqcup_{(\Ind_{\alpha}V)(x)} Y(x)\overset{\sigma(x)}\longrightarrow Z(x)\}_{d(x) = \alpha}$$
						belongs to $\sfC_W$ (see Remark \ref{des of fiber mor}), where $\sigma(x)$ appears in the following commutative diagram with the inner square a pushout:
						\begin{equation}\label{ast}
							\begin{gathered}
								\xymatrix@R=1cm@C=0.5cm{
									(\Ind_{\alpha}V)(x)
									\ar[d]_-{l^{\alpha}_Y(x)}
									\ar[r]^-{(\Ind_{\alpha}u)(x)}
									\ar@{}[rd]^(0.6)>>{\lrcorner}
									& (\Ind_{\alpha}W)(x)
									\ar[d]_{\varphi(x)}
									\ar@/^1.5pc/[ddr]^{l^{\alpha}_Z(x)}                          \\
									Y(x)
									\ar[r]^-{\theta(x)}
									\ar@/_1.5pc/[drr]_{f(x)}
									& (\Ind_{\alpha}W)(x) \sqcup_{(\Ind_{\alpha}V)(x)} Y(x)
									\ar[dr]|-{\, \sigma(x) \,}                              \\
									&& Z(x).}
							\end{gathered}
						\end{equation}
						Thus, the condition (c) is equivalent to
						\item[(c')] $\sigma(x)$ is monic with its cokernel contained in $\calC_x$ for each $x \in \Ob(\UC_{\alpha+1})$ with $d(x) = \alpha$.
					\end{prt}
				\end{rqm}
				
				We mention that no matter which class $f$ belongs to, one has $\Res_{\alpha}f = u\in\Mon(\Phi_{\alpha}(\calC))$ by Proposition \ref{useful}, so one gets that
				$$0 \to V \overset{u} \longrightarrow W \to \coker(u) \to 0$$
				is an exact sequence in $\UC_{\alpha}\Mod$ with $\coker(u) \in \Phi_{\alpha}(\calC)$. By Lemma \ref{ind keep exact}, for each $x \in \Ob(\UC_{\alpha+1})$ with $d(x) = \alpha$, the sequence
				\begin{equation*}
					0\to(\Ind_{\alpha}V)(x)
					\xra{(\Ind_{\alpha}u)(x)}
					(\Ind_{\alpha}W)(x)
					\to(\Ind_{\alpha}\coker(u))(x)
					\to 0
				\end{equation*}
				is exact. Since the inner square in (\ref{ast}) is a pushout, one gets the next commutative diagram with exact rows:
				\begin{equation*}
					\begin{gathered}
						\xymatrix@R=1cm@C=0.5cm{
							0
							\ar[r]
							&     (\Ind_{\alpha}V)(x)
							\ar[r]^-{(\Ind_{\alpha}u)(x)}
							\ar[d]_-{l^{\alpha}_Y(x)}
							\ar@{}[rd]^(0.6)>>{\lrcorner}
							&     (\Ind_{\alpha}W)(x)
							\ar[r]^{}
							\ar[d]_-{\varphi(x)}
							&     (\Ind_{\alpha}\coker(u))(x)
							\ar[r]
							\ar@{=}[d]
							&              0                         \\
							0
							\ar[r]
							&             Y(x)
							\ar[r]^-{\theta(x)}
							& (\Ind_{\alpha}W)(x) \sqcup_{(\Ind_{\alpha}V)(x)} Y(x)
							\ar[r]^-{}
							&             (\Ind_{\alpha}\coker(u))(x)
							\ar[r]
							&              0
						}
					\end{gathered}
				\end{equation*}
				Then one has the following commutative diagram:
				\begin{equation}\label{dagger}
					\begin{gathered}
						\xymatrix@R=1cm@C=0.5cm{
							0
							\ar[r]
							&               (\Ind_{\alpha}V)(x)
							\ar[r]^-{(\Ind_{\alpha}u)(x)}
							\ar[d]_-{l^{\alpha}_Y(x)}
							\ar@{}[rd]^(0.6)>>{\lrcorner}
							&               (\Ind_{\alpha}W)(x)
							\ar[r]^{}
							\ar[d]_-{\varphi(x)}
							&           (\Ind_{\alpha}\coker(u))(x)
							\ar[r]
							\ar@{=}[d]
							&                       0                               \\
							0
							\ar[r]
							&                      Y(x)
							\ar[r]^-{\theta(x)}
							\ar@{=}[d]
							& (\Ind_{\alpha}W)(x) \sqcup_{(\Ind_{\alpha}V)(x)} Y(x)
							\ar[r]^-{}
							\ar[d]_-{\sigma(x)}
							&                 (\Ind_{\alpha}\coker(u))(x)
							\ar[r]
							\ar[d]_-{\kappa}
							&                        0                            \\
							&                       Y(x)
							\ar[r]^-{f(x)}
							&                       Z(x)
							\ar[r]^-{}
							&                    \coker(f(x))
							\ar[r]
							&                        0
						}
					\end{gathered}
				\end{equation}
				Since $\Res_{\alpha}f = u$, it is easy to see that
				$(\Ind_{\alpha}\coker(u))(x) = (\Ind_{\alpha}\Res_{\alpha}\coker(f))(x)$. But $\sigma(x) \circ \varphi(x) = l^{\alpha}_Z(x)$ by the commutative diagram (\ref{ast}). By the universal property of cokernels, one can check that $\kappa = l^{\alpha}_{\coker(f)}(x)$. Applying the Snake Lemma to the lower part of (\ref{dagger}), one gets that
				\begin{equation}\label{9.2.1}
					\coker(\sigma(x)) \cong\coker(\kappa) = \coker(l^{\alpha}_{\coker(f)}(x)),\ \mathrm{and}
				\end{equation}
				\begin{equation}\label{9.2.2}
					\sigma(x)\ \mathrm{is\ monic\ if\ and\ only\ if}\ \kappa\ \mathrm{is\ monic\ if\ and\ only\ if}\ l^{\alpha}_{\coker(f)}(x)\ \mathrm{is\ monic}.
				\end{equation}
				
				Now we prove that $\Mon(\Phi_{\alpha+1}(\calC))=\sfC_{\UC_{\alpha+1}\Mod}$. Let $f$ be in $\Mon(\Phi_{\alpha+1}(\calC))$. Then by (\ref{9.2.1}) and (\ref{9.2.2}) one gets that the condition (c') holds. It is easy to see that the conditions (a) and (b) hold, so one gets that $f\in\sfC_{\UC_{\alpha+1}\Mod}$. Conversely, let $f$ be in $\sfC_{\UC_{\alpha+1}\Mod}$, that is, the conditions (a), (b) and (c') hold. Then by (\ref{9.2.1}) and (\ref{9.2.2}) one gets that the condition $(ii)$ holds. On the other hand, since $f(x) = \sigma(x) \circ \theta(x)$ and $\theta(x)$ is monic by the diagram (\ref{dagger}), one gets that the condition $(i)$ holds as $\sigma(x)$ is monic by (c').
			\end{bfhpg*}

We then prove that
$\sfD_{\UC_{\alpha+1}\Mod}
=\Epi(\Psi_{\alpha+1}(\calD))$.
Let $f:Y\to Z$ and let
$u=\Res_\alpha f:V\to W$.
For each $x$ with $d(x)=\alpha$,
the cartesian factorization identifies
$f^{\triangleleft}(x)$ with
\[
\rho(x):Y(x)\to
P_x:=Z(x)\times_{M_x(W)}M_x(V).
\]
Suppose that $u\in\Epi(\Psi_\alpha(\calD))$,
and let $K=\ker f$ and $L=\ker u$.
Then $L=\Res_\alpha K\in\Psi_\alpha(\calD)$.
By Lemma~\ref{matching-exact-on-psi},
\[
0\to M_x(L)\to M_x(V)\to M_x(W)\to0
\]
is exact. Hence the projection
$P_x\to Z(x)$ is epic with kernel $M_x(L)$.
By the definitions of the matching maps
and $\rho(x)$, one has
$\ker\rho(x)\cong\ker m_K^\alpha(x)$,
and there is an exact sequence
\[
0\to\coker(m_K^\alpha(x))
\to\coker(\rho(x))
\to\coker(f(x))\to0.
\]
Thus one gets that
$\rho(x)\in\Epi(\calD_x)$ if and only if
$f(x)$ and $m_K^\alpha(x)$ are epic and
$\ker m_K^\alpha(x)\in\calD_x$.
Since $u\in\Epi(\Psi_\alpha(\calD))$,
these conditions hold at all objects of degree
less than $\alpha$. The preceding equivalence
shows that they also hold at degree $\alpha$.
Thus one has
$\sfD_{\UC_{\alpha+1}\Mod}
=\Epi(\Psi_{\alpha+1}(\calD))$.
\end{bfhpg*}

\begin{bfhpg*}[The case that $\beta$ is a limit ordinal]
Since
$\UC_\beta=\bigcup_{\alpha<\beta}\UC_\alpha$,
compatible families of modules and morphisms
over the $\UC_\alpha$ glue uniquely to
modules and morphisms over $\UC_\beta$. Moreover, the latching and
matching maps at each object are determined
by its lower-degree subcategory.
By the construction in
Lemma~\ref{key lem for con wfs} and
Remark~\ref{compatibal-wfs}, the factorizations
and lifting morphisms can be chosen
compatibly at successive stages.
At limit stages, these compatible
morphisms glue uniquely. Since the classes
$\Mon(\Phi_\alpha(\calC))$ and
$\Epi(\Psi_\alpha(\calD))$ are detected
by restriction, the resulting factorization
and lifting properties give a weak
factorization system
\[
\bigl(\Mon(\Phi_\beta(\calC)),
      \Epi(\Psi_\beta(\calD))\bigr)
\]
in $\UC_\beta\Mod$, where the lifting-class
equalities follow from the usual retract
argument. The conclusion then follows from
Lemma \ref{bri lem}.
\qed
\end{bfhpg*}
\end{bfhpg}

	The next result shows that the induced complete cotorsion pair in $\UC \Mod$ inherits the hereditary property.
	
	\begin{proposition}\label{hereditary to here}
		Suppose that $\UC^+$ is a projective right $\UC^0$-module and $\UC^-$ is a projective left $\UC^0$-module. If $(\calC, \calD) = \{(\calC_x, \calD_x)\}_{x \in \Ob(\UC)}$ is a family of hereditary complete cotorsion pairs in those $A^0_x\Mod$'s, then $(\Phi(\calC), \Psi(\calD))$ forms a hereditary complete cotorsion pair in $\UC \Mod$.
	\end{proposition}
	
	\begin{prf*}
		Theorem \ref{ma re for ccp} asserts that $(\Phi(\calC), \Psi(\calD))$ is a complete cotorsion pair in $\UC \Mod$, so we only need to deal with the hereditary property. By Lemma \ref{test hereditary}, it suffices to show that it is resolving; that is, $\Phi(\calC)$ is closed under kernels of epimorphisms in $\UC \Mod$. Explicitly, given a short exact sequence
		\[
		0 \to Y \to Z \to T \to 0
		\]
		in $\UC \Mod$ with both $Z$ and $T$ in $\Phi(\calC)$, we need to show that
		$Y \in \Phi(\calC)$ as well.
		
		For each $x \in \Ob(\UC)$ with $d(x) = \alpha$, one gets a short exact sequence
		\[
		0 \to \Res_{\alpha}Y \to \Res_{\alpha}Z \to \Res_{\alpha}T \to 0
		\]
		in $\UC_{\alpha} \Mod$. But $\Res_{\alpha}T \in \Phi_{\alpha}(\calC)$ by Proposition \ref{useful} as $T \in \Phi(\calC)$, and so by Lemma \ref{ind keep exact}, it induces an exact sequence
		\[
		0 \to (\Ind_{\alpha}\Res_{\alpha}Y)(x)
		\to (\Ind_{\alpha}\Res_{\alpha}Z)(x)
		\to (\Ind_{\alpha}\Res_{\alpha}T)(x) \to 0
		\]
		in $A^0_x\Mod$, which fits into the following commutative diagram with exact rows:
		\[
		\xymatrix@R=1cm@C=0.5cm{
			0
			\ar[r]
			&     (\Ind_{\alpha}\Res_{\alpha}Y)(x)
			\ar[r]^-{}
			\ar[d]_-{l^{\alpha}_Y(x)}
			&     (\Ind_{\alpha}\Res_{\alpha}Z)(x)
			\ar[r]^{}
			\ar[d]_-{l^{\alpha}_Z(x)}
			&     (\Ind_{\alpha}\Res_{\alpha}T)(x)
			\ar[r]
			\ar[d]_-{l^{\alpha}_T(x)}
			&              0                                \\
			0
			\ar[r]
			&             Y(x)
			\ar[r]^-{}
			&             Z(x)
			\ar[r]^-{}
			&             T(X)
			\ar[r]
			&              0
		}
		\]
		Since $Z$ and $T$ belong to $\Phi(\calC)$, it follows that $l^{\alpha}_Z(x)$  and $l^{\alpha}_T(x)$ are monic, and both $\coker(l^{\alpha}_Z(x))$ and $\coker(l^{\alpha}_T(x))$ are in $\calC_x$. By the Snake Lemma one gets that $l^{\alpha}_Y(x)$ is monic, and obtain an exact sequence
		\[
		0 \to \coker(l^{\alpha}_Y(x))
		\to \coker(l^{\alpha}_Z(x))
		\to \coker(l^{\alpha}_T(x)) \to 0.
		\]
		of left $A_x^0$-modules. But $\calC_x$ is resolving as the cotorsion pair $(\calC_x, \calD_x)$ is hereditary by assumption, so one has $\coker(l^{\alpha}_Y(x)) \in \calC_x$, and hence $Y \in \Phi(\calC)$ as desired.
	\end{prf*}
	
	In the rest of the section, we consider two special types of generalized $k$-linear Reedy categories.
	
	\begin{definition}
		Let $\UC$ be a generalized $k$-linear Reedy category.
		\begin{rqm}
			\item
			If $\UC^+ = \UC$, then $\UC$ is called a \textit{generalized} $k$-\textit{linear direct category}; in this case one has $\UC^-=\UC^0$.
			
			\item
			If $\UC^- = \UC$, then $\UC$ is called a \textit{generalized} $k$-\textit{linear inverse category}; in this case one has $\UC^+=\UC^0$.
		\end{rqm}
	\end{definition}

Given a family
	$\calS = \{ \calS_x \}_{x \in \Ob(\UC)}$ of classes of objects in those $A_x^0 \Mod$'s, denote by $\UC\Mod_\calS$ the subcategory of $\UC \Mod$ consisting of $\UC$-modules $Y$ such that $Y(x)$ belongs to $\calS_x$ for each $x \in \Ob(\UC)$.
	
	\begin{lemma}\label{fact on direct cat}
		If $\UC$ is a generalized $k$-linear direct category, then for any family $\calS = \{ \calS_x \}_{x \in \Ob(\UC)}$ of classes of objects in those $A_x^0 \Mod$'s, there is an equality $\Psi(\calS) = \UC\Mod_\calS$.
\end{lemma}

\begin{prf*}
For each $\UC$-module $Y$ and any $x \in \Ob(\UC)$ with $d(x) = \alpha$, one has
		\begin{align*}
			(\coInd_{\alpha}\Res_{\alpha}Y)(x)
			& =      \Hom_{\UC_{\alpha}}(\UC 1_x, \Res_{\alpha}Y)            \\
			& \cong  \Hom_{\UC^-_{\alpha}}(\UC^- 1_x, \Res_{\alpha}Y)         \\
			& =      \Hom_{\UC^0_{\alpha}}(\UC^0 1_x, \Res_{\alpha}Y)         \\
			& =  0,
		\end{align*}
		where the first isomorphism holds by Theorem \ref{Cofinality}(b), the second equality holds as $\UC^-=\UC^0$, and the last equality follows from Lemma \ref{Reedy factorization}(a). Consequently, one has $\Psi(\calS) = \UC\Mod_\calS$.
	\end{prf*}

Dually, we have the next result.

	\begin{lemma}\label{fact on inverse cat}
		If $\UC$ is a generalized $k$-linear inverse category, then for any family $\calS = \{ \calS_x \}_{x \in \Ob(\UC)}$ of classes of objects in those $A_x^0 \Mod$'s, there is an equality $\Phi(\calS) = \UC\Mod_\calS$.
\end{lemma}
	
	In what follows, for any category $\calE$, denote by $\calP(\calE)$ (resp., $\calI(\calE)$) the class of all projective (resp., injective) objects in $\calE$. The following result gives a characterization of projective objects in $\UC\Mod$ when $\UC$ is a generalized $k$-linear direct category.
	
	\begin{corollary}\label{pro obj for direct}
		Let $\UC$ be a generalized $k$-linear direct category such that it is a projective right $\UC^0$-module. Then $\calP(\UC \Mod) = \Phi(\{\calP_x\}_{x\in\Ob(\UC)})$, where $\calP_x$ is the class of all projective left $A_x^0$-modules.
	\end{corollary}
\begin{prf*}
It is known that $(\calP_x, A_x^0 \Mod)$ is a complete cotorsion pair for each $x\in\Ob(\UC)$, and so it follows from Theorem \ref{ma re for ccp} and Lemma \ref{fact on direct cat} that $(\Phi(\{\calP_x\}_{x\in\Ob(\UC)}), \UC\Mod)$ is a complete cotorsion pair in $\UC\Mod$. Thus one has $\calP(\UC \Mod) = \Phi(\{\calP_x\}_{x\in\Ob(\UC)})$.
\end{prf*}
	
	A dual version of the above result gives a characterization of injective objects in $\UC\Mod$ when $\UC$ is a generalized $k$-linear inverse category.
	
	\begin{corollary}\label{inj obj for inverse}
		Let $\UC$ be a generalized $k$-linear inverse category such that it is a projective left $\UC^0$-module. Then $\calI(\UC \Mod) = \Psi(\{\calI_x\}_{x\in\Ob(\UC)})$, where $\calI_x$ is the class of all injective left $A_x^0$-modules.
	\end{corollary}
	
	
	\section{Lifting of abelian model structures}\label{lift HHT and Exa}
	\noindent
	In this section, as an application of Theorem \ref{ma re for ccp}, we show that under the extra (co)compatible condition, a family
	$\{(\calQ_x, \calW_x, \calR_x) \}_{x \in \Ob(\UC)}$ of (hereditary) Hovey triples in those $A_x^0 \Mod$'s induces a (hereditary) Hovey triple in $\UC \Mod$. We then give some examples in Gorenstein homological algebra.
	
	\begin{definition}\label{compatible}
		Let $\calS = \{ \calS_x \}_{x \in \Ob(\UC)}$ be a family with each $\calS_x$ a class of objects in $A_x^0 \Mod$. Then $\calS$ is called
		\begin{rqm}
			\item
			\emph{cocompatible} if $S \in \calS_y$ implies that $\UC^+(y, x) \otimes_{A_y^0} S \in \calS_x$ for all $x, y \in \Ob(\UC)$;
			
			\item
			\emph{compatible} if $S \in \calS_y$ implies that $\Hom_{A_y^0}(\UC^-(x, y), S) \in \calS_x$ for all $x, y \in \Ob(\UC)$.
		\end{rqm}
	\end{definition}
	
	\begin{remark}
		If $\UC$ is a $k$-linear Reedy category in the sense of \cite{GS}, then $A_x^0\cong k$ is a field for each $x \in \Ob(\UC)$. In this case, $\UC^+(y, x)$ and $\UC^-(x, y)$ are free $k$-modules for all $x, y \in \Ob(\UC)$. Thus, if $\calS$ is closed under coproducts then it is cocompatible, and if $\calS$ is closed under products then it is compatible.
	\end{remark}
	
	\begin{example}\label{exa-compatible}
		Let $\calP = \{ \calP_x \}_{x \in \Ob(\UC)}$ be a family with each $\calP_x$ the class of all projective objects in $A_x^0 \Mod$, and let  $\calF = \{ \calF_x \}_{x \in \Ob(\UC)}$ be a family with each $\calF_x$ the class of all flat objects in $A_x^0 \Mod$. If $\UC^+$ is a projective left $\UC^0$-module, then both $\calP$ and $\calF$ are cocompatible. Dually, let $\calI = \{ \calI_x \}_{x \in \Ob(\UC)}$ be a family with each $\calI_x$ the class of all injective objects in $A_x^0 \Mod$. If $\UC^-$ is a projective right $\UC^0$-module, then $\calI$ is compatible.
	\end{example}

	\begin{example}\label{dirct co-compatible}
		If $\UC$ is a generalized $k$-linear direct category, then any family $\calS = \{ \calS_x \}_{x \in \Ob(\UC)}$ of classes of objects in those $A_x^0 \Mod$'s is compatible. Indeed, for all $x, y \in \Ob(\UC)$ and $S \in \calS_y$,
		\begin{prt}
			\item[$\bullet$] if $x \neq y$, then
			$\Hom_{A_y^0}(\UC^-(x, y), S)
			= \Hom_{A_y^0}(\UC^0(x, y), S)
			= 0 \in \calS_x$;
			
			\item[$\bullet$] if $x = y$, then $\Hom_{A_y^0}(\UC^-(x, x), S)
			= \Hom_{A_x^0}(\UC^0(x, x), S)
			\cong S \in \calS_x$.
		\end{prt}
		
		Dually, if $\UC$ is a generalized $k$-linear inverse category, then any family $\calS = \{ \calS_x \}_{x \in \Ob(\UC)}$ of classes of objects in those $A_x^0 \Mod$'s is cocompatible.
	\end{example}
	
	Now we state the main result in this section.
	
	\begin{theorem}\label{HT TO HT}
		Suppose that $\UC^+$ is a projective right $\UC^0$-module and $\UC^-$ is a projective left $\UC^0$-module. Let $(\calQ, \calW, \calR) = \{(\calQ_x, \calW_x, \calR_x)\}_{x \in \Ob(\UC)}$ be a family with each $(\calQ_x, \calW_x, \calR_x)$ a Hovey triple in $A^0_x\Mod$. If $\calQ \cap \calW = \{ \calQ_x \cap \calW_x \}_{x \in \Ob(\UC)}$ is cocompatible and $\calW \cap \calR = \{ \calW_x \cap \calR_x \}_{x \in \Ob(\UC)}$ is compatible, then
		\[
		( \Phi(\calQ), \, \UC\Mod_\calW, \, \Psi(\calR) )
		\]
		forms a Hovey triple in $\UC \Mod$.
	\end{theorem}
	
	Before giving a proof, we need some auxiliary facts and results. The next result can be proved similarly as in Lemma \ref{filtration1}.
	
	\begin{lemma}\label{filtration2}
		Suppose that $\UC^+$ is a projective right $\UC^0$-module and $\UC^-$ is a projective left $\UC^0$-module. For each $x \in \Ob(\UC)$, the right $\UC$-module $1_x \UC$ admits a filtration of right $\UC$-modules with factors
		$$\bigoplus_{d(z) = \gamma} \UC^+(z, x) \otimes_{A^0_z} \Delta^z$$ for $0\leqslant\gamma\leqslant\lambda$, where $\Delta^z$ is the standard right $\UC$-module.
	\end{lemma}
	
	\begin{lemma}\label{each entry S}
		Suppose that $\UC^+$ is a projective right $\UC^0$-module and $\UC^-$ is a projective left $\UC^0$-module. Let $\calS = \{ \calS_x \}_{x \in \Ob(\UC)}$ be a cocompatible family of classes of objects in those $A_x^0 \Mod$'s such that $\calS_x$ is closed under filtrations and coproducts for each $x \in \Ob(\UC)$. Then one has $\Phi(\calS)\subseteq\UC\Mod_\calS$.
	\end{lemma}
	
	\begin{prf*}
	Let $Y$ be in $\Phi(\calS)$ and $x \in \Ob(\UC)$. Since $Y(x) \cong 1_x \UC \otimes_{\UC} Y$ and $\calS_x$ is closed under filtrations in $A_x^0 \Mod$, it suffices to show that $1_x \UC \otimes_{\UC} Y$ admits a filtration whose factors are in $\calS_x$. By Lemma \ref{filtration2}, the right $\UC$-module $1_x \UC$ has a filtration of right $\UC$-modules with factors $\bigoplus_{d(z) = \gamma} \UC^+(z, x) \otimes_{A^0_z} \Delta^z$ for $0\leqslant\gamma\leqslant\lambda$. Since $Y$ is contained in $\Phi(\calS)$, one has $\Tor^{\UC}_1 ( \Delta^z, Y ) = 0$. Hence,
		$\Tor^{\UC}_1 (\UC^+(z, x) \otimes_{A_z^0} \Delta^z, Y) = 0$ as $\UC^+(z, x)$ is a projective right $A_z^0$-module by Lemma \ref{projectivity}(a). It follows that
		\[
		\Tor^{\UC}_1 (\bigoplus_{d(z) = \gamma}\UC^+(z, x) \otimes_{A_z^0} \Delta^z, Y) = 0.
		\]
		Consequently, $1_x \UC \otimes_{\UC} Y$ admits a filtration with factors
		$(\bigoplus_{d(z) = \gamma} \UC^+(z, x) \otimes_{A^0_z} \Delta^z) \otimes_{\UC} Y$ for $0\leqslant\gamma\leqslant\lambda$.
		
		It remains to show that $( \UC^+(z, x) \otimes_{A_z^0} \Delta^z ) \otimes_{\UC} Y \in \calS_x$ as $\calS_x$ is closed under coproducts in $A_x^0\Mod$. This is clear. Indeed, by the definition of $\Phi(\calS)$, we know that
		$\Delta^z \otimes_{\UC} Y \in \calS_z$. Therefore, by the cocompatible assumption,
		$(\UC^+(z, x) \otimes_{A_z^0} \Delta^z) \otimes_{\UC} Y \cong \UC^+(z, x) \otimes_{A_z^0} (\Delta^z \otimes_{\UC} Y)$
		is in $\calS_x$, as desired.
	\end{prf*}
	
	The following result plays a key role in the proof of Theorem \ref{HT TO HT}.
	
	\begin{lemma}\label{tilQ = Q cap W}
		Suppose that $\UC^+$ is a projective right $\UC^0$-module and $\UC^-$ is a projective left $\UC^0$-module. Let $\calS = \{ \calS_x \}_{x \in \Ob(\UC)}$ and $\calW = \{ \calW_x \}_{x \in \Ob(\UC)}$ be families with both $\calS_x$ and $\calW_x$ classes of objects in $A_x^0 \Mod$. If the family $\calS \cap \calW = \{ \calS_x \cap \calW_x \}_{x \in \Ob(\UC)}$ is cocompatible, and for each $x\in\Ob(\UC)$ the following conditions hold:
		\begin{prt}
			\item
			$\calS_x \cap \calW_x$ is closed under filtrations and coproducts in $A_x^0 \Mod$, and
			
			\item
			$\calW_x$ is closed under cokernels of monomorphisms in $A_x^0 \Mod$,
		\end{prt}
		then one has $\Phi(\calS \cap \calW) = \Phi(\calS) \cap \UC\Mod_\calW$.
	\end{lemma}
	
	\begin{prf*}
		It is clear that $\Phi(\calS \cap \calW) \subseteq \Phi(\calS)$. On the other hand, by Lemma \ref{each entry S}, we have
		$\Phi(\calS \cap \calW) \subseteq \UC\Mod_{\calS \cap \calW} \subseteq \UC\Mod_{\calW}$.
		Thus, $\Phi(\calS \cap \calW) \subseteq \Phi(\calS) \cap \UC\Mod_\calW$.
		
		For the other inclusion, take a left $\UC$-module $Y \in \Phi(\calS) \cap \UC\Mod_\calW$. For each $x \in \Ob(\UC)$ with $d(x) = \alpha$, one gets that $l^{\alpha}_Y(x)$ is a monomorphism with $\Delta^x \otimes_{\UC} Y\in\calS_x$. Hence, to show $Y \in \Phi(\calS \cap \calW)$, it suffices to prove
		$\Delta^x \otimes_{\UC} Y \in \calW_x$. We check this fact by transfinite induction on $d(x) = \alpha$.
		
		If $\alpha = 0$, then $\Delta^x \cong 1_x \UC$, and so $\Delta^x \otimes_{\UC} Y \cong Y(x)\in \calW_x$ as $Y \in \UC\Mod_\calW$.
		
		Suppose that $\alpha > 0$ and $\Delta^z \otimes_{\UC} Y \in \calW_z$ for all $z \in \Ob(\UC)$ with $d(z) < \alpha$. Since $Y \in \Phi(\calS)$, one has $\Tor^{\UC}_1 (\Delta^x, Y) = 0$. So we have a short exact sequence
		\[
		0 \to 1_x \mathfrak{I}_{\alpha} \otimes_{\UC} Y
		\to 1_x \UC \otimes_{\UC} Y
		\to \Delta^x \otimes_{\UC} Y
		\to 0
		\]
		in $A_x^0 \Mod$. Note that $1_x \UC \otimes_{\UC} Y \cong Y(x) \in \calW_x$ as $Y \in \UC\Mod_\calW$ and  $\calW_x$ is closed under cokernels of monomorphisms in
		$A_x^0 \Mod$ by assumption. To check $\Delta^x \otimes_{\UC} Y \in \calW_x$, one has to prove $1_x \mathfrak{I}_{\alpha} \otimes_{\UC} Y \in \calW_x$.
		
		By Lemma \ref{key for iso}, one has
		\[
		1_x \mathfrak{I}_{\alpha} \otimes_{\UC} Y
		\cong  (\Ind_{\alpha}\Res_{\alpha}Y)(x)
		=    1_x \UC \otimes_{\UC_{\alpha}} \Res_{\alpha}Y.
		\]
	Note that $\Res_{\alpha}Y \in \Phi_{\alpha}(\calS)$ by Proposition \ref{useful} as $Y \in \Phi(\calS)$. Using the filtration of right $\UC_{\alpha}$-modules given in Remark \ref{filtration1} as well as the argument in the proof of the previous lemma, we deduce that $1_x \UC \otimes_{\UC_{\alpha}} \Res_{\alpha}Y$ admits a filtration with factors
		\[
		(\bigoplus_{d(z) = \gamma} \UC^+(z, x) \otimes_{A_z^0} \Delta^z ) \otimes_{\UC_{\alpha}} \Res_{\alpha}Y
		\]
		for $0 \leqslant \gamma < \alpha$. But $\calS_x \cap \calW_x$ is closed under filtrations and coproducts in $A_x^0 \Mod$, it suffices to show that
		\[
		(\UC^+(z, x) \otimes_{A_z^0} \Delta^z) \otimes_{\UC_{\alpha}} \Res_{\alpha}Y \in \calS_x \cap \calW_x
		\]
		for each $z$ with $d(z) < \alpha$.
		
		By the induction hypothesis, $\Delta^z \otimes_{\UC_{\alpha}} \Res_{\alpha}Y \cong \Delta^z \otimes_{\UC} Y$ is contained in $\calW_z$, and it is also contained in $\calS_z$ as $Y \in \Phi(\calS)$. Thus, $\Delta^z \otimes_{\UC_{\alpha}} \Res_{\alpha}Y \in \calS_z \cap \calW_z$. By the cocompatible assumption on $\calS \cap \calW$, one gets that
		\[
		( \UC^+(z, x) \otimes_{A_z^0} \Delta^z) \otimes_{\UC_{\alpha}} \Res_{\alpha}Y
		\cong
		\UC^+(z, x) \otimes_{A_z^0} (\Delta^z \otimes_{\UC_{\alpha}} \Res_{\alpha}Y )
		\]
		is contained in $\calS_x \cap \calW_x$, as desired.
	\end{prf*}	

\begin{lemma}\label{dual-tilQ-equality}
Suppose that $\UC^+$ is a projective right
$\UC^0$-module and $\UC^-$ is a projective left
$\UC^0$-module. Let
$\calS=\{\calS_x\}_{x\in\Ob(\UC)}$ and
$\calW=\{\calW_x\}_{x\in\Ob(\UC)}$
be families with both $\calS_x$ and $\calW_x$
classes of objects in $A_x^0\Mod$.
If the family
$\calS\cap\calW=
\{\calS_x\cap\calW_x\}_{x\in\Ob(\UC)}$
is compatible, and for each $x\in\Ob(\UC)$
the following conditions hold:
\begin{prt}
    \item
    $\calS_x\cap\calW_x$ is the right-hand
    class of a complete cotorsion pair
    in $A_x^0\Mod$, and

    \item
    $\calW_x$ is closed under kernels of
    epimorphisms in $A_x^0\Mod$,
\end{prt}
then one has
$\Psi(\calS\cap\calW)
=\Psi(\calS)\cap\UC\Mod_\calW$.
\end{lemma}

\begin{prf*}
Let $\calT=\calS\cap\calW$, and let
$H_x(Y)=\Hom_\UC(\Delta_x,Y)$.
We first observe that, for $Y\in\Psi(\calS)$
and $d(x)=\alpha$, the matching object
$M_x(\Res_\alpha Y)$ admits a cofiltration
whose successive kernels are products of
modules of the form
$\Hom_{A_z^0}(\UC^-(x,z),H_z(Y))$ with
$d(z)<\alpha$.
Indeed, this follows by applying
$\Hom_{\UC_\alpha}(-,\Res_\alpha Y)$
to the filtration of the restricted
representable used in
Lemma~\ref{matching-exact-on-psi}.
The transition maps are epimorphisms
as $\Ext^1_\UC(\Delta_z,Y)=0$.

Since $\calT$ is compatible and each
$\calT_x$ is the right-hand class of a
cotorsion pair, the dual Eklof lemma
shows that if $H_z(Y)\in\calT_z$ for $d(z)<\alpha$,
then one gets that
$M_x(\Res_\alpha Y)\in\calT_x$.
Now let $Y\in\Psi(\calT)$.
For every $x$, the exact sequence
\[
0\to H_x(Y)\to Y(x)
\to M_x(\Res_\alpha Y)\to0
\]
and the preceding observation imply
$Y(x)\in\calT_x$ as $\calT_x$
is closed under extensions.
Thus one has
$$\Psi(\calT)\subseteq
\Psi(\calS)\cap\UC\Mod_\calW.$$

Conversely, let
$Y\in\Psi(\calS)\cap\UC\Mod_\calW$.
We prove $H_x(Y)\in\calT_x$
by induction on $d(x)=\alpha$.
If this holds at all lower degrees,
the preceding observation gives
$M_x(\Res_\alpha Y)\in\calT_x
\subseteq\calW_x$.
Since $Y(x)\in\calW_x$ and the
matching map is epic, the closure
of $\calW_x$ under kernels of
epimorphisms yields
$H_x(Y)\in\calW_x$.
But $H_x(Y)\in\calS_x$, so
$H_x(Y)\in\calT_x$.
Thus $Y\in\Psi(\calT)$.
\end{prf*}

	Now we give a proof of Theorem \ref{HT TO HT}.
	
	\begin{bfhpg}[Proof of Theorem \ref{HT TO HT}]
		The given family $\{(\calQ_x, \calW_x, \calR_x)\}_{x \in \Ob(\UC)}$ of Hovey triples in those $A^0_x\Mod$'s gives rise to two families
		\[
		\{(\calQ_x \cap \calW_x, \calR_x)\}_{x \in \Ob(\UC)}
		\quad \text{ and } \quad
		\{(\calQ_x, \calW_x \cap \calR_x)\}_{x \in \Ob(\UC)}
		\]
		of complete cotorsion pairs in them. By Theorem \ref{ma re for ccp}, one obtains two complete cotorsion pairs
		\[
		( \Phi(\calQ \cap \calW), \Psi(\calR) )
		\quad \text{ and } \quad
		( \Phi(\calQ), \Psi(\calW \cap \calR) )
		\]
		in $\UC \Mod$. Note that $\calQ \cap \calW$ is cocompatible by assumption, $\calQ_x \cap \calW_x$ as the left part of a complete cotorsion pair
		is closed under filtrations and coproducts in $A_x^0 \Mod$ for each $x \in \Ob(\UC)$, and $\calW_x$ is closed under cokernels of monomorphisms in $A_x^0 \Mod$. It follows that $\Phi(\calQ \cap \calW) = \Phi(\calQ) \cap \UC\Mod_\calW$ by Lemma \ref{tilQ = Q cap W}, and $\Psi(\calW \cap \calR) = \UC\Mod_\calW \cap \Psi(\calR)$ by Lemma~\ref{dual-tilQ-equality}. Consequently, one gets two complete cotorsion pairs
		\[
		( \Phi(\calQ) \cap \UC\Mod_\calW, \Psi(\calR) )
		\quad \text{ and } \quad
		( \Phi(\calQ), \UC\Mod_\calW \cap \Psi(\calR) )
		\]
		in $\UC \Mod$. Since $\UC\Mod_\calW$ is clearly a thick subcategory of $\UC \Mod$, the conclusion follows from Hovey's correspondence; see Theorem \ref{Hovey cor}. \qed
	\end{bfhpg}
	
	As an immediate consequence of Theorem \ref{HT TO HT} and Lemma \ref{hereditary to here}, we obtain the following result.
	
	\begin{proposition}\label{hHT TO hHT}
		Under the assumptions of Theorem \ref{HT TO HT}, if each Hovey triple $(\calQ_x, \calW_x, \calR_x)$ in $A^0_x\Mod$ is hereditary, then so is the Hovey triple $(\Phi(\calQ), \UC\Mod_\calW, \Psi(\calR))$.
	\end{proposition}
	
	\begin{corollary}\label{direct hHT}
		Let $\UC$ be a generalized $k$-linear direct category such that it is a projective right $\UC^0$-module, and let
		$(\calQ, \calW, \calR) = \{(\calQ_x, \calW_x, \calR_x)\}_{x \in \Ob(\UC)}$
		be a family with each $(\calQ_x, \calW_x, \calR_x)$ a (hereditary) Hovey triple in $A^0_x\Mod$. If $\{ \calQ_x \cap \calW_x \}_{x \in \Ob(\UC)}$ is cocompatible, then
		$(\Phi(\calQ), \, \UC\Mod_\calW, \, \UC\Mod_\calR)$
		forms a (hereditary) Hovey triple in $\UC \Mod$.
	\end{corollary}
	\begin{prf*}
		It follows from Lemma \ref{fact on direct cat} that $\Psi(\calR)=\UC\Mod_\calR$. Combining Theorem \ref{HT TO HT}, Proposition \ref{hHT TO hHT} and Example \ref{dirct co-compatible}, one gets that $(\Phi(\calQ), \UC\Mod_\calW, \UC\Mod_\calR)$
		forms a (hereditary) Hovey triple in $\UC \Mod$.
	\end{prf*}
	
	A dual result is:
	
	\begin{corollary}\label{inverse hHT}
		Let $\UC$ be a generalized $k$-linear inverse category such that it is a projective left $\UC^0$-module, and let
		$(\calQ, \calW, \calR) = \{(\calQ_x, \calW_x, \calR_x)\}_{x \in \Ob(\UC)}$
		be a family with each $(\calQ_x, \calW_x, \calR_x)$ a (hereditary) Hovey triple in $A^0_x\Mod$. If $\{ \calW_x \cap \calR_x \}_{x \in \Ob(\UC)}$ is compatible, then
		$( \UC\Mod_\calQ, \, \UC\Mod_\calW, \, \Psi(\calR) )$
		forms a (hereditary) Hovey triple in $\UC\Mod$.
	\end{corollary}
	
	We end this section with some examples of Theorem \ref{HT TO HT}; we obtain several abelian model structures in $\UC \Mod$ related to Gorenstein homological modules.
	
	\begin{bfhpg}[Gorenstein homological modules]
		Recall from Enochs and Jenda \cite{rha} that a left $\UC$-module $M$ is called \emph{Gorenstein injective} if there is an exact sequence
		\[
		\cdots \to E_{1} \to E_0 \to E_{-1} \to \cdots
		\]
		of injective left $\UC$-modules such that it remains exact after applying the functor $\Hom_{\UC}(I,-)$ for every injective left $\UC$-module $I$, and $M \cong \coker{(E_1 \to E_0)}$. A left $\UC$-module $M$ is called \emph{Gorenstein flat} if there is an exact sequence
		\[
		\cdots \to F_{1} \to F_0 \to F_{-1} \to \cdots
		\]
		of flat left $\UC$-modules such that it remains exact after applying the functor $I\otimes_{\UC}-$ for every injective right $\UC$-module $I$, and $M \cong \coker{(F_1 \to F_0)}$. Recently, \v{S}aroch and \v{S}t'ov\'{\i}\v{c}ek introduced \emph{projectively coresolved Gorenstein flat} left $\UC$-modules in \cite{SS20} via replacing flat modules in the above exact sequence by projective modules.
	\end{bfhpg}
	
	Let $\calGI = \{ \calGI_x \}_{x \in \Ob(\UC)}$ (resp., $\calGF = \{ \calGF_x \}_{x \in \Ob(\UC)}$, $\calPGF = \{ \calPGF_x \}_{x \in \Ob(\UC)}$) be  families with $\calGI_x$ (resp., $\calGF_x$, $\calPGF_x$) the class of all Gorenstein injective (resp., Gorenstein flat, projectively coresolved Gorenstein flat) modules in $A_x^0 \Mod$.
	
	Recall that a left $A_x^0$-module $N$ is called \textit{cotorsion} if $\Ext_{A_x^0}^1(F, N)=0$ for each flat left $A_x^0$-module $F$. We let $\calCot = \{ \calCot_x \}_{x \in \Ob(\UC)}$ be the family with $\calCot_x$ the class of all cotorsion modules in $A_x^0 \Mod$. The next result can be found in \cite{SS20}.
	
	\begin{proposition}\label{SS20}
		For all $x\in \Ob(\UC)$, the triples
		\[
		(\calPGF_x, {{\calPGF_{x}}^\perp}, A_x^0 \Mod), \quad (\calGF_x, {{\calPGF_{x}}^\perp}, {\calCot_x}), \quad (A_x^0 \Mod, {^\perp\calGI_x}, \calGI_x)
		\]
		are hereditary Hovey triples in  $A_x^0 \Mod$ with
		\[
		\calPGF_x \cap {{\calPGF_{x}}^\perp} = \calP_x, \quad \calGF_x \cap {{\calPGF_{x}}^\perp} = \calF_x, \quad {^\perp\calGI_x}\cap\calGI_x=\calI_x.
		\]
	\end{proposition}
	
	Now we are ready to obtain several abelian model structures on $\UC \Mod$.
	
	\begin{corollary}\label{PGF}
		Let $\UC$ be a generalized $k$-linear direct category such that it is projective both as a left and a right $\UC^0$-module. Then
		$( \Phi(\calPGF), \, \UC\Mod_{\calPGF^\perp}, \, \UC\Mod )$
		forms a hereditary Hovey triple in $\UC \Mod$, where ${\calPGF^\perp}= \{{\calPGF_{x}}^\perp\}_{x\in\Ob(\UC)}$.
	\end{corollary}
	
	\begin{prf*}
		We mention that $(\calPGF_x, {{\calPGF_{x}}^\perp}, A_x^0 \Mod)$ is a hereditary Hovey triple with $\calPGF_x\cap{{\calPGF_{x}}^\perp}=\calP_x$ by Proposition \ref{SS20}. The family $\{\calP_x \}_{x \in \Ob(\UC)}$ is cocompatible as $\UC^+$ is a projective left $\UC^0$-module; see Example \ref{exa-compatible}. The conclusion then follows from Corollary \ref{direct hHT}.
	\end{prf*}
	
	Similarly, we have the following result.
	
	\begin{corollary}\label{GF}
		Let $\UC$ be a generalized $k$-linear direct category such that it is projective both as a left and a right $\UC^0$-module. Then
		$( \Phi(\calGF), \, \UC\Mod_{\calPGF^\perp}, \, \UC\Mod_{\calCot} )$
		forms a hereditary Hovey triple in $\UC \Mod$.
	\end{corollary}
	
	The next result can be proved dually using Corollary \ref{inverse hHT}.
	
	\begin{corollary}\label{GI}
		Let $\UC$ be a generalized $k$-linear inverse category such that it is projective both as a left and a right $\UC^0$-module. Then
		$( \UC\Mod, \, \UC\Mod_{^{\perp}{\calGI}}, \, \Psi(\calGI) )$
		forms a hereditary Hovey triple in $\UC \Mod$, where $^\perp\calGI= \{^\perp{\calGI_{x}}\}_{x\in\Ob(\UC)}$.
	\end{corollary}

	We mention that the $k$-linearizations of quite a few combinatorial categories described before satisfy the conditions specified in these corollaries.
	
	\begin{example}\label{direct}
		Let $k$ be a field of characteristic 0. Then the $k$-linearizations of skeletons of the following combinatorial categories satisfy conditions specified in Corollaries \ref{PGF} and \ref{GF}:
		\begin{rqm}
			\item the category of finite sets and injections;
			\item the category of finite dimensional vector spaces over a finite field and injective linear maps;
			\item the category of finite totally ordered sets and order-preserving injections;
			\item the category of cyclically ordered sets and order-preserving injections.
		\end{rqm}
		
		Dually, via replacing injections in the definitions of the above categories by surjections, we obtain examples satisfying conditions specified in Corollary \ref{GI}.
	\end{example}
	
	
	
	\def\soft#1{\leavevmode\setbox0=\hbox{h}\dimen7=\ht0\advance \dimen7
		by-1ex\relax\if t#1\relax\rlap{\raise.6\dimen7
			\hbox{\kern.3ex\char'47}}#1\relax\else\if T#1\relax
		\rlap{\raise.5\dimen7\hbox{\kern1.3ex\char'47}}#1\relax \else\if
		d#1\relax\rlap{\raise.5\dimen7\hbox{\kern.9ex \char'47}}#1\relax\else\if
		D#1\relax\rlap{\raise.5\dimen7 \hbox{\kern1.4ex\char'47}}#1\relax\else\if
		l#1\relax \rlap{\raise.5\dimen7\hbox{\kern.4ex\char'47}}#1\relax \else\if
		L#1\relax\rlap{\raise.5\dimen7\hbox{\kern.7ex
				\char'47}}#1\relax\else\message{accent \string\soft \space #1 not
			defined!}#1\relax\fi\fi\fi\fi\fi\fi}
	\providecommand{\MR}[1]{\mbox{\href{http://www.ams.org/mathscinet-getitem?mr=#1}{#1}}}
	\renewcommand{\MR}[1]{\mbox{\href{http://www.ams.org/mathscinet-getitem?mr=#1}{#1}}}
	\providecommand{\arxiv}[2][AC]{\mbox{\href{http://arxiv.org/abs/#2}{\sf
				arXiv:#2 [math.#1]}}} \def\cprime{$'$}
	\providecommand{\bysame}{\leavevmode\hbox to3em{\hrulefill}\thinspace}
	\providecommand{\MR}{\relax\ifhmode\unskip\space\fi MR }
	\providecommand{\MRhref}[2]{%
		\href{http://www.ams.org/mathscinet-getitem?mr=#1}{#2}
	}
	\providecommand{\href}[2]{#2}

\end{document}